\documentclass[10pt]{article}
\usepackage[top=1in, bottom=1in, left=1in, right=1in]{geometry}
\usepackage{tikz,graphicx,amssymb,epstopdf,float,amsmath,amsthm,bm,color,epsfig,enumerate,enumitem,caption,hyperref,booktabs,multirow,array,lineno,comment}
\usetikzlibrary{shapes,arrows,matrix,patterns,positioning}
\usepackage[norelsize,boxed,linesnumbered,vlined,algo2e]{algorithm2e}

\newtheorem{theorem}{Theorem}[section]

\newtheorem{definition}[theorem]{Definition}
\newtheorem{lemma}[theorem]{Lemma}

\renewcommand{\appendix}[1]{
\section*{Appendix: #1}
}
\newcommand{\norm}[1]{\left\lVert#1\right\rVert}

\renewcommand{\O}[1]{O\left(#1\right)}
\renewcommand{\L}[1]{\log\left(#1\right)}

\newcommand{\bbC}{\mathbb{C}}
\newcommand{\bbR}{\mathbb{R}}

\renewcommand{\i}{i}
\renewcommand{\j}{j}
\newcommand{\calI}{\mathcal{I}}

\newcommand{\rID}{{\it rID}}
\newcommand{\cID}{{\it cID}}
\newcommand{\round}{\text{round}}
\newcommand{\length}{\text{length}}
\newcommand{\sz}{\text{size}}
\newcommand{\svd}{\text{svd}}
\newcommand{\qd}{\text{;\quad}}
\SetKwProg{Fn}{Function}{ }{ }
\SetKwFunction{rvone}{RecoveryVector1}
\SetKwFunction{rmone}{RecoveryMatrix1}
\SetKwFunction{pibone}{Partition1}
\SetKwFunction{rpm}{RecoveryPath}
\SetKwFunction{rvm}{RecoveryVector2}
\SetKwFunction{rmm}{RecoveryMatrix2}
\SetKwFunction{pibm}{Partition2}
\SetKwFunction{dt}{delaunayTriangulation}
\SetKwFunction{mst}{minspantree}
\SetKwFunction{bfs}{bfsearch}
\SetKwFunction{lrf}{LowRankFactorization}
\SetKwFunction{rsvd}{randomizedSVD}
\SetKwFunction{rp}{randperm}
\SetKwFunction{qr}{qr}
\SetKwFunction{lq}{lq}

\makeatletter
\newcommand*{\extendadd}{
  \mathbin{
    \mathpalette\extend@add{}
  }
}
\newcommand*{\extend@add}[2]{
  \ooalign{
    $\m@th#1\leftrightarrow$%
    \vphantom{$\m@th#1\updownarrow$}
    \cr
    \hfil$\m@th#1\updownarrow$\hfil
  }
}
\makeatother

\begin{document}

\title{Multidimensional Phase Recovery and\\ Interpolative Decomposition Butterfly Factorization}
\author{Ze Chen \\ Department of Mathematics, National University of Singapore, Singapore
    \and Juan Zhang \\ Department of Mathematics and Computational Science, Xiangtan University, China
    \and Kenneth L. Ho \\ Center for Computational Mathematics, Flatiron Institute, USA
    \and Haizhao Yang \\ Department of Mathematics, Purdue University, USA}

\maketitle

\begin{abstract}
    This paper focuses on the fast evaluation of the matrix-vector multiplication (matvec) $g=Kf$ for $K\in \mathbb{C}^{N\times N}$, which is the discretization of a multidimensional oscillatory integral transform $g(x) = \int K(x,\xi) f(\xi)d\xi$ with a kernel function $K(x,\xi)=e^{2\pi\i \Phi(x,\xi)}$, where $\Phi(x,\xi)$ is a piecewise smooth phase function with $x$ and $\xi$ in $\mathbb{R}^d$ for $d=2$ or $3$. A new framework is introduced to compute $Kf$ with $\O{N \L{N}}$ time and memories complexity in the case that only indirect access to the phase function $\Phi$ is available. This framework consists of two main steps: 1) an $\O{N \L{N}}$ algorithm for recovering the multidimensional phase function $\Phi$ from indirect access is proposed; 2) a multidimensional interpolative decomposition butterfly factorization (MIDBF) is designed to evaluate the matvec $Kf$ with an $\O{N \L{N}}$ complexity once $\Phi$ is available. Numerical results are provided to demonstrate the effectiveness of the proposed framework.
\end{abstract}


\section{Introduction} \label{sec:intro} 

This paper is concerned with the efficient evaluation of multidimensional oscillatory integral transforms. After discretization with $N$ grid points in each variable, the integral transform is reduced to a dense matrix-vector multiplication (matvec) as follows:
\begin{equation}\label{eq:kernel}
    g(x) = \sum_{\xi\in \Omega} K(x,\xi)f(\xi)=   \sum_{\xi\in \Omega}   e^{2\pi\i \Phi(x,\xi)} f(\xi),\quad x\in X,
\end{equation}
where $X$ and $\Omega$ are typically point sets in $\bbR^d$ for $d > 1$, $K(x,\xi) = e^{2\pi\i \Phi(x,\xi)}$ is a kernel function, $\Phi(x,\xi)$ is a piecewise smooth phase function with $\O{1}$ discontinuous points in $x$ and $\xi$, $f(\xi)$ is a given function, and $g(x)$ is a target function.

When the explicit formula of the kernel function is known, the direct computation of matvec in \eqref{eq:kernel} takes $\O{N^2}$ operations and is prohibitive in large-scale computation. There has been an active research line aiming at a nearly linear-scaling matvec for evaluating \eqref{eq:kernel}. In the case of uniformly distributed point sets $X$ and $\Omega$, the fast Fourier transform (FFT) \cite{doi:10.1137/1.9781611970999} can evaluate \eqref{eq:kernel} when $\Phi(x,\xi)=x\cdot\xi$ in $\O{N \L{N}}$ operations. When the point sets are non-uniform, the non-uniform FFT (NUFFT) algorithms in \cite{NUFFT,Alex2} are able to evaluate \eqref{eq:kernel} when $\Phi(x,\xi)=x\cdot\xi$ in $\O{N \L{N}}$ operations. For more general kernel functions, the butterfly factorization (BF) \cite{BF,Butterfly1,Butterfly2,IDBF} can factorize the dense matrix $e^{2\pi\i \Phi(x,\xi)}$ as a product of $\O{\L{N}}$ sparse matrices, each of which has only $\O{N}$ non-zero entries. Hence, storing and applying $e^{2\pi\i \Phi(x,\xi)}$ via the BF for evaluating \eqref{eq:kernel} take only $\O{N \L{N}}$ complexity. 

However, for multidimensional kernel functions, existing algorithms are efficient only when the explicit formula of the phase function $\Phi$ is known \cite{Gang,Alex2,Butterfly1,FIO09,IBF,BF,HSSBF,IDBF}. {The case of indirect access of the kernel function is illustrated in Table~\ref{tab:sc} for a list of different scenarios. When $\O{1}$ rows and columns of the phase matrices are available by solving PDE's, Scenario 3, as well as Scenario 1, are special cases of Scenario 2.  Therefore, we will focus more on Scenario 2 in this paper and will discuss the relationship between there Scenarios in detail. In fact, it is hard to evaluate any arbitrary entry of the kernel matrix directly in $\O{1}$ operations in Scenario 2. Therefore, the computational challenge in the case of indirect access of the kernel function motivates a series of new algorithms in this paper.}

\begin{table}[!ht]
    \begin{center}
        {\renewcommand\arraystretch{1.25}
        \begin{tabular}{|l|l|l|} \hline
            Scenario $1$ :& \multicolumn{2}{p{12cm}|}{\raggedright There exists an algorithm for evaluating an arbitrary entry of the kernel matrix $K$ in $\O{1}$ operations \cite{James:2017,Bremer201815,BF,Butterfly2}.} \\ \hline
            Scenario $2$ :& \multicolumn{2}{p{12cm}|}{\raggedright There exists an $\O{N \L{N}}$ algorithm for applying the kernel matrix $K$ and its transpose to a vector \cite{LUBF,BF,HSSBF,precon2}.} \\ \hline
            Scenario $3$ :& \multicolumn{2}{p{12cm}|}{\raggedright The phase functions $\Phi$ are solutions of partial differential equations (PDE's) \cite{Yingwave}. $\O{1}$ rows and columns of the phase matrices are available by solving PDE's.} \\ \hline
        \end{tabular}}
    \end{center}
    \caption{Three scenarios of the indirect access of the phase functions.}
    \label{tab:sc}
\end{table}

As the first main contribution of this paper, in the case of indirect access, a nearly linear scaling algorithm is proposed to recover multidimensional phase matrices in the form of low-rank matrix factorization. In scientific computing, several important problems require the construction of low-rank phase matrices \cite{James:2017,Bremer201815,precon3,precon1,precon2,Unw2,Unw3,Unw1,LUBF}. Previously, a nearly linear scaling algorithm has been proposed in \cite{NUFFTorBF} to recover the low-rank phase matrix with uniform discretization grid points in 1D. However, the 1D algorithm in \cite{NUFFTorBF} is problematic in the case of high-dimensional nonuniform discretization grid points. In this paper, we address the problem in multidimensional cases via Delaunay triangulation (DT) and minimum spanning tree (MST) construction. {Assuming the geometric coordinates of the discretization grids are given, and the indirect access of the phase functions is known, such as Scenario 2 in Table~\ref{tab:sc}. The phase matrices will be recovered to piecewise smoothness matrices by a fast MST algorithm based on DT. Then, low-rank approximations of the recovered phase matrices will be constructed.}

Secondly, when low-rank constructions of the phase matrices have been recovered, a new BF, multidimensional interpolative decomposition butterfly factorization (MIDBF), is proposed for the matvec $Kf$ with an $\O{N \L{N}}$ complexity for both precomputation and application. The MIDBF is a generalization of the interpolative decomposition butterfly factorization (IDBF) \cite{IDBF} in multidimensional cases especially when the discretization grid points are non-uniform. These two contributions lead to the first framework for multidimensional fast oscillatory integral transforms in the case of indirect access with non-uniform grid points.

The rest of the paper is organized as follows. In Section~\ref{sec:pre}, we revisit and generalize existing low-rank phase matrix factorization techniques, and propose a new low-rank matrix factorization in the case of indirect access. Next, the MIDBF will be introduced in Section~\ref{sec:MIDBF}. Finally, we provide several numerical examples to demonstrate the efficiency of the proposed framework in Section~\ref{sec:results}. For simplicity, we adopt MATLAB notations for the algorithm described in this paper: {given row and column index sets $I$ and $J$, $K(I,J)$ is the submatrix with entries from rows in $I$ and columns in $J$; the index set for an entire row or column is denoted as $``:"$.}

\section{Low-rank phase matrix factorization}
\label{sec:pre}

This section introduces a new low-rank phase matrix factorization for indirect access, which is the first main step in the proposed framework. We begin with a brief review of existing techniques and introduce a new algorithm afterward. These low-rank factorization methods will be applied repeatedly.

\subsection{Low-rank approximation by randomized sampling}
\label{sec:LRF}

Let us revisit an existing low-rank matrix factorization with linear complexity. For $A \in \bbC^{m\times n}$, a rank-$r$ approximate singular value decomposition (SVD) of $A$ is defined as
\begin{equation}\label{eqn:SVD}
    A \approx U \Sigma V^T,
\end{equation}
where $U\in \bbC^{m\times r}$ is orthogonal, $\Sigma\in \bbR^{r\times r}$ is diagonal, and $V\in \bbC^{n\times r}$ is orthogonal, and $r=\O{1}$ independent of the matrix size $m$ and $n$ with a prefactor depending only on the approximation error $\epsilon$. Previously, \cite{Butterfly4,randomSVD} have proposed efficient randomized tools to compute approximate SVDs for numerically low-rank matrices. The method in \cite{Butterfly4} is more attractive because it only requires $\O{1}$ randomly sampled rows and columns of $A$ for constructing \eqref{eqn:SVD} with $\O{m+n}$ operations and memories complexity, {and it is observed that $|A(i,j) - (U \Sigma V^T)(i,j)| = \O{\epsilon}$ in a probabilistic sense, where $1 \le i \le m$ and $1 \le j \le n$.}

The method is denoted as Function \rsvd and is presented in Algorithm~\ref{alg:rSVD}. Assuming the whole low-rank matrix $A$ is known, the input of Function \rsvd is $A$, $\O{1}$ randomly sampled row indices $\mathcal{R}$ and column indices $\mathcal{C}$, as well as a rank parameter {$r_\epsilon$ based on the error $\epsilon$}. Equivalently, it can also be assumed that $A(\mathcal{R},:)$ and $A(:,\mathcal{C})$ are known as the inputs. {Let $r$ be an empirical estimation of $r_\epsilon$,} then the outputs are three matrices $U\in \bbC^{m \times r}$, $\Sigma\in \bbR^{r \times r}$, and $V\in \bbC^{n\times r}$ satisfying \eqref{eqn:SVD}. {In Function \rsvd, for simplicity, given any matrix $K \in \bbC^{s\times t}$, Function \qr{K} performs a pivoted QR decomposition $K(:,P) = QR$, where $P$ is a permutation vector of the $t$ columns, $Q$ is a unitary matrix, and $R$ is an upper triangular matrix with positive diagonal entries in decreasing order.} Function \rp{m,r} denotes an algorithm that randomly selects $r$ different samples in the set $\{1,2,\dots,m\}$. If necessary, we can add an over sampling parameter $q$ such that we sample $rq$ rows and columns and only generate a rank $r$ truncated SVD in Line~\ref{alg:rSVD:svd} in Algorithm~\ref{alg:rSVD}. Larger $q$ results in better stability of Algorithm~\ref{alg:rSVD}.

\begin{algorithm2e}
    
    \Fn{$\left[U, \Sigma, V\right] \leftarrow$ \rsvd{$A, \mathcal{R}, \mathcal{C}, r$}}{
    
    $\left[m,n\right] \leftarrow \sz(A)$
    
    {$P \leftarrow $ \qr{$A(\mathcal{R},:)$} \qd $\Pi_{col} \leftarrow P(1:r)$
    \tcp*[f]{$A(\mathcal{R},P) = QR$}}
    
    {$P \leftarrow $ \qr{$A(:,\mathcal{C})^T$} \qd $\Pi_{row} \leftarrow P(1:r)$
    \tcp*[f]{$A(P,\mathcal{C}) = R^T Q^T$}}
    
    {$Q \leftarrow $ \qr{$A(:,\Pi_{col})$} \qd $Q_{col} \leftarrow Q(:,1:r)$
    \tcp*[f]{$A(P,\Pi_{col}) = QR$}}
    
    {$Q \leftarrow $ \qr{$A(\Pi_{row},:)^T$} \qd $Q_{row} \leftarrow Q(:,1:r)$
    \tcp*[f]{$A(\Pi_{row},P) = R^T Q^T$}}
    
    $S_{row} \leftarrow $ \rp{$m,r$} \qd $I \leftarrow [\Pi_{row}, S_{row}]$
    
    $S_{col} \leftarrow $ \rp{$n,r$} \qd $J \leftarrow [\Pi_{col}, S_{col}]$
    
    $M \leftarrow \left( Q_{col}(I,:) \right) ^\dagger A(I,J) \left( Q_{row}^T(:,J) \right) ^\dagger$
    \tcp*[f]{$(\cdot)^\dagger:$ pseudo-inverse}
    
    $\left[U_M, \Sigma_M, V_M \right]\leftarrow \svd(M)$ \label{alg:rSVD:svd}
    
    $U \leftarrow Q_{col} U_M$ \qd $\Sigma \leftarrow \Sigma_M$ \qd $V \leftarrow Q_{row} V_M$
    }
    
    \caption{Randomized sampling for a rank-$r$ approximate SVD with $\O{m+n}$ operations, such that $A \approx U \Sigma V^T$.}
    \label{alg:rSVD}
    
\end{algorithm2e}

\subsection{One-dimensional phase matrix factorization with indirect access}
\label{sec:1LR}

A nearly linear scaling algorithm for constructing the low-rank factorization of the phase matrix $\Phi\in\mathbb{R}^{N \times N}$ in \eqref{eq:kernel} has been proposed in \cite{NUFFTorBF} when only $\O{1}$ selected rows and columns of a 1D kernel matrix $K=e^{2\pi\i \Phi}$ with uniform discretization grid points are available as Scenario $2$ in Table~\ref{tab:sc}. In this subsection, we revisit the algorithms in \cite{NUFFTorBF} as a motivation for the multidimensional case proposed in this paper. The introduction of the 1D algorithms also helps to clarify the difficulties in the multidimensional case.

The difficulty of reconstructing $\Phi$ from $K=e^{2\pi\i\Phi}$ comes from the fact that
\[
    \frac{1}{2\pi} \Im\left(\L{K(i,j)})\right) =\frac{1}{2\pi}\Im\left(\L{e^{2\pi\i \Phi(i,j)}}\right) =  \frac{1}{2\pi}\arg\left( e^{2\pi\i \Phi(i,j)}\right)= \bmod(\Phi(i,j),1),
\]
where $\Im(\cdot)$ returns the imaginary part of the complex number, and $\arg(\cdot)$ returns the argument of a complex number. Thus, $\Phi$ is only known up to modular $1$. 

Since the point sets of the 1D kernel matrix are uniformly distributed, the main idea of \cite{NUFFTorBF} is to recover $\Phi$ by looking for the solution of the following combinatorial constrained $TV^3$-norm\footnote{{The $TV^3$-norm of a vector $v\in\mathbb{R}^N$ is defined as $\|v\|_{TV^3}:= \sum_{i=4}^{N} |v_{i}-3v_{i-1}+3v_{i-2}-v_{i-3}|$ in this paper.}} minimization problem: 

\begin{equation}
    \label{eqn:mintv}
        \begin{split}
        \smash{\displaystyle\min_{\Phi\in\mathbb{R}^{N\times N}}} &\quad \sum_{i\in\mathcal{R}}\|\Phi(i,:)\|_{TV^3}+ \sum_{j\in\mathcal{C}}\|\Phi(:,j)\|_{TV^3} \\
        \text{subject to} &\quad \bmod(\Phi(i,j),1) = \frac{1}{2\pi} \Im\left(\L{K(i,j)}\right) \text{ for }i\in \mathcal{R} \text{ or } j\in\mathcal{C},
    \end{split}
\end{equation}
where $\mathcal{R}$ and $\mathcal{C}$ are row and column index sets with $\O{1}$ randomly selected indices, respectively. The optimization problem above is appealing because it only requires the knowledge of $\O{1}$ rows and columns of $K$ and the computational cost in each iteration takes $\O{N}$ operations and memories. If the optimization problem can be solved in $\O{1}$ iterations, then the recovered rows and columns of $\Phi$ can be used to compute the low-rank factorization of $\Phi$ by Function \rsvd in Algorithm \eqref{alg:rSVD}. The final computational cost is nearly linear in $N$. However, due to the non-convexity of \eqref{eqn:mintv}, $\O{1}$ iterations are almost impossible to give a good solution unless a very good initial guess is available. This motivates \cite{NUFFTorBF} to design an empirical $\O{N}$ algorithm to provide a good initial guess to the optimization problem in \eqref{eqn:mintv}. 

The main algorithms of \cite{NUFFTorBF} are revisited and summarized in Algorithm~\ref{alg:rv1} and Algorithm~\ref{alg:rm1} in this paper for the preparation of higher dimensional cases. Algorithm~\ref{alg:rm1} relies on the repeated application of Algorithm~\ref{alg:rv1}, which adjusts the values of phase vectors by minimizing the absolute value of the third-order derivative, to provide an empirical solution to \eqref{eqn:mintv}. The functions in these two algorithms are denoted as \rvone and \rmone, respectively.  In fact, the algorithms presented in this paper are slightly different from those in \cite{NUFFTorBF} for robustness against discontinuity detection, which relies on a class of vectors $C_\tau$ with a threshold $\tau$ defined via: 
\begin{equation}\label{smooth}
    C_\tau = \left\{ u \in \mathbb{R}^n: |u(i)-3u(i-1)+3u(i-2)-u(i-3)|<\tau, \forall i \in \left\{ 4,5,\dots,n \right\} \right\}.
\end{equation}
Essentially, $C_\tau$ consists of vectors with a small absolute value of the third order derivative controlled by $\tau$ in the sense of finite difference. In our algorithms, if $|u(i)-3u(i-1)+3u(i-2)-u(i-3)| \ge \tau$, we will consider the original function that generates $u$ to be discontinuous at the location corresponding to $u(i)$. With this definition ready, we are able to explain our algorithms as follows.

For Function \rvone in Algorithm~\ref{alg:rv1}, input variables are a vector $u$ of length $N$, a discontinuity detection parameter $\tau$, and a parameter $flag$ which indicates whether $u$ will be recovered from the first entry or the fourth entry. Then, the outputs are a smooth vector $v$ satisfying $\bmod(v,1)=\bmod(u,1)$ and a vector of indices $\mathcal{D}$ for discontinuity locations.

\begin{algorithm2e}
    
    \Fn{$\left[v,\mathcal{D}\right] = $ \rvone{$u,\tau,flag$}}{
    
    $N \leftarrow \length(u) \qd v \leftarrow u \qd \mathcal{D} \leftarrow [1] \qd n \leftarrow 1 \qd c \leftarrow 1$
    
    \While{$c\leq n$}{
        
        $st \leftarrow \mathcal{D}(c)$
        
        \If{$flag \sim= 1$ \rm{\textbf{or}} $st \sim=1 $}{
        
            $v(st+1) \leftarrow u(st+1) - \round(u(st+1) - v(st))$
            
            $v(st+2) \leftarrow u(st+2) - \round(u(st+2) - 2v(st+1) + v(st))$}
        
        \For{$a = st+3 : N$}{
        
            $v(a) \leftarrow u(a) - \round(u(a) - 3v(a-1) + 3v(a-2) - v(a-3))$
            
            \If{$|v(a)-3v(a-1)+3v(a-2)-v(a-3)| \ge \tau$ \rm{\textbf{and}} $a \le N-3$}{
            
                $\mathcal{D} \leftarrow \left[\mathcal{D}, a \right] \qd n \leftarrow n+1$
                \tcp*[f]{detect discontinuous locations}
                
                $v(a) \leftarrow u(a) - \round(u(a) - v(a-1))$
            
                \textbf{Break}
            
            }
            }
            
            $c \leftarrow c + 1$
        
        }
    }
    
    \caption{An $\O{N}$ algorithm for recovering a vector $v$ from the observation $u=\bmod(v,1)$. The locations of discontinuity in $v$ are automatically detected. A vector $v$ is identified via empirically minimizing the magnitude of the absolute value of its third-order derivative.}
    \label{alg:rv1}
    
\end{algorithm2e}

In Function \rmone in Algorithm~\ref{alg:rm1}, one of the input variables is a function handle $\Phi$, which can evaluate an arbitrary row or column of the phase matrix. The other inputs are a vector $\mathcal{R}$ and a vector $\mathcal{C}$ as the row and column index sets indicating $\O{1}$ randomly selected rows and columns of the phase matrix, as well as a discontinuity detection parameter $\tau$. 

\begin{figure}[!ht]
    \centering
    \begin{tabular}{cccc}
        \includegraphics[height=0.6in]{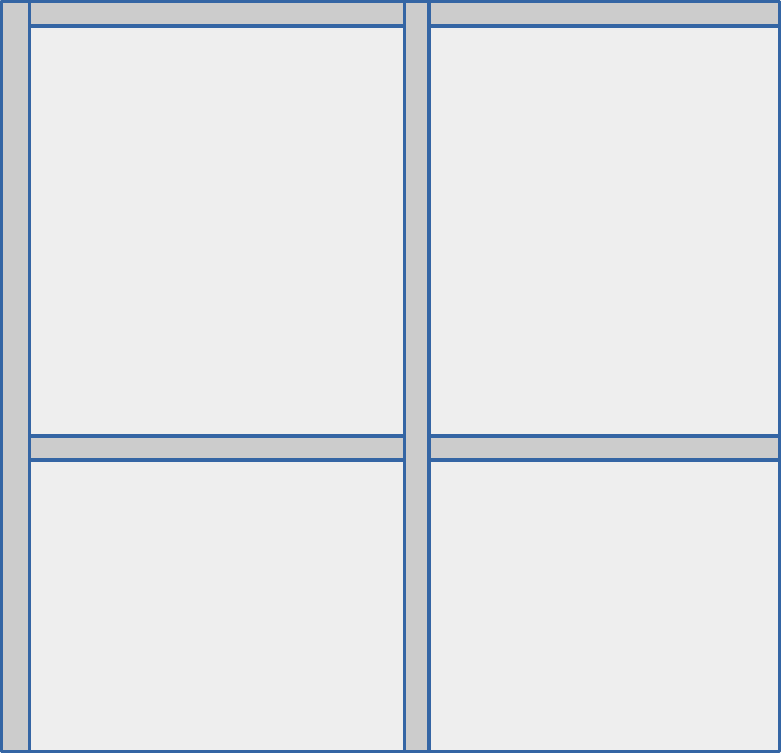} & \includegraphics[height=0.6in]{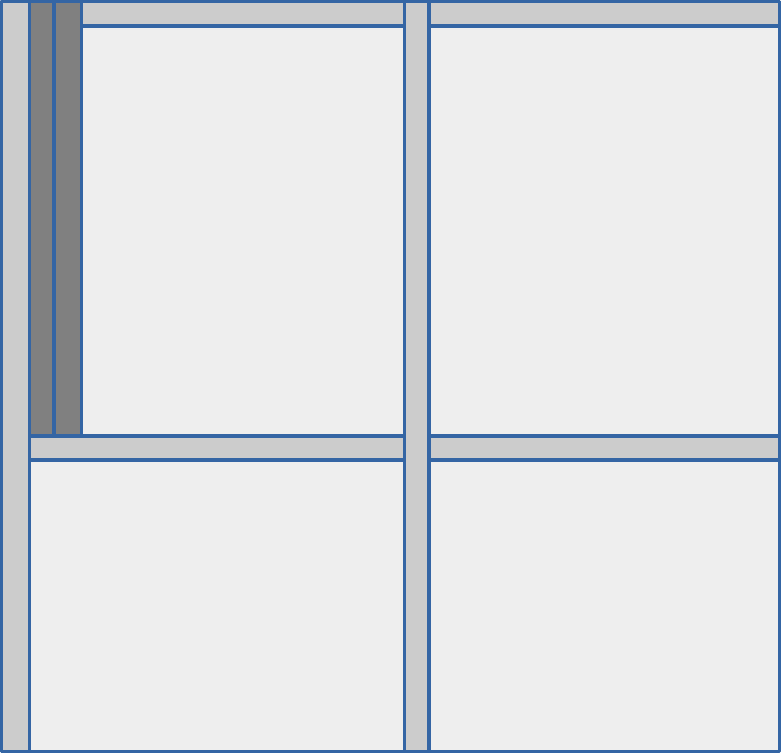} & \includegraphics[height=0.6in]{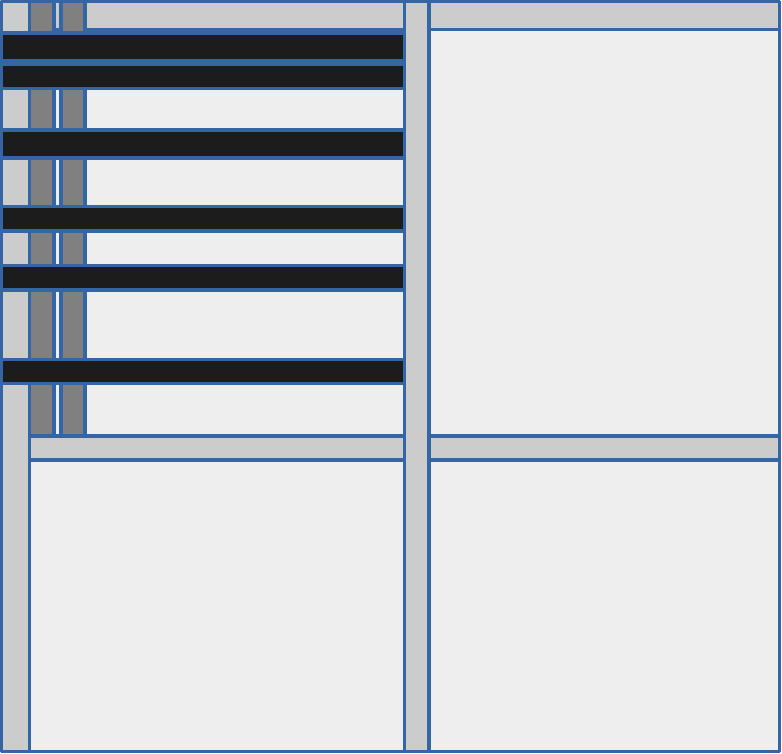} & \includegraphics[height=0.6in]{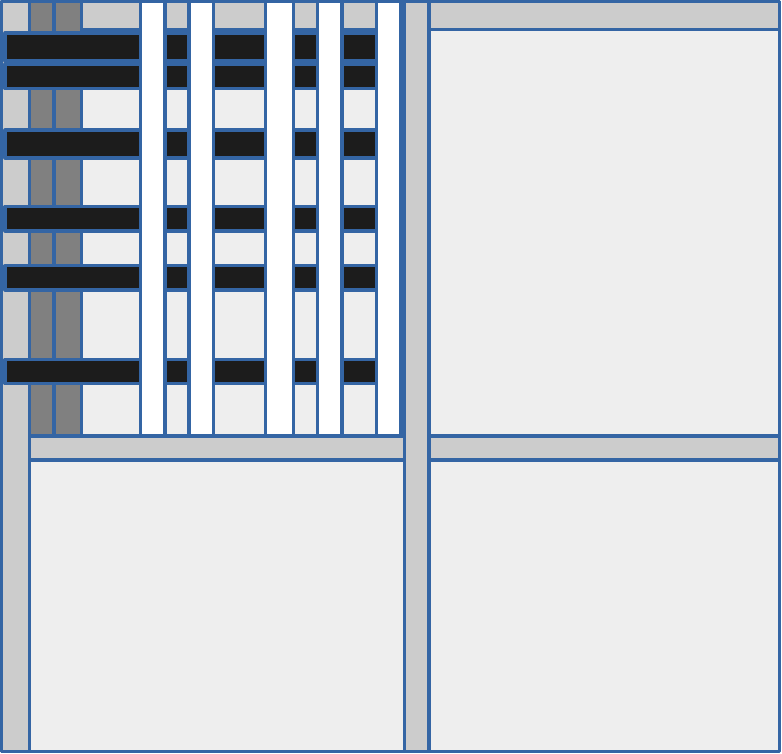}\\
        (a) & (b) & (c) & (d)
    \end{tabular}
    \caption{An illustration of the low-rank matrix recovery for a 1D phase matrix in Algorithm~\ref{alg:rm1}. (a) Line~\ref{alg:rm1:pib} partitions the phase matrix into submatrices such that there is no discontinuity along rows and columns in each submatrix. Then, Line~\ref{alg:rm1:rv11}-\ref{alg:rm1:rv12} recovers the first row and column of each submatrix. (b) Next, Line~\ref{alg:rm1:rv12} recovers the second and the third columns for each submatrix. (c) Next, Line~\ref{alg:rm1:rv13}-\ref{alg:rm1:rv14} recovers $\O{1}$ rows (including the second row and the third row) of each submatrix. (d) Finally, Line~\ref{alg:rm1:rv15} recovers $\O{1}$ columns of each submatrix.}
    \label{fig:trace1}
\end{figure}

Because it is more convenient to apply Algorithm~\ref{alg:rv1} to recover a vector representing a continuous function, we first apply Algorithm~\ref{alg:rv1} with $\tau$ to identify the sets of discontinuous points $\mathcal{D}_r$ and $\mathcal{D}_c$, each of which contains the first index $1$. Next, the phase matrix is partitioned into $n_r \times n_c$ blocks, each of which is denoted as $\Phi.\mathcal{B}_s \mathcal{B}_t$ representing a continuous piece of the phase function, where $n_r$ is the cardinality of $\mathcal{D}_r$, $n_c$ is the cardinality of $\mathcal{D}_c$, $s = 1, 2, \dots, n_r$, and $t = 1, 2, \dots, n_c$. This procedure is referred to as the Function \pibone in Line~\ref{alg:rm1:pib} in Algorithm~\ref{alg:rm1}. Similarly, $\mathcal{R}$ and $\mathcal{C}$ are partitioned into $n_r$ and $n_c$ parts by $\mathcal{D}_r$ and $\mathcal{D}_c$, and saved as $\mathcal{R}.\mathcal{B}_s$ and $\mathcal{C}.\mathcal{B}_t$ respectively. For example, Panel (a) in Figure~\ref{fig:trace1} visualizes an example when the phase function contains only 4 continuous blocks: $\Phi.\mathcal{B}_1 \mathcal{B}_1$, $\Phi.\mathcal{B}_1 \mathcal{B}_2$, $\Phi.\mathcal{B}_2 \mathcal{B}_1$, $\Phi.\mathcal{B}_2 \mathcal{B}_2$. Panel (c) and (d) in Figure~\ref{fig:trace1} visualize the randomly selected rows $\mathcal{R}.\mathcal{B}_1$ and columns $\mathcal{C}.\mathcal{B}_1$ in $\Phi.\mathcal{B}_1 \mathcal{B}_1$.

Finally, the selected rows and columns are recovered by Algorithm~\ref{alg:rv1} with a carefully designed order in Line~\ref{alg:rm1:rv11}-\ref{alg:rm1:rv15} in Algorithm~\ref{alg:rm1}. The parameter for detecting discontinuous points is set to $1$ since there is no need to detect discontinuity anymore. Note that there is no uniqueness for recovering a smooth vector from its values after $\bmod$ $1$. Hence, we introduce the specially designed order in Line~\ref{alg:rm1:rv11}-\ref{alg:rm1:rv15} to guarantee that each recovered row and column at their intersection share the same value, as long as the discontinuous points in the phase function are well distinguished by a parameter $\tau$ from continuous points, which can be shown by Lemma~\ref{1LRthm} below.

\begin{algorithm2e}
    
    \Fn{$\left[ \Phi,\mathcal{R},\mathcal{C} \right] = $ \rmone{$\Phi,\mathcal{R},\mathcal{C},\tau$}}{
    
    $\mathcal{D}_r \leftarrow$ \rvone{$\Phi(:,\mathcal{C}(1)),\tau,0$}
    \tcp*[f]{$\mathcal{D}_r:$ discontinuous point set}
    
    $\mathcal{D}_c \leftarrow$ \rvone{$\Phi(\mathcal{R}(1),:),\tau,0$}
    \tcp*[f]{$\mathcal{D}_c:$ discontinuous point set}
    
    $\mathcal{R} \leftarrow \left[\mathcal{R}, \mathcal{D}_r \right]$ \qd $\mathcal{C} \leftarrow \left[\mathcal{C}, \mathcal{D}_c \right]$
    
    $n_r \leftarrow $ \length($\mathcal{D}_r$) \qd $n_c \leftarrow $ \length($\mathcal{D}_c$)
    
    $\left[\Phi, \mathcal{R}, \mathcal{C}\right] \leftarrow $ \pibone{$\Phi, \mathcal{R}, \mathcal{C}, \mathcal{D}_r, \mathcal{D}_c$} \label{alg:rm1:pib}
    
    \For{$s = 1:n_r$}{
    
        \For{$t = 1:n_c$}{
            $\Phi.\mathcal{B}_s\mathcal{B}_t(1,:) \leftarrow$ \rvone{$\Phi.\mathcal{B}_s\mathcal{B}_t(1,:), 1, 0$} \label{alg:rm1:rv11}
            
            $\Phi.\mathcal{B}_s\mathcal{B}_t(:,k) \leftarrow$ \rvone{$\Phi.\mathcal{B}_s\mathcal{B}_t(:,k), 1, 0$} for $k = 1,2,3$ \label{alg:rm1:rv12}
            
            $\Phi.\mathcal{B}_s\mathcal{B}_t(k,:) \leftarrow$ \rvone{$\Phi.\mathcal{B}_s\mathcal{B}_t(k,:), 1, 1$} for $k = 2,3$ \label{alg:rm1:rv13}
            
            $\Phi.\mathcal{B}_s\mathcal{B}_t(\mathcal{R}.\mathcal{B}_s(k),:) \leftarrow$ \rvone{$\Phi.\mathcal{B}_s\mathcal{B}_t(\mathcal{R}.\mathcal{B}_s(k),:), 1, 1$} for all $k$ \label{alg:rm1:rv14}
            
            $\Phi.\mathcal{B}_s\mathcal{B}_t(:,\mathcal{C}.\mathcal{B}_t(k)) \leftarrow$ \rvone{$\Phi.\mathcal{B}_s\mathcal{B}_t(:,\mathcal{C}.\mathcal{B}_t(k)), 1, 1$} for all $k$ \label{alg:rm1:rv15}
        
        }
    }
    }
    
    \caption{An $\O{N}$ algorithm for the approximate solution of the $TV^3$-norm minimization when the phase function $\Phi(x,\xi)$ is defined on $\mathbb{R}\times \mathbb{R}$.} 
    \label{alg:rm1}
    
\end{algorithm2e}

\begin{lemma}
    \label{1LRthm}
    {Given $\bmod(\phi,1)\in \mathbb{R}^{n \times m}$ and the recovered values of $\phi(1:3,1:3)$, where $\phi$ is a one-dimensional phase matrix.}
    Assuming that all rows and columns of $\phi$ belong to the class $C_\tau$ with a threshold $\tau \le \frac{1}{16}$, then the intersection of each recovered row and column by Algorithm~\ref{alg:rm1} share the same value.
\end{lemma}

The proof of Lemma~\ref{1LRthm} can be found in the appendix. The correct $\tau$ depends on the phase function and is not known a priori. In practice, $\tau$ is set as $\frac{1}{16}$ according to Lemma~\ref{1LRthm} and it performs good enough to identify $\O{1}$ discontinuous points with $\O{N}$ operations.

Once the phase function recovery algorithm in Algorithm~\ref{alg:rm1} is ready, following the idea of low-rank matrix factorization via randomized sampling in Algorithm~\ref{alg:rSVD}, we can obtain a nearly linear scaling algorithm to construct the low-rank factorization of the phase matrix.

\subsection{Multidimensional phase matrix factorization with indirect access}
\label{sec:nLR}

\subsubsection{Overview}

In this subsection, a nearly linear scaling algorithm for constructing the low-rank factorization of the multidimensional phase matrix $\Phi\in\mathbb{R}^{N \times N}$ will be introduced when we only know the kernel matrix $K=e^{2\pi\i \Phi}$ with non-uniform discretization grid points through Scenario $2$ in Table~\ref{tab:sc}. In the multidimensional case, the coordinates of $N \times N$ discretization grid points will be required for our methods, where $N = n^d$ is the number of points in a $d$-dimensional domain, $d = 2$ or $3$, and $n$ is the number of points in each dimension. Recall that the main purpose of our algorithm is to recover $\O{1}$ randomly selected rows and columns of $\Phi$, and construct the low-rank factorization in the end.

In Scenario $2$, applying the kernel matrix $K$ and its transpose to $\O{1}$ randomly selected natural basis vectors in $\mathbb{R}^{N}$ can obtain the rows and columns of $K$ in $\O{N \L{N}}$ operations. Notice that Scenario $1$ is a special case of Scenario $2$, we only focus on Scenario $2$ for phase recovery.

Similar to the 1D case, instead of recovering the exact $\Phi$ that generates $K$, our primary purpose is to find a low-rank matrix $\Psi$ such that 
\begin{equation}\label{eqn:tr}
    \bmod(\Psi,1) = \frac{1}{2\pi} \Im\left(\L{K}\right).
\end{equation} 

Based on the piecewise smoothness of the multidimensional phase function, a recovery algorithm similar to the 1D case can be proposed to recover the rows and columns of $\Phi$ up to an additive error matrix $E$ that is numerically low-rank, i.e., the method returns a matrix $\Psi=\Phi+E$ such that $e^{2\pi\i \Psi}=e^{2\pi\i \Phi}$ and $E$ is numerically low-rank. However, the discretization of the integral operator especially in the case of non-uniform grid points can introduce ``artificial'' discontinuity along the rows and columns of the phase matrix. Hence, it is impossible to apply the vector class $C_\tau$ and the algorithms in the 1D case. Although informally the recovery problem can be stated as
\begin{equation}
    \label{eqn:smphi}
    \begin{split}
        \text{Find} &\quad \text{ piecewise smooth } \Psi(i,:) \text{ and } \Psi(:,j) \\
        \text{subject to} &\quad  \bmod(\Psi(i,j),1) = \frac{1}{2\pi} \Im\left(\L{K(i,j)}\right) \text{ for }i\in \mathcal{R} \text{ or } j\in\mathcal{C}.
    \end{split}
\end{equation}

Notice that the vectors $\Psi(i,:)$ and $\Psi(:,j)$ are not ``smooth'' at the location when adjacent entries are corresponding to non-adjacent points in the high-dimensional spatial domain in $\mathbb{R}^d$. In other words, the definition of the smoothness of these vectors should rely on the smoothness of the phase function in the original domain in $\mathbb{R}^d$ instead of the difference of adjacent entries as in \eqref{smooth}. 

How to recover such piecewise smooth vectors is the main difficulty of the extension of the 1D algorithm to high-dimensional cases. A naive algorithm is to identify the value according to the adjacent point with the smallest distance through all points. However, this takes $\O{N^2}$ operations to find the adjacent point. In other words, how to solve this difficulty with nearly linear computational complexity is the main challenge for us.

\subsubsection{Vector recovery}

Let us use the example of a vector recovery in the high-dimensional case to illustrate the ideas to conquer the difficulty mentioned above. Suppose $v$ is the discretization of a piecewise smooth function $\phi(x)$ with $N$ (possibly nonuniform) grid points in $[0,1]^d$ and $\O{1}$ pieces of domains in which $\phi(x)$ is smooth. The spatial locations of the $N$ grid points are stored in a matrix $\mathcal{X} \in \mathbb{R}^{N \times d}$, i.e., $\mathcal{X}(i,:)$ is the location of the $i$-th entry of $v$. Assume that $k$ is a vector representing $e^{2\pi i \phi(x)}$ using the same discretization.  Informally, the vector recovery problem is to find a ``piecewise smooth'' vector  $v$ subject to $\bmod(v,1) = \frac{1}{2\pi} \Im\left(\L{k}\right)$.

To conquer the difficulty of artificial discontinuity, the entry values of $v$ are identified via minimizing the variation of $\phi(x)$ using physically adjacent locations in $\mathbb{R}^d$. For this purpose, we introduce a special recovery path matrix $P\in\mathbb{Z}^{(N-1)\times 2}$ with a beginning Node $q$ such that $P(:,2)$ is a permutation of $\{1,2,\dots,N\}\setminus q$, and $(P(i,1),P(i,2))$ is a pair of indices of $v$ with corresponding spatial locations adjacent to each other in $\mathbb{R}^d$, i.e., $\mathcal{X}(P(i,1),:)$ is an adjacent grid point of $\mathcal{X}(P(i,2),:)$ in $\mathbb{R}^d$. 

If the recovery path matrix $P$ and a set of indices for discontinuous locations $\mathcal{D}$ are given, the recovery of $v$ can be solved via the optimization problem:
\begin{equation}
    \label{eqn:mv2}
    \begin{split}
        \smash{\displaystyle\min_{ v\in\mathbb{R}^{N}}} &\quad \sum_{i\in \{1,\dots,N-1\}\setminus \mathcal{D}}|v(P(i,2))-v(P(i,1))|\\
        \text{subject to} &\quad \bmod(v,1) = \frac{1}{2\pi} \Im\left(\L{k}\right).
    \end{split}
\end{equation}

We will introduce the construction of $P$ later and focus on the construction of $\mathcal{D}$ and a nearly linear scaling empirical solution to \eqref{eqn:mv2} first. Similarly to the 1D case, to detect discontinuity of the piecewise smooth function automatically, we define a class of vectors $C_{\tau,P}$ for a threshold $\tau$ and a recovery path matrix $P$ via:
\[
    C_{\tau,P} = \left\{ v \in \mathbb{R}^n: |v(P(i,2))-v(P(i,1))|<\tau, \forall i \in \left\{ 1,2,\dots,n-1 \right\} \right\}.
\]
$C_{\tau,P}$ consists of vectors with a small absolute value of the first order derivative controlled by $\tau$ in the sense of finite difference. In our assumption, if $|v(P(i,2))-v(P(i,1))| \ge \tau$, we will consider the original function that generates $v$ to be discontinuous at the location $\mathcal{X}(P(i,2),:)$, which will be justified by our method afterwards.

Function \rvm in Algorithm~\ref{alg:rv2} below identifies a piecewise smooth vector $v$ from a given vector $u = \frac{1}{2\pi} \Im\left(\L{k}\right)$ via empirically minimizing $|v(P(i,2)) - v(P(i,1))|$ such that $\bmod(v(P(i,2)),1)=u(P(i,2))$, for each $i=1,2,\dots,N$ (corresponding to Line~\ref{alg:rv2:rec} in Algorithm~\ref{alg:rv2}). Each smooth piece of $v$ belongs to $C_{\tau,P}$. The discontinuity location $i$ will be detected and assigned to the discontinuity location set $\mathcal{D}$ if $|v(P(i,2)) - v(P(i,1))| \ge \tau$. It is clear that the complexity of Algorithm~\ref{alg:rv2} to empirically solve \eqref{eqn:mv2} and detect discontinuity is $\O{N}$. Note that Function \rvm in Algorithm~\ref{alg:rv2} is based on the first-order derivative of the phase function while Function \rvone in Algorithm~\ref{alg:rv1} is based on the third-order derivative. It is a simple extension to apply higher order derivative in Algorithm~\ref{alg:rv2} using the high-order finite difference schemes in \cite{JIANCHUN1995162,VASILYEV2000746}, which is left as future work if necessary.

\begin{algorithm2e}
    
    \Fn{$\left[v,\mathcal{D}\right] = $ \rvm{$u, \tau, P$}}{
    
    $N \leftarrow \length(u) \qd \mathcal{D} \leftarrow [1] \qd v \leftarrow u$
    
    \For{$c = 1:N-1$}{
        
        $bg \leftarrow P(c, 1) \qd ed \leftarrow P(c, 2)$

        $v(ed) \leftarrow u(ed) - \round(u(ed) - v(bg))$ \label{alg:rv2:rec}
        
        \If{$|v(ed)-v(bg)| \ge \tau$}{
            $\mathcal{D} \leftarrow \left[\mathcal{D}, ed \right]$
            \tcp*[f]{detect discontinuous locations}
        }
    }
    }
    
    \caption{An $\O{N}$ algorithm for recovering a vector $v$ from the observation $u = \bmod(v,1)$ and detecting discontinuity using the recovery path matrix $P$.} 
    \label{alg:rv2}
    
\end{algorithm2e}

\subsubsection{Recovery path}

The main challenge of vector recovery is to identify a recovery path matrix $P$ efficiently. Recall that the naive algorithm to identify an adjacent point of a given location is to traverse all other points, compute distances, and pick up the smallest one, which needs $\O{N^2}$ operations to construct $P$ for $N$ points.

\begin{figure}[!ht]
    \centering
    \begin{tabular}{ccc}
        \fbox{\includegraphics[height=3cm,trim={8cm 8cm 6cm 6cm},clip]{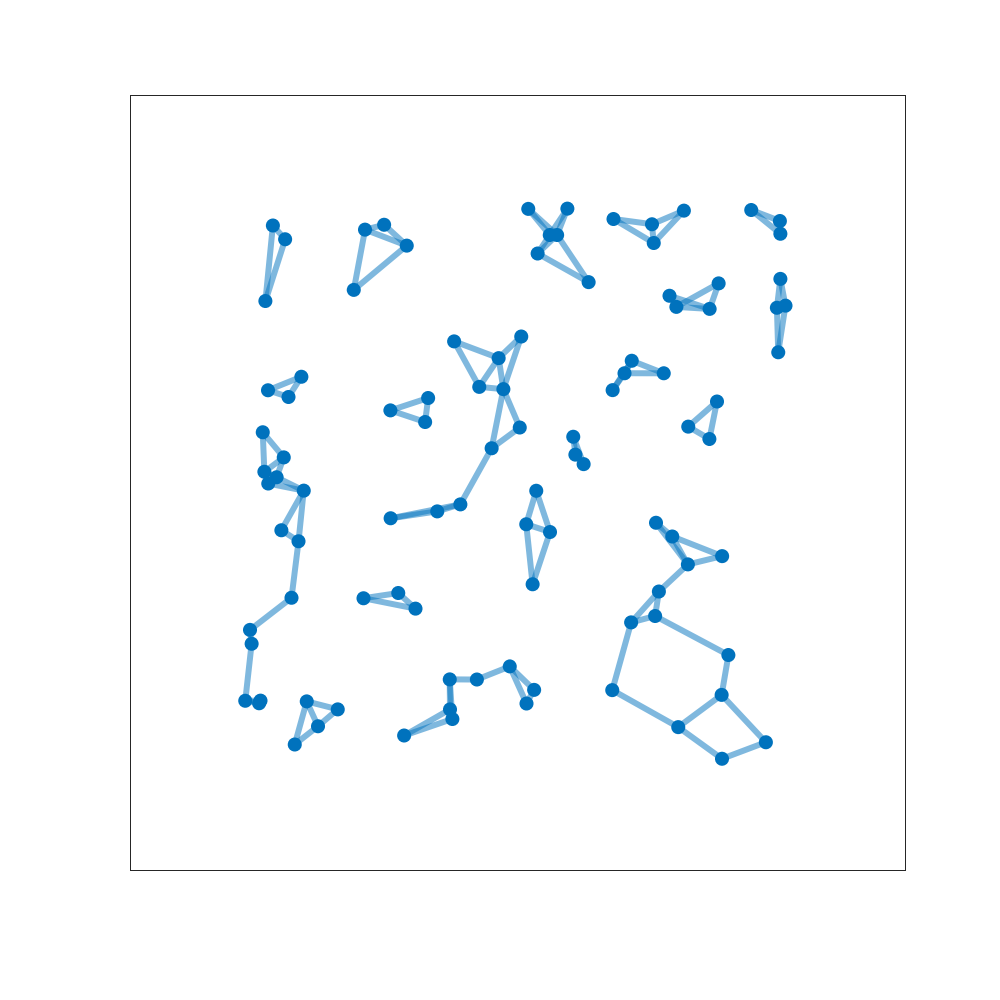}} & &        \fbox{\includegraphics[height=3cm,trim={8cm 8cm 6cm 6cm},clip]{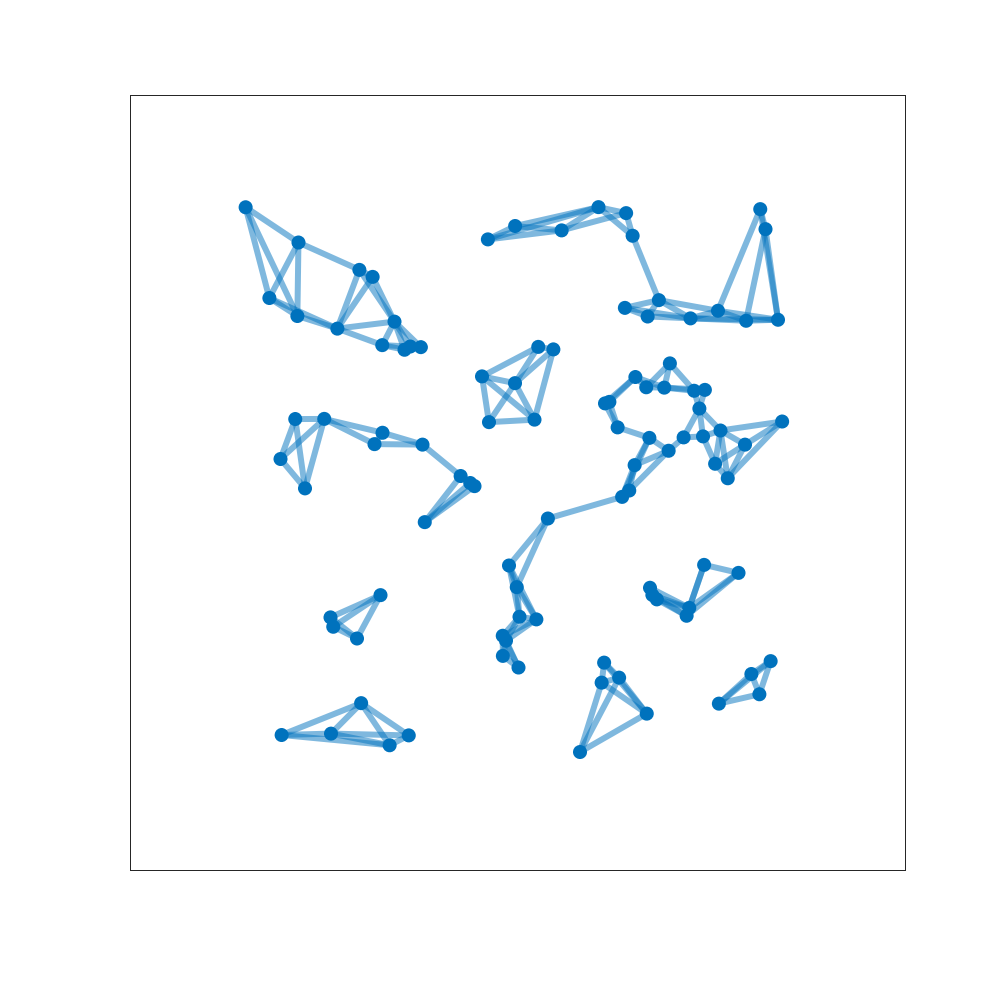}} \\
        (a) & & (b)
    \end{tabular}
    \caption{(a) 100 randomly generated points connected with their 2-nearest neighbors. (b) 100 randomly generated points connected with their 3-nearest neighbors.}
    \label{fig:KNN}
\end{figure}

First of all, we consider an algorithm for constructing a recovery path matrix based on $k$-nearest neighbors algorithm in $\O{N \L{N}}$ operations \cite{KNN}. When $k$-nearest neighbors of each point are found, the recovery path can be constructed by the edges between each point and its $k$-nearest neighbors. However, it is not efficient to find an integrated recovery path through all points. For example, in Figure~\ref{fig:KNN} (a), for a row vector $v$ of a phase matrix, 100 points as the locations of $v$ are randomly generated and connected with their $2$-nearest neighbors. Then, the result shows that this graph is split to $19$ connected components. If we recover $v$ for each component, at least $19$ column indices of the phase matrix should be selected as initialization. Another similar example of a graph for connecting $3$-nearest neighbors is illustrated in Figure~\ref{fig:KNN} (b). In addition, for a graph of $N$ points connected with their $k$-nearest neighbors, the largest number of connected components is $\O{\frac{N}{k+1}}$. Thus, this method may not be robust compared to our assumption: only $\O{1}$ rows and columns of the kernel matrix can be used for recovery.

Secondly, we also consider an algorithm based on a radius search in $\O{N \L{N}}$ operations \cite{RadiusSearch}. By this method, a graph of recovery path can be generated by connecting each point with their neighbors no further apart than a search radius. Unfortunately, this graph may also be split to a number of connected components, which depends on the selection of the search radius. Otherwise, how to choose a search radius and detect discontinuity will become new challenges.

Therefore, in the rest of this subsection, we propose an algorithm based on the Delaunay triangulation (DT) and the minimum spanning tree (MST) with nearly linear computational complexity instead of $k$-nearest neighbors and radius search algorithm to conquer the main difficulty of vector recovery.

\begin{definition}\label{def:dt} 
For a set of points in the $d$-dimensional Euclidean space with locations $\mathcal{X}\in\mathbb{R}^{N\times d}$, a {\bf Delaunay triangulation} is a triangulation DT($\mathcal{X}$) such that no point in this set is inside the circum-hypersphere of any $d$-simplex in DT($\mathcal{X}$). 
\end{definition}

\begin{definition}\label{def:mst}
A {\bf minimum spanning tree} (MST) $\mathcal{T}$ is a subset of the edges of a connected, edge-weighted undirected graph $\mathcal{G}$ that connects all the vertices, without any cycle and with the minimum possible total edge weight.
\end{definition}

DTs are widely used in scientific computing in many diverse applications. The Delaunay criterion is the fundamental property of DTs, which is often called as the empty circumcircle criterion in the case of 2D triangulations. In other words, a Delaunay triangulation of a set of points in 2D ensures the circumcircle associated with each triangle containing no other point in its interior. This property can be extended to higher dimensions. For instance, in 3D cases, the triangulation of a set of points is composed of tetrahedra. Then, the circumspheres of all tetrahedra also satisfy the empty circumsphere criterion.

In our problem, given the location matrix $\mathcal{X}\in\mathbb{R}^{N\times d}$ of $N$ points in $\mathbb{R}^d$, DT($\mathcal{X}$) can be treated as a fully connected undirected graph $\mathcal{G}$ with edges weighted by the Euclidean distance of two connected points. Due to the property of DT, useless long edges between $\mathcal{X}$ can be eliminated efficiently. Since a DT is a planar graph, and there are no more than three times as many edges as vertices in any planar graph, DT($\mathcal{X}$) will generate only $\O{N}$ edges. Moreover, it has been a standard routine to identify DT($\mathcal{X}$) with an expected runtime bounded by $\O{N \L{N}}$ for $d=2$ or $3$ (e.g., see \cite{dt1,dt2,dt3}). 

Based on the fact in \cite{DTMST} that the set of edges of DT($\mathcal{X}$) contains an MST for $\mathcal{X}$, we can use an MST $\mathcal{T}(\mathcal{X})$ as an efficient representation of the graph $\mathcal{G} = \text{DT} (\mathcal{X})$. Since there are $\O{N}$ edges in DT($\mathcal{X}$), any of the standard minimum spanning tree algorithms is able to find $\mathcal{T}(\mathcal{X})$ with an $\O{N \L{N}}$ complexity such as the Prim's algorithm \cite{Prim}.

Finally, a recovery path matrix $P$ can be identified following the order of nodes in $\mathcal{T}(\mathcal{X})$. Breadth-first search algorithm \cite{BFS} can be applied for traversing $\mathcal{T}(\mathcal{X})$ starting from the root $q$ and exploring all of the neighbor nodes at the present depth prior to moving on to the nodes at the next depth level. It is an efficient method for constructing $P$ with an $\O{N}$ complexity. Otherwise, the definition of the recovery path matrix $P$ is modified according to $\mathcal{T}$ as follows. 

\begin{definition}\label{def:rpm}
Given an MST $\mathcal{T}$ with $N$ nodes and the root at Node $q$, a {\bf recovery path matrix} $P\in\mathbb{Z}^{(N-1)\times 2}$ associated to $\mathcal{T}$ is a matrix such that 1) $P(:,2)$ is a permutation vector of $\{1,2,\dots,N\}\setminus q$; 2) the depth of Node $P(i,2)$ is less than or equal to that of Node $P(j,2)$ if $i\leq j$; 3) Node $P(i,1)$ is the predecessor node of Node $P(i,2)$ in $\mathcal{T}$ for all $i = 1, 2, \dots, N-1$.
\end{definition}

\begin{figure}[!ht]
    \centering
    \begin{tabular}{ccccccc}
        \includegraphics[height=0.6in,trim={6.5cm 9.5cm 5.5cm 8.5cm},clip] {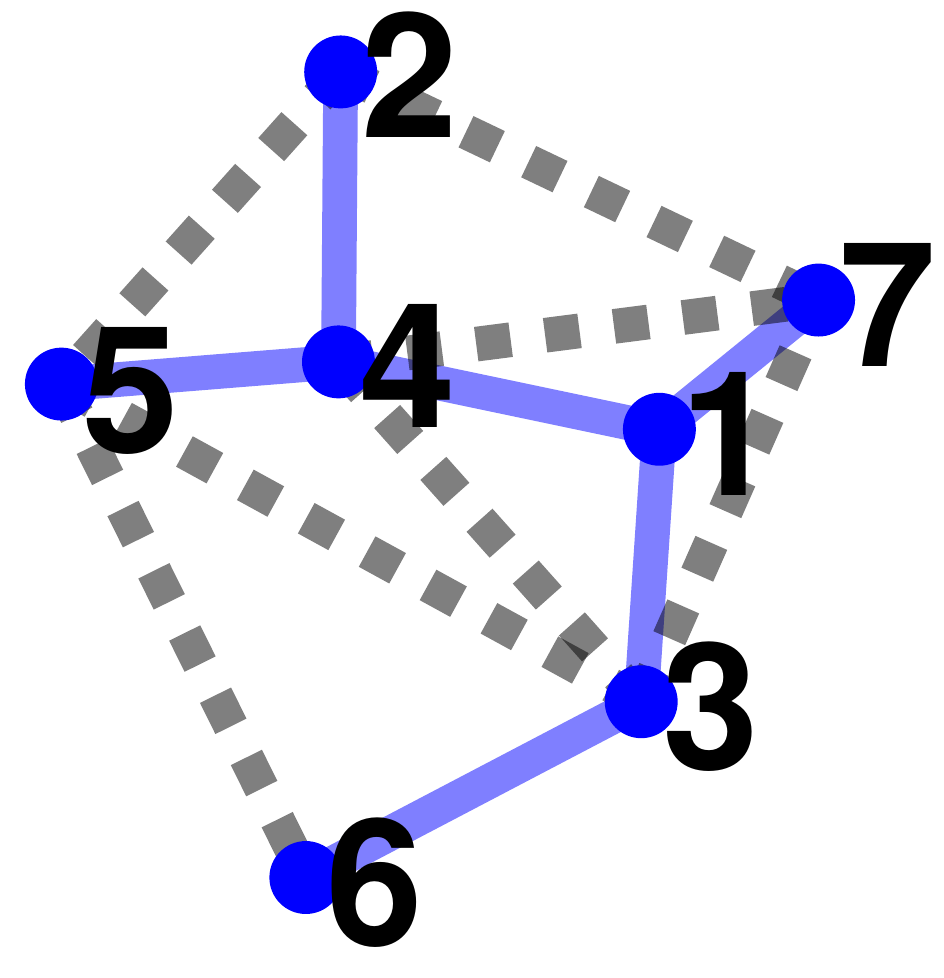} &  \includegraphics[height=0.6in,trim={6.5cm 9.5cm 5.5cm 8.5cm},clip] {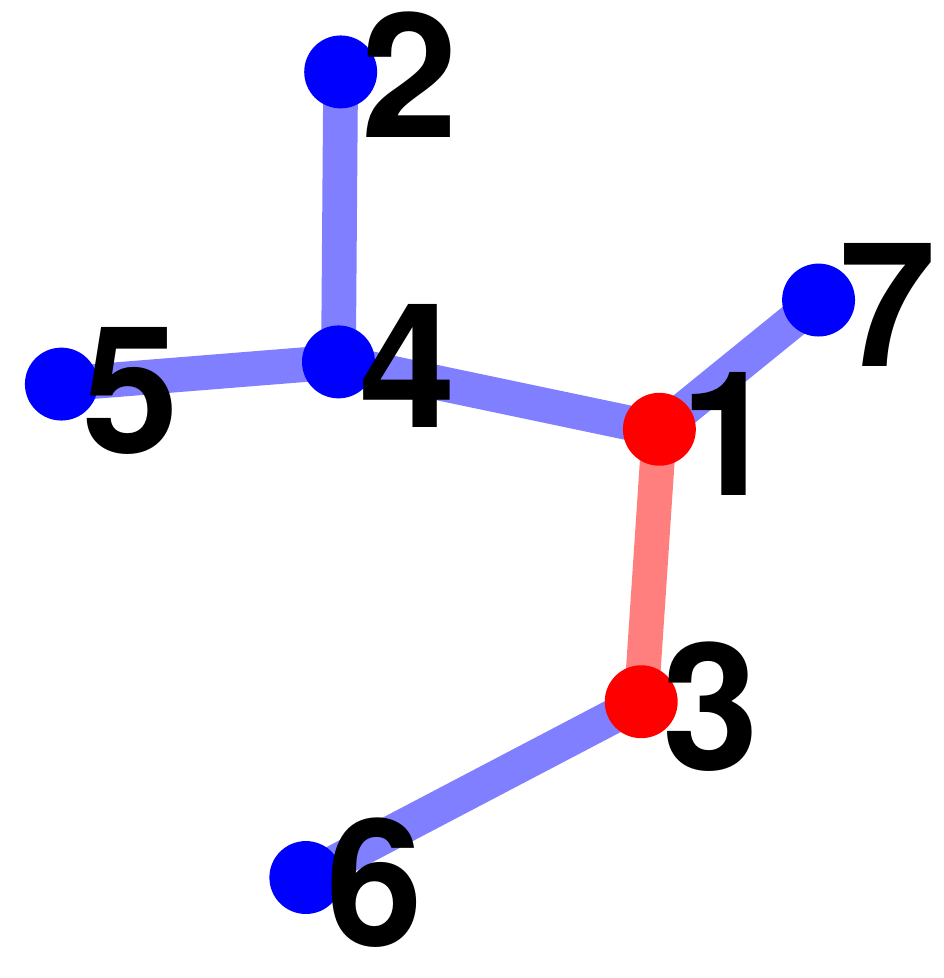} &  \includegraphics[height=0.6in,trim={6.5cm 9.5cm 5.5cm 8.5cm},clip] {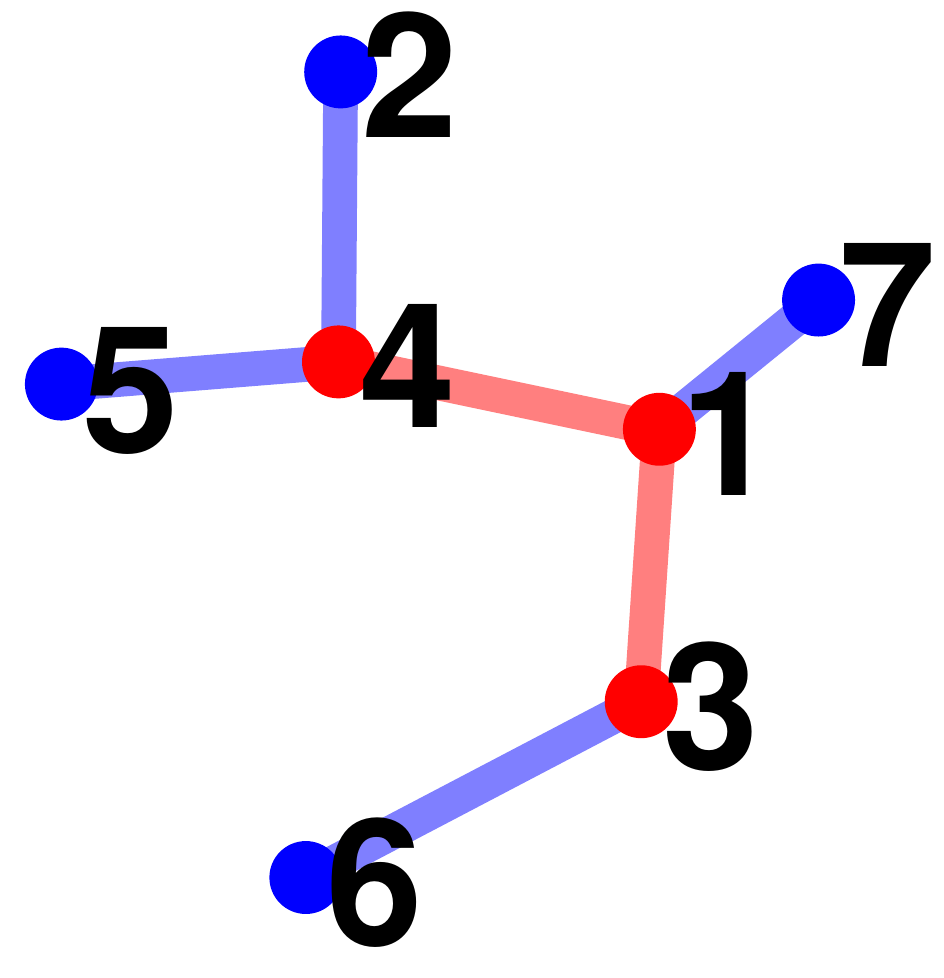} &  \includegraphics[height=0.6in,trim={6.5cm 9.5cm 5.5cm 8.5cm},clip] {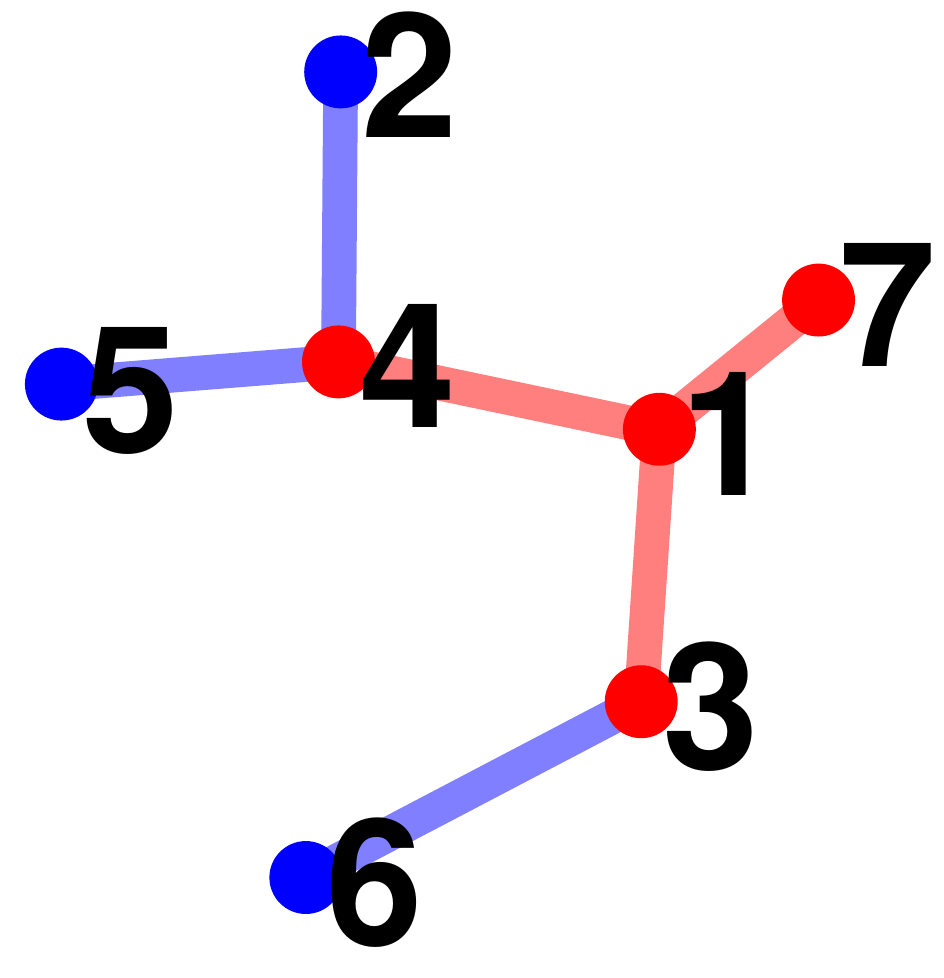} &    \includegraphics[height=0.6in,trim={6.5cm 9.5cm 5.5cm 8.5cm},clip] {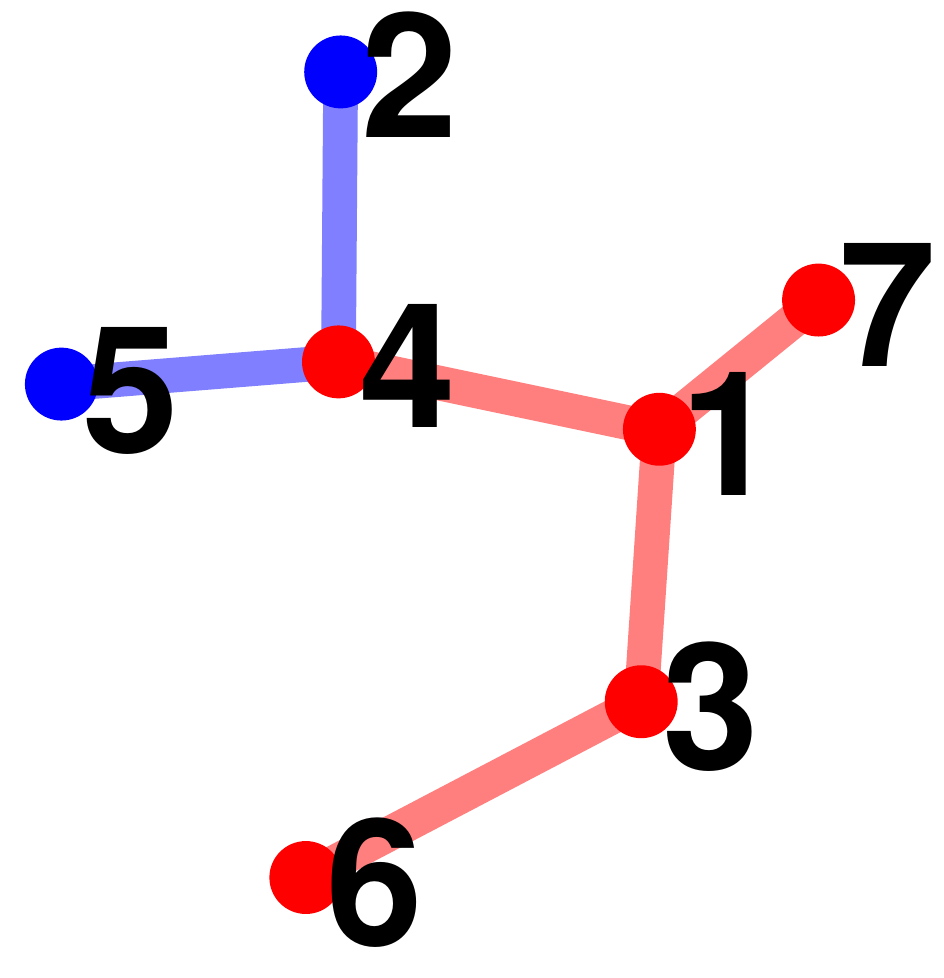} &  \includegraphics[height=0.6in,trim={6.5cm 9.5cm 5.5cm 8.5cm},clip] {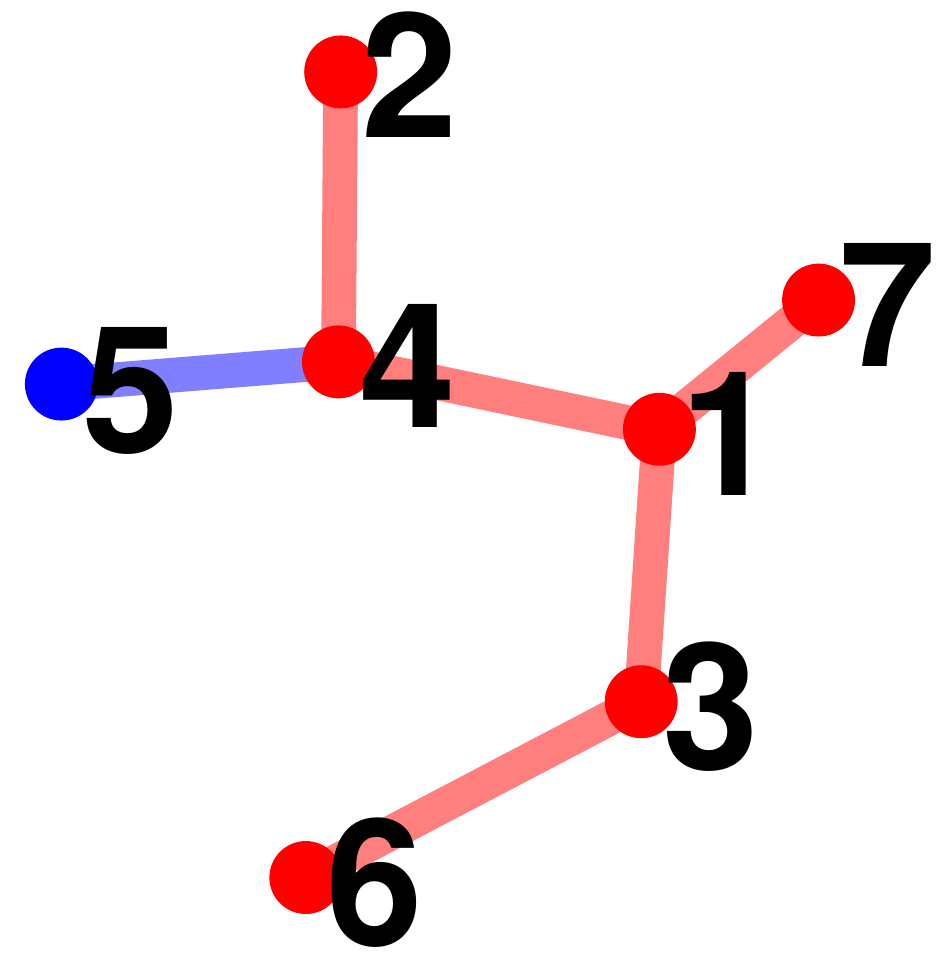} &  \includegraphics[height=0.6in,trim={6.5cm 9.5cm 5.5cm 8.5cm},clip] {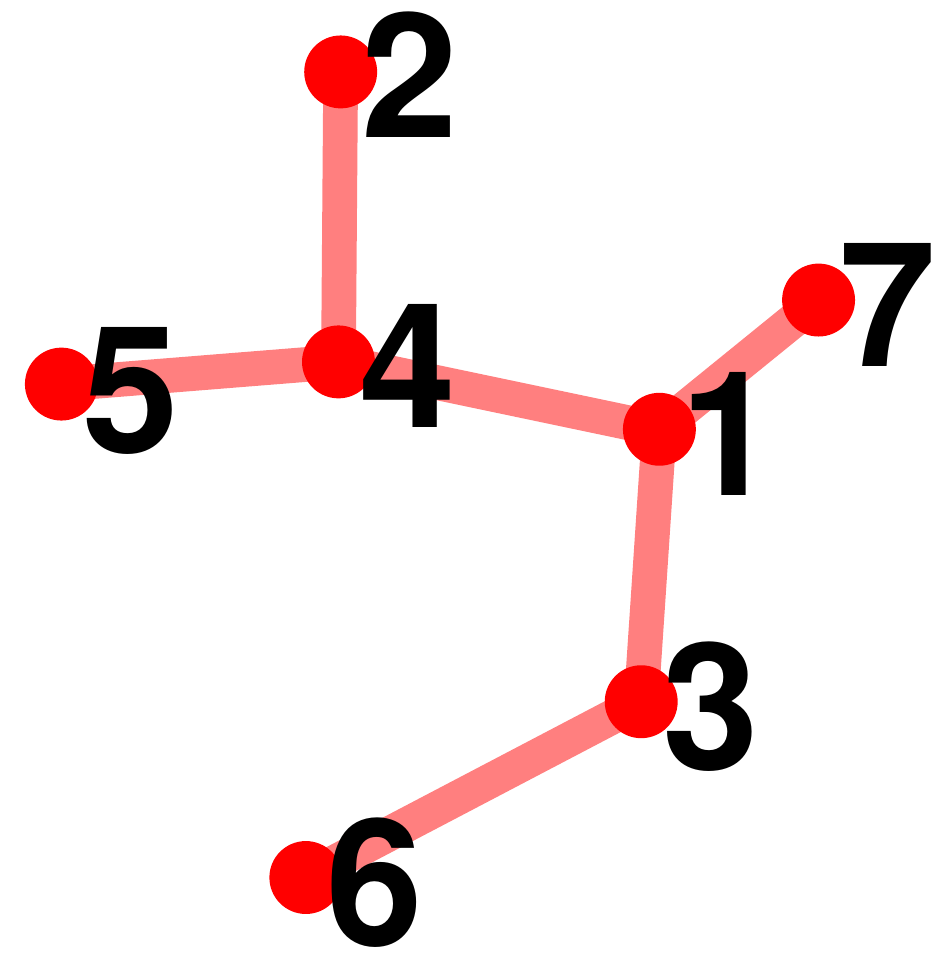} \\
        (a) & (b) & (c) & (d) & (e) & (f) & (g)
    \end{tabular}
    \caption{An illustration of DT($\mathcal{X}$), $\mathcal{T}(\mathcal{X})$, and the corresponding $P$ for $\mathcal{X}\in\mathbb{R}^{7\times 2}$. (a) DT($\mathcal{X}$) (black dash line) and $\mathcal{T}(\mathcal{X})$ (in blue). (b) Starting from the root (Node $1$), find the first undiscovered node, e.g., Node $3$ with depth $1$, then let $P = [1, 3]$. (c) Add $[1, 4]$ to $P$. (d) Add $[1, 7]$ to $P$. (e) find the first undiscovered node, e.g., Node $6$ with depth $2$, then add $[3, 6]$ to $P$. (f) Add $[4, 2]$ to $P$. (g) Add $[4, 5]$ to $P$. Finally, a recovery path matrix $P \in \mathbb{R}^{6 \times 2}$ is set to be $P=[1, 3; 1, 4; 1, 7; 3, 6; 4, 2; 4, 5]$.}
    \label{fig:order}
\end{figure}

Figure~\ref{fig:order} visualizes an example of DT($\mathcal{X}$) and $\mathcal{T}(\mathcal{X})$ for $\mathcal{X}\in\mathbb{R}^{7\times 2}$. The process of constructing $P$ by the Breadth-First search algorithm is illustrated as well. The whole algorithm is summarized in Algorithm~\ref{alg:rpm}.

\begin{algorithm2e}
    
    \Fn{$P = $ \rpm{$\mathcal{X}$}}{
    
    $\mathcal{G} \leftarrow $ \dt{$\mathcal{X}$};
    
    $\mathcal{T} \leftarrow $ \mst{$\mathcal{G}$};
    
    $P \leftarrow $ \bfs{$\mathcal{T}$};
    }
    
    \caption{An $\O{N \L{N}}$ algorithm for generating a recovery path matrix $P$.}
    \label{alg:rpm}
    
\end{algorithm2e}

\subsubsection{Matrix recovery}
\label{sec:mr}

When the vector recovery algorithms in Algorithm~\ref{alg:rv2} and Algorithm~\ref{alg:rpm} are ready, we apply them to design a matrix recovery algorithm. Recall that the main idea is to identify piecewise smooth rows and columns of $\Psi$ satisfying \eqref{eqn:tr} as summarized in an informal problem statement in \eqref{eqn:smphi}. Let $\mathcal{X}_1$ and $ \mathcal{X}_2 \in \mathbb{R}^{N \times d}$ store the spatial locations of the $N$ grid points for the discretization of $\Phi(x,\xi)$ in $x$ and $\xi$, respectively. 

First, Algorithm~\ref{alg:rpm} is applied to construct the recovery path matrices $P_1$ and $P_2$ corresponding to $\mathcal{X}_1$ and $ \mathcal{X}_2$, respectively. Then the matrix recovery problem can be formally stated as
\begin{equation}
    \label{eqn:mintv2}
    \begin{split}
        \smash{\displaystyle\min_{\Phi\in\mathbb{R}^{N\times N}}} &\quad \sum_{i\in\mathcal{R}}\sum_{s\in \{1,\dots,N-1\}\setminus \mathcal{D}_c}|\Phi(i,P_2(s,2))-\Phi(i,P_2(s,1))| \\
        &\quad + \sum_{j\in\mathcal{C}}\sum_{t\in \{1,\dots,N-1\}\setminus \mathcal{D}_r}|\Phi(P_1(t,2),j)-\Phi(P_1(t,1),j)|\\
        \text{subject to} &\quad \bmod(\Phi(i,j),1) = \frac{1}{2\pi} \Im\left(\L{K(i,j)}\right) \text{ for }i\in \mathcal{R} \text{ or } j\in\mathcal{C},
    \end{split}
\end{equation}
where $\mathcal{D}_c$ and $\mathcal{D}_r$ are index sets indicating the discontinuous locations of $\Phi$ along columns and rows,  $\mathcal{R}$ and $\mathcal{C}$ are row and column index sets with $\O{1}$ randomly selected indices, respectively.

Next, Algorithm~\ref{alg:rv2} is applied with $\tau$ to identify the sets of discontinuous points $\mathcal{D}_r$ and $\mathcal{D}_c$ to make \eqref{eqn:mintv2} self-contained. Similarly to the 1D case, we can partition the phase matrix into (usually non-contiguous) submatrices corresponding to the domains in which the phase matrix is continuous, which is equivalent to dividing the MST $\mathcal{T}(\mathcal{X})$ into subtrees whenever an edge connects a predecessor node considered as a discontinuous point. Correspondingly, the recovery path matrix is partitioned into submatrices associated with these subtrees. Figure~\ref{fig:segment} below visualizes an example when an MST $\mathcal{T}(\mathcal{X})$ is partitioned into two MSTs at the discontinuity location at Node $4$. 

The partition procedure is denoted as Function \pibm in Algorithm~\ref{alg:rm2}, resulting in  $n_r \times n_c$ submatrices of the phase matrix denoted as $\Phi.\mathcal{B}_s \mathcal{B}_t$, $n_r$ submatrices of the recovery path matrix $P_1$, and $n_s$ submatrices of the recovery path matrix $P_2$, for $s = 1, 2, \dots, n_r$, and $t = 1, 2, \dots, n_c$. The random samples of the row and column indices in the submatrices are denoted as $\mathcal{R}.\mathcal{B}_s$ and $\mathcal{C}.\mathcal{B}_t$, respectively. For example, Panel (a) in Figure~\ref{fig:trace2} visualizes an example when the phase function contains only 4 continuous submatrices (from light color to dark color): $\Phi.\mathcal{B}_1 \mathcal{B}_1$, $\Phi.\mathcal{B}_1 \mathcal{B}_2$, $\Phi.\mathcal{B}_2 \mathcal{B}_1$, $\Phi.\mathcal{B}_2 \mathcal{B}_2$. Panel (b) in Figure~\ref{fig:trace2} visualizes the root row and the root column of each submatrix. Panel (c) and (d) in Figure~\ref{fig:trace2} visualize the randomly selected rows $\mathcal{R}.\mathcal{B}_1$ and columns $\mathcal{C}.\mathcal{B}_1$ in $\Phi.\mathcal{B}_1 \mathcal{B}_1$.

\begin{minipage}{1\linewidth}
    \hspace{-0.025\linewidth}
    \begin{minipage}{0.325\linewidth}
        \begin{figure}[H]
            \centering
            \begin{tabular}{ccc}
                \includegraphics[height=0.6in,trim={6.5cm 9.5cm 5.5cm 8.5cm},clip] {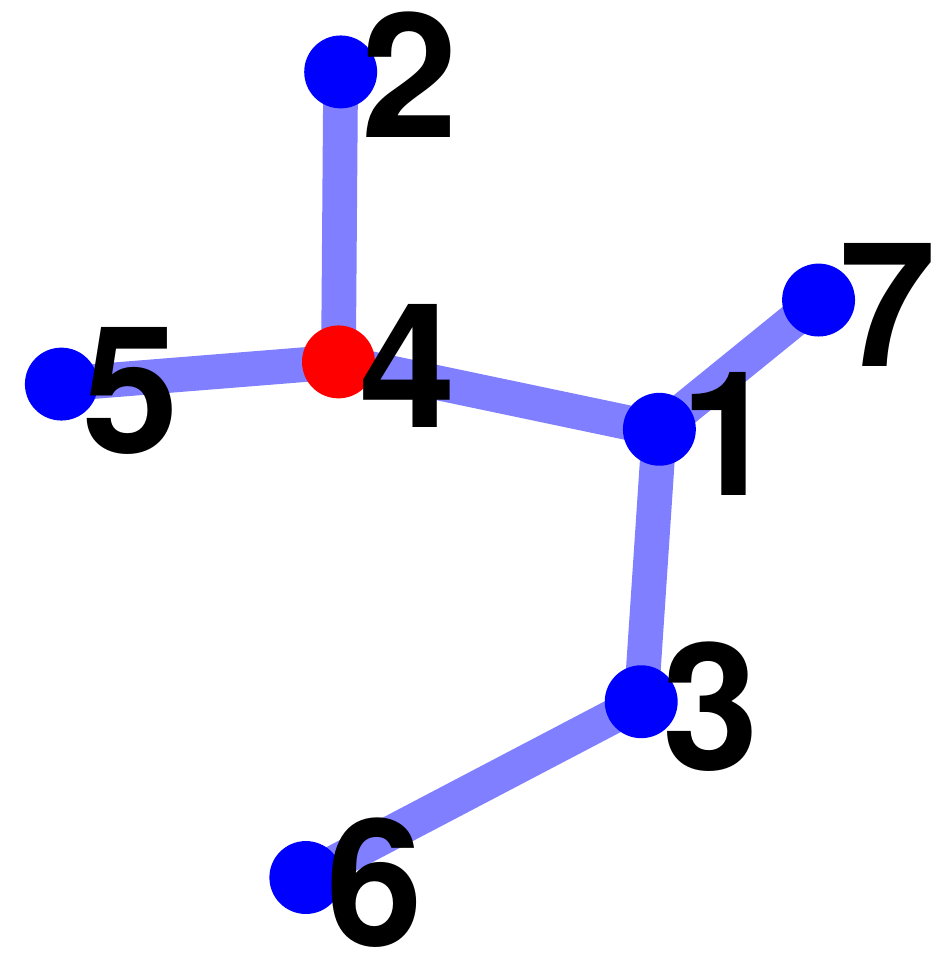} &  \includegraphics[height=0.6in,trim={6.5cm 9.5cm 5.5cm 8.5cm},clip] {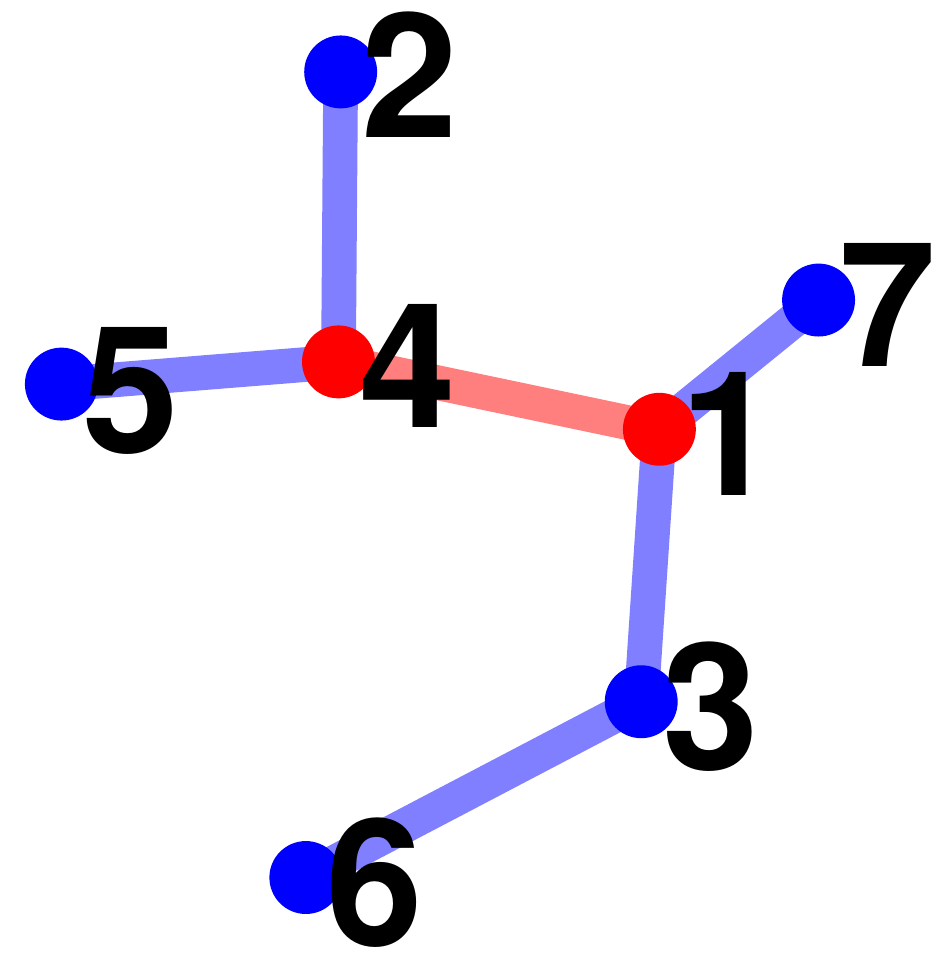} &  \includegraphics[height=0.6in,trim={6.5cm 9.5cm 5.5cm 8.5cm},clip] {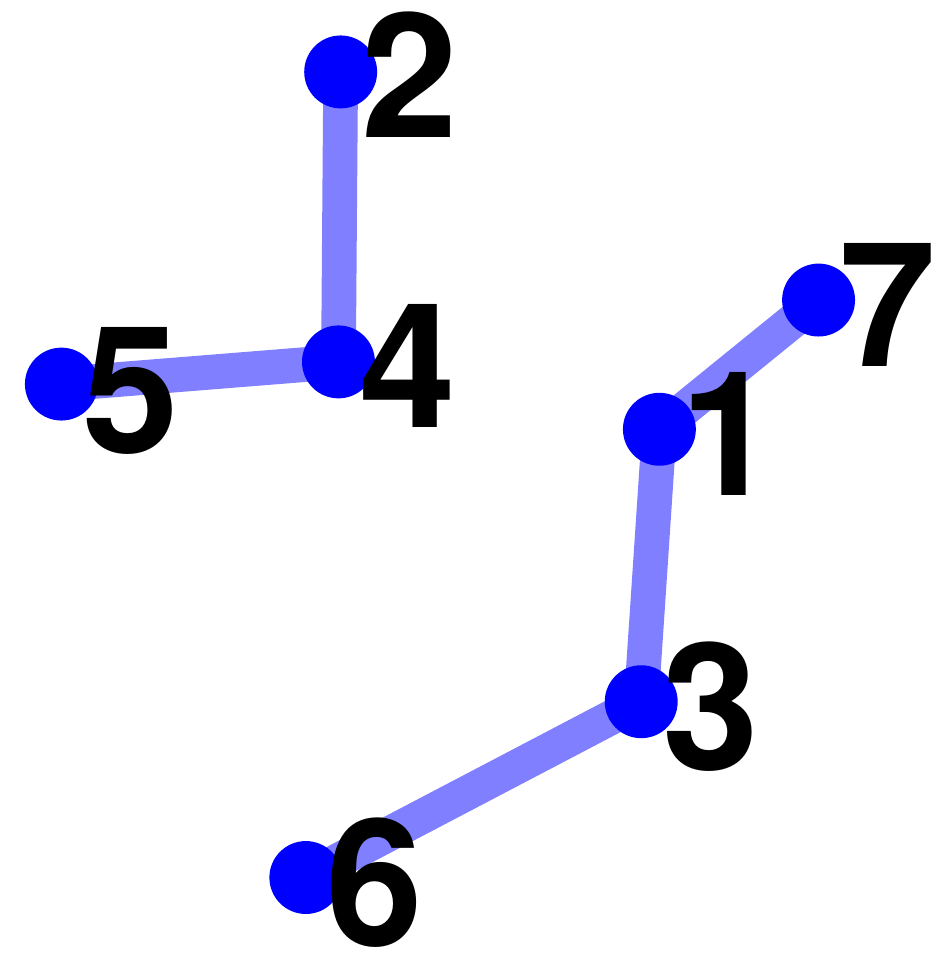} \\
                (a) & (b) & (c)
            \end{tabular}
            \caption{(a) An MST $\mathcal{T}(\mathcal{X})$ with a discontinuity location at Node $4$. (b) Separate $\mathcal{T}(\mathcal{X})$ at the edge between Node $4$ and its predecessor Node $1$. (c) Two resulting subtrees and the corresponding recovery path matrices $[1, 3; 1, 7;3, 6]$ and $[4, 2; 4, 5]$.}
            \label{fig:segment}
        \end{figure}
    \end{minipage}
    \hspace{0.075\linewidth}
    \begin{minipage}{0.55\linewidth}
        \hspace{-0.05\linewidth}
        \begin{figure}[H]
            \begin{tabular}{ccccccc}
                \includegraphics[height=0.6in]{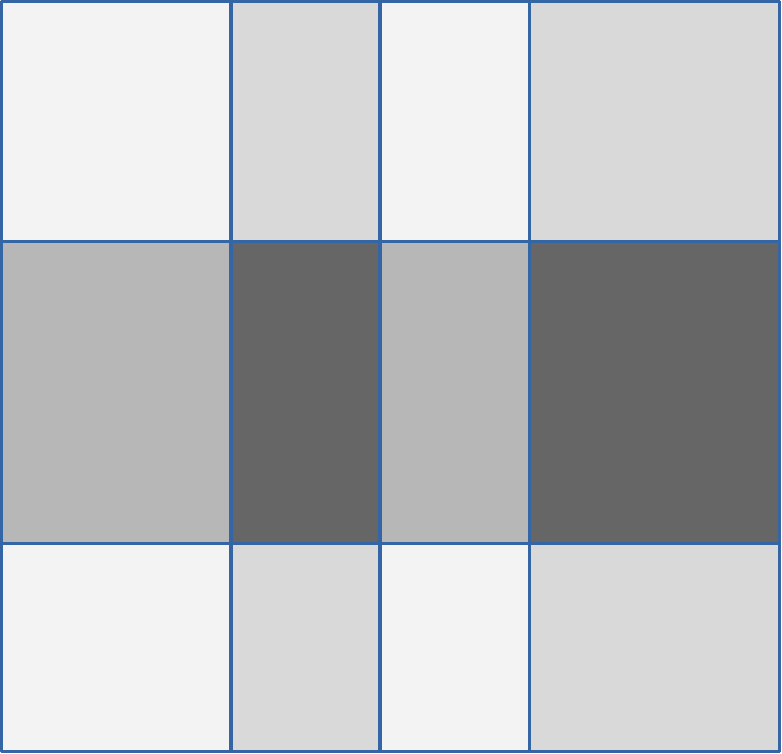} & & \includegraphics[height=0.6in]{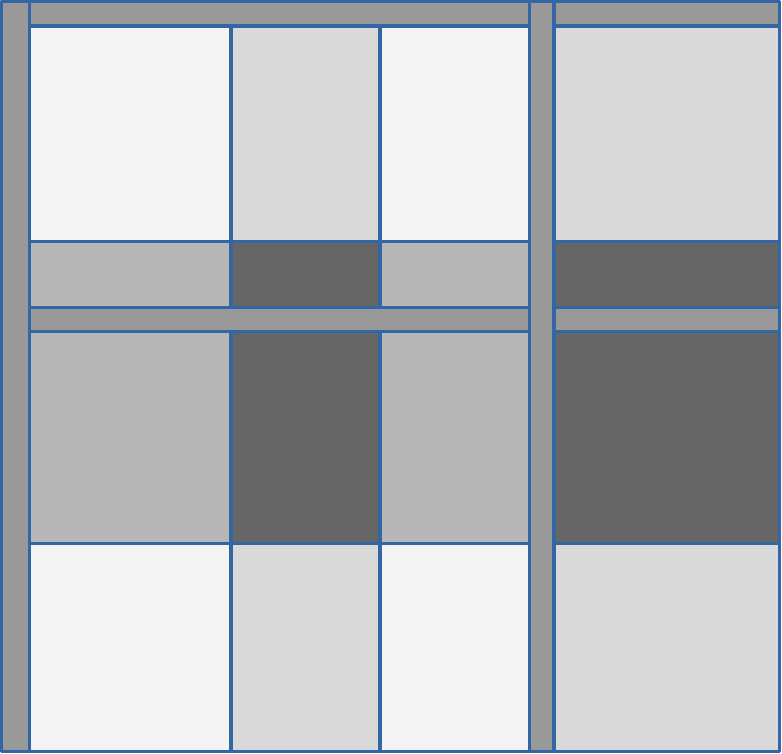} & & \includegraphics[height=0.6in]{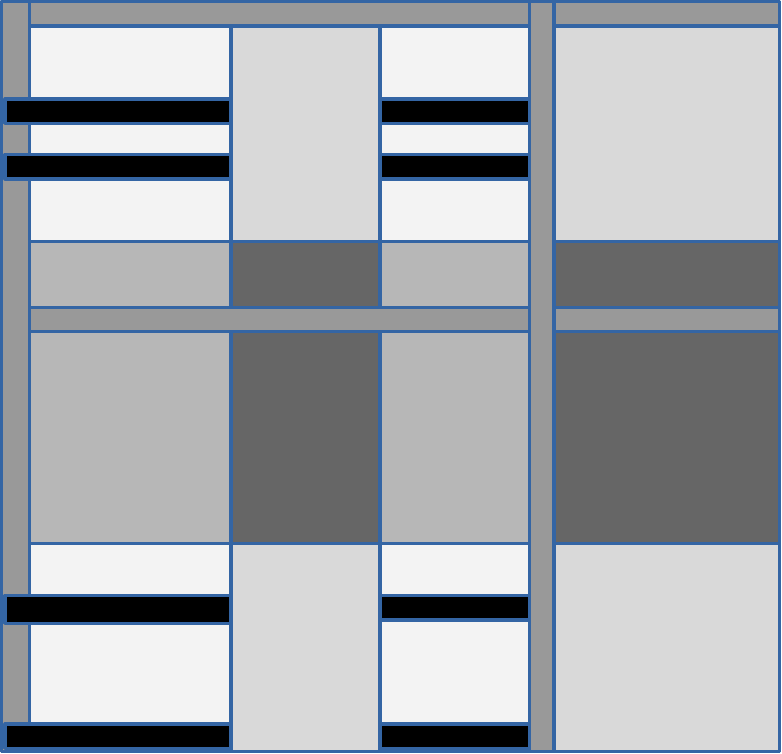} & & \includegraphics[height=0.6in]{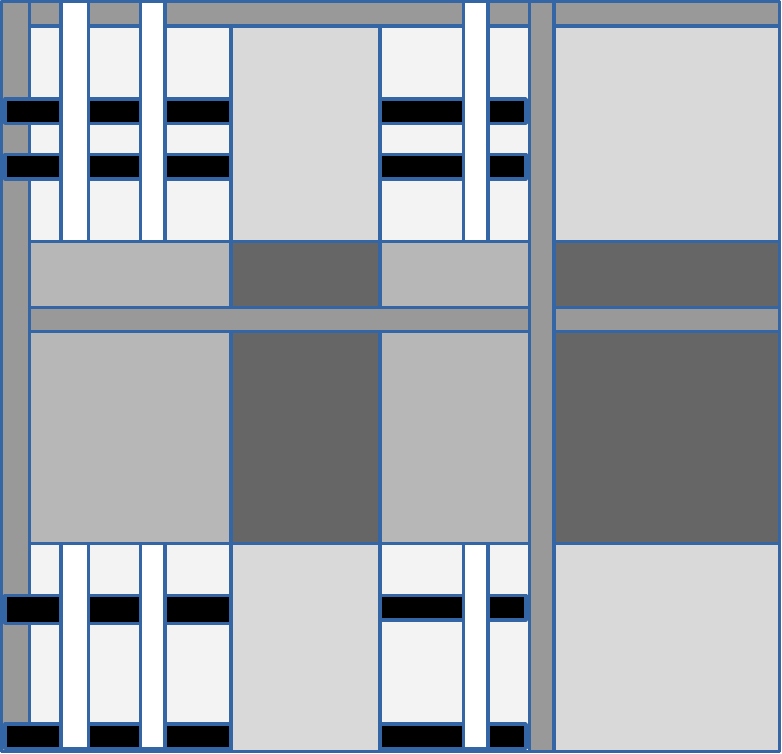} \\
                (a) & & (b) & & (c) & & (d)
            \end{tabular}
            \caption{An illustration of the low-rank matrix recovery for multidimensional phase matrix in Algorithm~\ref{alg:rm2}. (a) Line~\ref{alg:rm2:pib} partitions the phase matrix into 4 submatrices in 4 kinds of color such that there is no discontinuity along rows and columns in each submatrix. (b) Line~\ref{alg:rm2:rv21}-\ref{alg:rm2:rv22} recovers the row and column of each submatrix corresponding to the root node of each sub-MST. (c) Line~\ref{alg:rm2:rv23} recovers $\O{1}$ rows of each submatrix. (d) Line~\ref{alg:rm2:rv24} recovers $\O{1}$ columns of each submatrix.}
            \label{fig:trace2}
        \end{figure}
    \end{minipage}
\end{minipage}
\vspace{0.3cm}

Finally, we apply Algorithm~\ref{alg:rv2} again to recover each submatrix. The parameter for detecting discontinuity is set to $1$ since there is no need to detect discontinuity. The specially designed order also guarantees that each recovered row and column at their intersection share the same value, as long as the discontinuous points in the phase function have already been well distinguished, as proved by Lemma~\ref{nLRthm} below.

\begin{lemma}
    \label{nLRthm}
    Given $\bmod(\phi,1) \in \mathbb{R}^{n \times m}$, {where $\phi$ is a $d$-dimensional phase matrix, $d = 2$ or $3$. }
    Assuming that all rows and columns of $\phi$ belong to the class $C_{\tau,P}$ with a threshold $\tau \le \frac{1}{4}$, then the intersection of each recovered row and column by Algorithm~\ref{alg:rm2} share the same value.
\end{lemma}

The proof of Lemma~\ref{nLRthm} is simple and similar to Lemma~\ref{1LRthm}. For simplicity, we leave the proof to the reader.

Recall that the correct $\tau$ depends on the phase function and is not known a priori in one-dimensional cases. In practice, $\tau$ can be set as $\frac{1}{4}$ for identifying discontinuous point, which can guarantee that the intersection of each recovered row and column share the same value. When the number of discontinuous points is too large, $\tau+\epsilon$ is used to identify discontinuous points, i.e. $\epsilon = \frac{1}{40}$. This procedure can be repeated until $\O{1}$ discontinuous points have been detected and it takes at most $\O{N}$ operations to obtain a reasonable $\tau$. When $\tau$ increases to $\frac{1}{2}$, no more discontinuous point will be detected.

In fact, if $\tau$ is set larger than $\frac{1}{4}$, the consistency of the intersection of each recovered row and column should be checked manually instead of by Lemma~\ref{nLRthm}. As previously said, our method is based on the first-order derivative of the phase function, the extension of Algorithm~\ref{alg:rv2} using the high-order finite difference schemes in \cite{JIANCHUN1995162,VASILYEV2000746} is left as future work if the recovered intersection values are not consistent. In our numerical tests for multidimensional cases, $\tau = \frac{1}{4}$ is good enough for all numerical examples.

When $\O{1}$ discontinuous points have been detected, Algorithm~\ref{alg:rv2} will recover $\O{1}$ randomly selected rows and columns of the phase matrix with nearly linear computational complexity.

Algorithm~\ref{alg:rm2} below summarizes the above steps and the whole process is illustrated in Figure~\ref{fig:trace2}.

\begin{algorithm2e}
    
    \Fn{$\Phi = $ \rmm{$\Phi,\mathcal{R},\mathcal{C},\mathcal{X}_1,\mathcal{X}_2, \tau$}}{
    
    $P_1 \leftarrow $ \rpm{$\mathcal{X}_1$} \qd 
    $P_2 \leftarrow $ \rpm{$\mathcal{X}_2$}
    
    $\mathcal{D}_r \leftarrow$ \rvm{$\Phi(:,1),\tau,P_1$}
    \tcp*[f]{$\mathcal{D}_r:$ discontinuous point set}
    
    $\mathcal{D}_c \leftarrow$ \rvm{$\Phi(1,:),\tau,P_2$}
    \tcp*[f]{$\mathcal{D}_c:$ discontinuous point set}
    
    $\mathcal{R} \leftarrow \left[\mathcal{R}, \mathcal{D}_r \right]$ \qd $\mathcal{C} \leftarrow \left[\mathcal{C}, \mathcal{D}_c \right]$
    
    $n_r \leftarrow $ \length($\mathcal{D}_r$) \qd $n_c \leftarrow $ \length($\mathcal{D}_c$)
    
    $\left[\Phi, \mathcal{R}, \mathcal{C}, P_1, P_2 \right] \leftarrow $ \pibm{$\Phi, \mathcal{R}, \mathcal{C}, P_1, P_2,  \mathcal{D}_r, \mathcal{D}_c$} \label{alg:rm2:pib}
    
    \For{$s = 1:n_r$}{
    
        \For{$t = 1:n_c$}{
        
            $\Phi.\mathcal{B}_s\mathcal{B}_t(1,:) \leftarrow$ \rvm{$\Phi.\mathcal{B}_s\mathcal{B}_t(1,:), 1, P_2.\mathcal{B}_t$} \label{alg:rm2:rv21}
            
            $\Phi.\mathcal{B}_s\mathcal{B}_t(:,1) \leftarrow$ \rvm{$\Phi.\mathcal{B}_s\mathcal{B}_t(:,1), 1, P_1.\mathcal{B}_s$} \label{alg:rm2:rv22}
            
            $\Phi.\mathcal{B}_s\mathcal{B}_t(\mathcal{R}.\mathcal{B}_s(k),:) \leftarrow$ \rvm{$\Phi(\mathcal{R}.\mathcal{B}_s(k),:),1,P_2.\mathcal{B}_t$} for all $k$ \label{alg:rm2:rv23}
    
            $\Phi.\mathcal{B}_s\mathcal{B}_t(:,\mathcal{C}.\mathcal{B}_t(k)) \leftarrow$ \rvm{$\Phi(:,\mathcal{C}.\mathcal{B}_t(k)),1,P_1.\mathcal{B}_s$} for all $k$ \label{alg:rm2:rv24}
        }
    }
    }
    
    \caption{An $\O{N \L{N}}$ algorithm for the solution of matrix recovery problem \eqref{eqn:smphi} when the phase function $\Phi(x,\xi)$ is defined on $\mathbb{R}^d \times \mathbb{R}^d$.} 
    \label{alg:rm2}
    
\end{algorithm2e}

\subsubsection{Phase matrix factorization}
\label{sec:nmLR}

Once the phase function recovery algorithm in Algorithm~\ref{alg:rm2} is ready, following the idea of low-rank matrix factorization via randomized sampling in Algorithm~\ref{alg:rSVD}, we can introduce a nearly linear scaling algorithm to construct the low-rank factorization of the phase matrix as summarized in Algorithm~\ref{alg:lrf}. In particular, Algorithm~\ref{alg:lrf} constructs a low-rank factorization $UV^T$, where $U\in\mathbb{C}^{N\times r}$ and $V\in\mathbb{C}^{N\times r}$, such that $e^{2\pi\i UV^T}\approx e^{2\pi\i  \Phi}$ when we only know the kernel matrix $K=e^{2\pi\i \Phi}$ through Scenarios $1$ and $2$ in Table~\ref{tab:sc}.

In Algorithm~\ref{alg:lrf}, $K$ (and $\Phi$) is a function handle for evaluating an arbitrary entry of the kernel matrix, or evaluating an arbitrary row or column of $K$ (and $\Phi$). Two coordinate matrices $\mathcal{X}_1, \mathcal{X}_2 \in \mathbb{R}^{N \times d}$, a rank parameter $r$, an over-sampling parameter $q$, and the matrix size $N$ are also inputs. We randomly select $rq$ rows and columns of the kernel matrix and use \rmm to obtain the corresponding rows and columns of $\Psi$ such that $e^{2\pi\i \Psi} \approx K$. Finally, apply Function \rsvd in Algorithm~\ref{alg:rSVD} in Subsection~\ref{sec:LRF} to evaluate the low-rank factorization of $\Psi\approx UV^T$ such that $e^{2\pi\i UV^T}\approx K=e^{2\pi\i \Phi}$. The reconstructed phase matrix can be set as an initial guess to the optimization problem in \eqref{eqn:mintv2} and it takes $\O{1}$ iterations for sub-gradient descent methods to converge.

\begin{algorithm2e}
    
    \Fn{$\left[ U, V \right] = $ \lrf{$K, \mathcal{X}_1, \mathcal{X}_2, r, q, N$}}{
    
    $\mathcal{R} \leftarrow \rp{N,rq}$ \qd $\mathcal{C} \leftarrow \rp{N,rq}$
    
    $\Phi \leftarrow \frac{1}{2\pi} \Im\left(\L{K}\right)$\tcp*[f]{generate a function handle for the evaluation of $\Phi$}
    
    $\Psi \leftarrow \rmm{$\Phi,\mathcal{R},\mathcal{C},\mathcal{X}_1,\mathcal{X}_2$}$ \\
    \tcp*[f]{generate a function handle for the evaluation of $\Psi$} \label{alg:lrf:rm2}
    
    $\left[U,\Sigma,V\right] \leftarrow$ \rsvd{$\Psi, \mathcal{R}, \mathcal{C}, r$} \label{alg:lrf:rSVD}
    
    $V \leftarrow V\Sigma$
    
    }
    
    \caption{An $\O{N \L{N}}$ algorithm for low-rank matrix factorization of phase functions in the case of indirect access.}
    \label{alg:lrf}
    
\end{algorithm2e}

\subsubsection{Summary}

Before moving to the next algorithm, let us summarize how those algorithms in Subsection~\ref{sec:nLR} can be applied to construct the low-rank matrix factorization of the multidimensional phase functions with nearly linear computational complexity.

For a general kernel function $K(x,\xi)=e^{2\pi\i \Phi(x,\xi)}$, suppose we discretize $\Phi(x,\xi)$ with $N$ grid points in each variable to obtain the phase matrix $\Phi$. When the explicit formulas of $\Phi(x,\xi)$ are known, it takes $\O{N}$ operations to evaluate one column or one row of $\Phi$. Then, the randomized SVD in Subsection~\ref{sec:LRF} is able to construct the low-rank matrix factorization of $\Phi$ in $\O{N}$ operations.

When the explicit formulas are unknown such as in Scenario $2$, it takes $\O{N \L{N}}$ operations to evaluate one column or one row of the kernel matrix $K$. Hence, the phase recovery and the low-rank factorization of $\Phi$ can be constructed by Algorithm~\ref{alg:rm2} and Algorithm~\ref{alg:lrf} in $\O{N \L{N}}$ operations. 

In the case of indirect access in Scenario $3$, $\O{1}$ columns and rows of and phase functions are available by solving certain PDE's. For example, in practical applications like solving wave equations \cite{Yingwave}, each column or row can be obtained via interpolating the solution of the PDE on a coarse grid of size independent of $N$. Thus, the phase recovery algorithm is not required for Scenario $3$, it only needs to construct a low-rank factorization of $\Phi$ by Algorithm~\ref{alg:lrf} in $\O{N \L{N}}$ operations.

For Scenario $1$, which is a special case included in Scenario $2$, any arbitrary entry of the kernel matrix is available in $\O{1}$ operations. Therefore, it can be applied directly to the next algorithm.

Since Line~\ref{alg:lrf:rm2} in Algorithm~\ref{alg:lrf} identifies $\O{1}$ rows and columns of a low-rank matrix $\Psi$ such that $\bmod(\Psi(i,j),1) = \frac{1}{2\pi} \Im\left(\L{K(i,j)}\right)$ for $i\in \mathcal{R}$ or $j\in\mathcal{C}$, there is not any error generated in this step. The approximation error of Algorithm~\ref{alg:lrf} is $\O{\epsilon}$, which is caused by the low-rank approximation algorithm (Line~\ref{alg:lrf:rSVD}).

\section{Multidimensional Interpolative Decomposition Butterfly Factorization (MIDBF)} \label{sec:MIDBF}

This section will introduce the multidimensional interpolative decomposition butterfly factorization for a matrix $K = (K(x,\xi))_{x\in X,\xi\in \Omega}$ satisfying a complementary low-rank property \cite{BF}, where $X$ and $\Omega$ contain $\O{N}$ points possibly non-uniformly distributed in $[0,1)^d$ and $d$ is the dimension of the domain. As a special example, the kernel matrix $K(x,\xi)=e^{2\pi i \Phi(x,\xi)}$ satisfies the complementary low-rank property. Hence, once the phase function $\Phi$ in Scenarios $2$ and $3$ has been recovered by Algorithm~\ref{alg:lrf} in Subsection~\ref{sec:nmLR} in the form of low-rank factorization, we can construct a function handle to evaluate an arbitrary entry of $K$ in $\O{1}$ operations. Especially, in Scenario $1$, this kind of function handle is known directly. Then, the MIDBF can construct the butterfly factorization of $K$ for nearly linear scaling fast matvec, when the function handle is given.

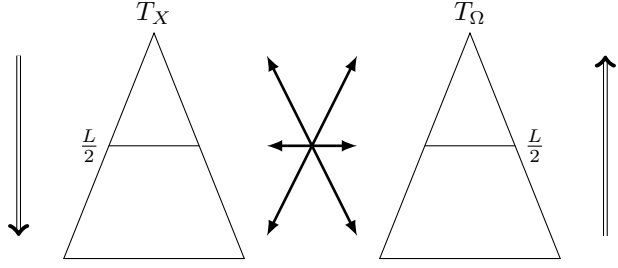
\begin{figure}[!ht]
\centering
\begin{tikzpicture}[scale=0.3]

\draw [->,double,double distance=1pt] (-6,9) -> (-6,1);

\draw (0,10) -- (-4,0) -- (4,0) -- cycle;

\draw [->,double,double distance=1pt] (20,1) -> (20,9);

\draw (14,10) -- (10,0) -- (18,0) -- cycle;

\draw [latex-latex,line width=1pt] (9,5) -> (5,5);
\draw [latex-latex,line width=1pt] (9,9) -> (5,1);
\draw [latex-latex,line width=1pt] (9,1) -> (5,9);

\coordinate [label=above:$T_X$] (X) at (0,10);
\coordinate [label=above:$T_\Omega$] (X) at (14,10);

\draw (-2,5) -- (2,5);
\draw (12,5) -- (16,5);

\coordinate [label=right:$\frac{L}{2}$] (X) at (16,5);
\coordinate [label=left:$\frac{L}{2}$] (X) at (-2,5);
\end{tikzpicture}
\caption{Trees of the row and column indices.
    Left: $T_X$ for the row indices $X$.
    Right: $T_\Omega$ for the column indices $\Omega$.
    The interaction between $A\in T_X$ and $B\in T_\Omega$
    starts at the root of $T_X$ and the leaves of $T_\Omega$. }
\label{fig:domain-tree-BF}
\end{figure}

Let us recall the definition of complementary low-rank matrices in \cite{BF}. For such a matrix, we construct two trees $T_X$ and $T_\Omega$ for point sets $X$ and $\Omega$, respectively, assuming that both trees have
the same depth $L=\O{\L{N}}$, with the top-level being level $0$ and
the bottom one being level $L$ (see Figure~\ref{fig:domain-tree-BF} for an illustration). Such a matrix $K$ of size $N\times N$
is said to satisfy the {\bf complementary low-rank property} if for
any level $\ell$, any node $A$ in $T_X$ at level $\ell$, and any node
$B$ in $T_\Omega$ at level $L-\ell$, the submatrix $K(A,B)$, obtained
by restricting $K$ to the rows indexed by the points in $A$ and the
columns indexed by the points in $B$, is numerically low-rank.

\subsection{Notations and overall structure}

The notation of the 1D IDBF introduced in \cite{IDBF} will be adopted and adjusted to the multidimensional case in this paper. With no loss of generality, we focus on the 2D case with uniform point distributions first. The notations and overall structure discussed below are similar to that in \cite{MBF,IDBF}.

\begin{figure}[!ht]
\begin{minipage}{.33\textwidth}
\centering
\begin{tikzpicture}[scale=3]
\fill[white] (0,0) rectangle (1,1);
\draw[black] (0,0) rectangle (1,1);
\draw (0.5,0.5) node[rectangle] {$A^0_{{0}}$};
\draw[->,thick] (0,1.1) -> (1,1.1);
\draw (0.5,1.1) node[rectangle,above] {$x_2$};
\draw[->,thick] (-0.1,1) -> (-0.1,0);
\draw (-0.1,0.5) node[rectangle,left] {$x_1$};
\end{tikzpicture}
\end{minipage}%
\begin{minipage}{.33\textwidth}
\centering
\begin{tikzpicture}[scale=3]
\fill[white] (0,0) rectangle (0.5,0.5);
\draw[black] (0,0) rectangle (0.5,0.5);
\draw (0.25,0.25) node[rectangle] {$A^1_{2}$};
\fill[white] (0,0.5) rectangle (0.5,1);
\draw[black] (0,0.5) rectangle (0.5,1);
\draw (0.25,0.75) node[rectangle] {$A^1_{0}$};
\fill[white] (0.5,0) rectangle (1,0.5);
\draw[black] (0.5,0) rectangle (1,0.5);
\draw (0.75,0.25) node[rectangle] {$A^1_{3}$};
\fill[white] (0.5,0.5) rectangle (1,1);
\draw[black] (0.5,0.5) rectangle (1,1);
\draw (0.75,0.75) node[rectangle] {$A^1_{1}$};
\draw[lightgray,thick] (0.25,0.75) -- (0.75,0.75) -- (0.25,0.25) -- (0.75,0.25);
\draw[->,thick] (0,1.1) -> (1,1.1);
\draw (0.5,1.1) node[rectangle,above] {$x_2$};
\draw[->,thick] (-0.1,1) -> (-0.1,0);
\draw (-0.1,0.5) node[rectangle,left] {$x_1$};
\end{tikzpicture}
\end{minipage}%
\begin{minipage}{.33\textwidth}
\centering
\begin{tikzpicture}[scale=3]
\fill[white] (0,0.75) rectangle (0.25,1);
\draw[black] (0,0.75) rectangle (0.25,1);
\draw (0.125,0.875) node[rectangle] {$A^2_{0}$};
\fill[white] (0.25,0.75) rectangle (0.5,1);
\draw[black] (0.25,0.75) rectangle (0.5,1);
\draw (0.375,0.875) node[rectangle] {$A^2_{1}$};
\fill[white] (0,0.5) rectangle (0.25,0.75);
\draw[black] (0,0.5) rectangle (0.25,0.75);
\draw (0.125,0.625) node[rectangle] {$A^2_{2}$};
\fill[white] (0.25,0.5) rectangle (0.5,0.75);
\draw[black] (0.25,0.5) rectangle (0.5,0.75);
\draw (0.375,0.625) node[rectangle] {$A^2_{3}$};
\fill[white] (0.5,0.75) rectangle (0.75,1);
\draw[black] (0.5,0.75) rectangle (0.75,1);
\draw (0.625,0.875) node[rectangle] {$A^2_{4}$};
\fill[white] (0.75,0.75) rectangle (1,1);
\draw[black] (0.75,0.75) rectangle (1,1);
\draw (0.875,0.875) node[rectangle] {$A^2_{5}$};
\fill[white] (0.5,0.5) rectangle (0.75,0.75);
\draw[black] (0.5,0.5) rectangle (0.75,0.75);
\draw (0.625,0.625) node[rectangle] {$A^2_{6}$};
\fill[white] (0.75,0.5) rectangle (1,0.75);
\draw[black] (0.75,0.5) rectangle (1,0.75);
\draw (0.875,0.625) node[rectangle] {$A^2_{7}$};
\fill[white] (0,0.25) rectangle (0.25,0.5);
\draw[black] (0,0.25) rectangle (0.25,0.5);
\draw (0.125,0.375) node[rectangle] {$A^2_{8}$};
\fill[white] (0.25,0.25) rectangle (0.5,0.5);
\draw[black] (0.25,0.25) rectangle (0.5,0.5);
\draw (0.375,0.375) node[rectangle] {$A^2_{9}$};
\fill[white] (0,0) rectangle (0.25,0.25);
\draw[black] (0,0) rectangle (0.25,0.25);
\draw (0.125,0.125) node[rectangle] {$A^2_{10}$};
\fill[white] (0.25,0) rectangle (0.5,0.25);
\draw[black] (0.25,0) rectangle (0.5,0.25);
\draw (0.375,0.125) node[rectangle] {$A^2_{11}$};
\fill[white] (0.5,0.25) rectangle (0.75,0.5);
\draw[black] (0.5,0.25) rectangle (0.75,0.5);
\draw (0.625,0.375) node[rectangle] {$A^2_{12}$};
\fill[white] (0.75,0.25) rectangle (1,0.5);
\draw[black] (0.75,0.25) rectangle (1,0.5);
\draw (0.875,0.375) node[rectangle] {$A^2_{13}$};
\fill[white] (0.5,0) rectangle (0.75,0.25);
\draw[black] (0.5,0) rectangle (0.75,0.25);
\draw (0.625,0.125) node[rectangle] {$A^2_{14}$};
\fill[white] (0.75,0) rectangle (1,0.25);
\draw[black] (0.75,0) rectangle (1,0.25);
\draw (0.875,0.125) node[rectangle] {$A^2_{15}$};
\draw[lightgray, thick]
                 (0.125,0.875) -- (0.375,0.875) -- (0.125,0.625) -- (0.375,0.625)
              -- (0.625,0.875) -- (0.875,0.875) -- (0.625,0.625) -- (0.875,0.625)
              -- (0.125,0.375) -- (0.375,0.375) -- (0.125,0.125) -- (0.375,0.125)
              -- (0.625,0.375) -- (0.875,0.375) -- (0.625,0.125) -- (0.875,0.125);

\draw[->,thick] (0,1.1) -> (1,1.1);
\draw (0.5,1.1) node[rectangle,above] {$x_2$};
\draw[->,thick] (-0.1,1) -> (-0.1,0);
\draw (-0.1,0.5) node[rectangle,left] {$x_1$};
\end{tikzpicture}
\end{minipage}%
\caption{An illustration of Z-order curve across levels. The
  superscripts indicate the different levels while the subscripts
  indicate the index in the Z-ordering. The light gray lines show the
  ordering among the subdomains on the same level.  Left: The root at
  level $0$.  Middle: At level $1$, the domain $A^0_{0}$ is divided
  into $2\times 2$ subdomains $A^1_\i$ with $\i\in \calI^1=\{0,1,2,3\}$.
  These $4$ subdomains are ordered according to the Z-ordering. Right:
  At level $2$, the domain $A^0_0$ is divided into $4\times 4$
  subdomains $A^2_\i$ with $\i\in \calI^2=\{0,1,\ldots,15\}$. These $16$
  subdomains are ordered similarly.}
\label{fig:domain-order-2D}
\end{figure}
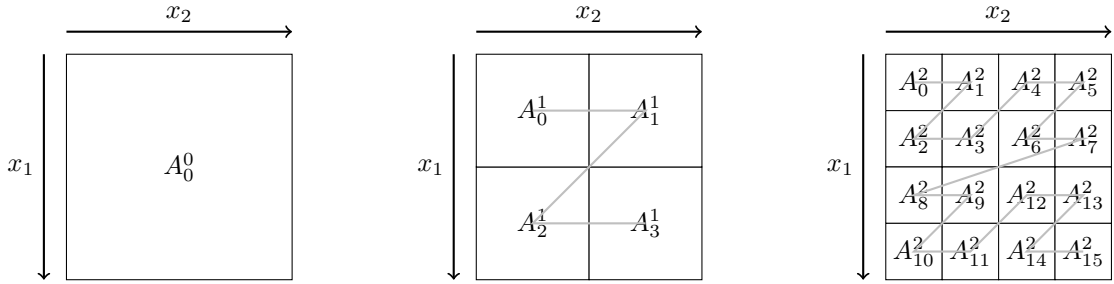

Recall that $n$ is the number of grid points on each dimension, $N = n^2 = 4^{L} n_0$ is the total number of points, $n_0 = \O{1}$ is the number of row or column indices in a leaf in the quadtrees of row and column spaces and, without loss of generality, $L$ is an even integer, i.e. $T_X$ and $T_\Omega$ with $L$ levels.  For a fixed level $\ell$ between $0$ and $L$, the quadtree $T_X$ has $4^\ell$ nodes at level $\ell$. By defining $\calI^\ell = \{0,1,\ldots,4^\ell-1\}$, we denote these nodes by $A^\ell_i$ with $i\in \calI^\ell$. These $4^\ell$ nodes at level $\ell$ are further ordered according to a Z-order curve (or Morton order) as illustrated in Figure~\ref{fig:domain-order-2D}. Based on this Z-ordering, the node $A^\ell_i$ at level $\ell$ has four child nodes denoted by $A^{\ell+1}_{4i + t}$ with $t=0,\dots,3$. The nodes plotted in Figure~\ref{fig:domain-order-2D} for $\ell = 1$ (middle) and $\ell=2$ (right) illustrate the relationship between the parent node and its child nodes. Similarly, in the quadtree $T_\Omega$, the nodes at level $L-\ell$ are denoted as $B^{L-\ell}_j$ for $j\in\calI^{L-\ell}$.

For any level $\ell$ between $0$ and $L$, the kernel matrix $K$ can be partitioned into $\O{N}$ submatrices $K(A^\ell_i,B^{L-\ell}_j) := (K(x,\xi))_{x\in A^\ell_i,\xi\in B^{L-\ell}_j}$ for $i\in\calI^\ell$ and $j\in\calI^{L-\ell}$. For simplicity, we shall denote $K(A^\ell_i,B^{L-\ell}_j)$ as $K^{\ell}_{i,j}$, where the superscript $\ell$ denotes the level in the quadtree $T_X$. Because of the complementary low-rank property, every submatrix $K^\ell_{i,j}$ is numerically low-rank with the rank bounded by a uniform constant $r$ independent of $N$.

The multidimensional interpolative decomposition butterfly factorization for $K$ is a product of $\O{\L{N}}$ sparse matrices, each of which contains $\O{\frac{k^2}{n_0}N}$ nonzero entries as follows:
\begin{equation}
    K \approx U^{L}U^{L-1}\cdots U^{h} S^h V^{h}\cdots V^{L-1}V^{L},
\end{equation}
where $k$ is a local rank parameter, $h=\frac{L}{2}$, and the level $L$ is assumed to be even.

\subsection{Linear scaling Interpolative Decomposition (ID)}

This subsection introduces the linear scaling ID method in \cite{IDBF}. Suppose $K\in \bbC^{m \times n}$ has a numerical rank $k_\epsilon \ll \min\{m,n\}$, i.e., $K$ admits a rank $k_\epsilon$ factorization with $\epsilon$ relative approximation accuracy. Let $s$ be an index set containing $tk$ rows of $K$ chosen from the Mock-Chebyshev grids as in \cite{NUFFTorBF,Mock1,Mock2}, $t$ is an oversampling parameter, and $k$ is an empirical estimation of $k_\epsilon$. $s$ is empirically selected and gradually increased if not large enough. We apply the rank revealing thin QR to $K(s,:)$:
\[
K(s,:)\Lambda = QR = Q[R_{1} \ R_{2}]\qquad\text{with}\qquad R_1\in\bbC^{tk\times tk}\text{ and } R_2\in\bbC^{tk\times (n-tk)}.
\]
Define \[T = (R_{1}(1:{k},1:{k}))^{-1}[R_1(1:{k},{k+1:kt}) \ R_{2}(1:{k},:)]\in \mathbb{C}^{{k\times (n-k)}},\] and $V = [I\ T]\Lambda^* \in  \mathbb{C}^{{k}\times n}$. Let $q$ be the index set with $|q| = k$ such that \[K(s,q)=QR_1(1:{k},1:{k}),\] then $q$ and $V$ will satisfy
\begin{equation}
    K(s,:) \approx K(s,q)V
\end{equation} 
with an approximation error by the QR truncation. By the approximation power of Lagrange interpolation with Mock-Chebyshev points if $K$ is the discretization of a smooth function, we have 
\begin{equation}
    K \approx K(:,q)V
\end{equation} 
with an approximation error coming from the QR truncation and the Lagrange interpolation. Hence, $K(:,q)$ are important columns of $K$ such that they can be ``interpolated'' back to $K$ via a {\it column interpolation matrix} $V$. In this sense, $q$ is called the {\it skeleton} index set, and the rest of indices are called {\it redundant} indices. This column ID requires only $\O{nk^2}$ operations and $\O{nk}$ memories and is denoted as $\cID$ for short.

Similarly, a row ID with $\O{mk^2}$ operations and $\O{mk}$ memories, denoted as $\rID$, can be constructed via
\begin{equation}
    K \approx \Lambda [I \ T]^* K(q,:):=UK(q,:)
\end{equation}
with a {\it row interpolation matrix} $U$.

\subsection{Leaf-root complementary skeletonization (LRCS)}
\label{sec:lrskeleton}

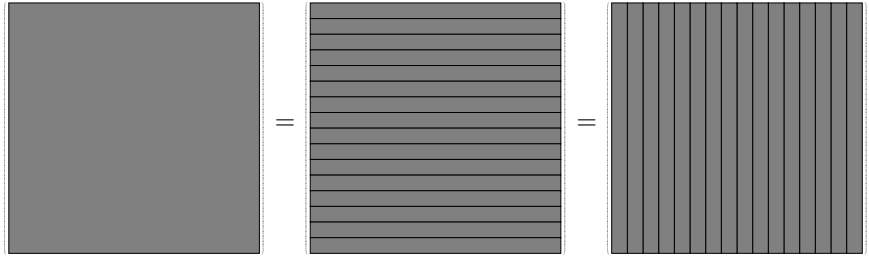
\begin{figure}[!ht]
\begin{minipage}{\textwidth}
\centering
\resizebox{3.5cm}{!}{
\begin{tikzpicture}[baseline=-0.5ex]
      \tikzset{every left delimiter/.style={xshift=-1ex},every right delimiter/.style={xshift=1ex}}
      \matrix (mat) [left delimiter=(, right delimiter=)] {
\draw[fill=gray] (32,0) rectangle (0,32);
\\
      };
\end{tikzpicture}
}
=
\resizebox{3.5cm}{!}{
\begin{tikzpicture}[baseline=-0.5ex]
      \tikzset{every left delimiter/.style={xshift=-1ex},every right delimiter/.style={xshift=1ex}}
      \matrix (mat) [left delimiter=(, right delimiter=)] {
\draw;
\foreach \i in {0,2,...,30}{
   \draw[fill=gray] (0,32-\i) rectangle (32,32-\i-2);
}\\
      };
\end{tikzpicture}
}
=
\resizebox{3.5cm}{!}{
\begin{tikzpicture}[baseline=-0.5ex]
      \tikzset{every left delimiter/.style={xshift=-1ex},every right delimiter/.style={xshift=1ex}}
      \matrix (mat) [left delimiter=(, right delimiter=)] {
\draw;
\foreach \i in {0,2,...,30}{
   \draw[fill=gray] (32-\i,0) rectangle (32-\i-2,32);
}\\
      };
\end{tikzpicture}
}\\
\end{minipage}
\caption{The left figure is a complementary two-dimensional low-rank kernel matrix $K$. Assume that the depth of the quadtrees of column and row spaces is $3$. The middle figure illustrates the root-leaf partitioning that divides the row index set into $16$ subsets as $16$ leaves. The right one is for the leaf-root partitioning that divides the column index set into $16$ subsets as $16$ leaves.}
\label{fig:partition-K}
\end{figure}

This subsection introduces the LRCS of a 2D complementary low-rank kernel matrix $K$, $K\approx USV$, via $\cID$s of the submatrices corresponding to the leaf-root levels of the column-row quadtrees (e.g., see the associated matrix partition in Figure~\ref{fig:partition-K} (right)), and $\rID$s of the submatrices corresponding to the root-leaf levels of the column-row quadtrees (e.g., see the associated matrix partition in Figure~\ref{fig:partition-K} (middle)). Assume 
$k_\epsilon$ is constant in all IDs for low-rank approximations and is denoted by $k$ for simplicity.

Assume that the row index set $r$ and the column index set $c$ of $K$ can be partitioned into leaves $\{r_i\}_{i \in \calI^L}$ and $\{c_j\}_{j \in \calI^L}$ at the leaf level of the row and column quadtrees as follows
\begin{equation}
    \label{eq:partrc}
    r = [r_{0},r_{1},\cdots,r_{m-1}] \qquad (\text{and } c = [c_{0},c_{1},\cdots,c_{m-1}]),
\end{equation}
with $|r_i|=n_0$ (and $|c_j|=n_0$) for all $0\leq i,j \leq m-1$, where $m = 4^L = \frac{N}{n_0}, L = \log_4 \left(N\right)- \log_4 \left(n_0\right)$, and $L+1$ is the depth of quadtrees $T_{X}$ and $T_{\Omega}$. See an example of row and column quadtrees with $m=16$ in Figure~\ref{fig:partition-K}.

Apply $\rID$ to each $K(r_{i},:)$ to obtain the row interpolation matrix $U_{i}$ and the associated skeleton indices $\hat{r}_{i}\subset r_{i}$ for all $0 \le i \le m-1$. Then, after denoting $K(\hat{r},:)$ as the important skeleton of $K$, where
\begin{equation}
    \label{eq:partrc2}
    \hat{r} = [\hat{r}_{0},\hat{r}_{1},\cdots,\hat{r}_{m-1}],
\end{equation}
we have
\[
K\approx
\begin{pmatrix}
    U_{1} & & & \\
    & U_{2} & & \\
    & & \ddots & \\
    & & & U_{m} 
\end{pmatrix}
\begin{pmatrix}
    K(\hat{r}_{0},c_{0}) & K(\hat{r}_{0},c_{1}) & \hdots & K(\hat{r}_{0},c_{m-1}) \\
    K(\hat{r}_{1},c_{0}) & K(\hat{r}_{1},c_{1}) & \hdots & K(\hat{r}_{1},c_{m-1}) \\
    \vdots & \vdots & \ddots & \vdots \\
    K(\hat{r}_{m-1},c_{0}) & K(\hat{r}_{m-1},c_{1}) & \hdots & K(\hat{r}_{m-1},c_{m-1})  
\end{pmatrix}
:=UM.
\]

Similarly, apply $\cID$ to each $K(\hat{r},c_{j})$ to obtain the column interpolation matrix $V_{j}$ and the skeleton indices $\hat{c}_{j}\subset c_{j}$ for all $0 \le j \le m-1$. Then, the LRCS of $K$ will be formed as
\begin{equation}
    \label{facK}
    \begin{split}
        K\approx & \;
        \begin{pmatrix}
            U_{1} & & & \\
            & U_{2} & & \\
            & & \ddots & \\
            & & & U_{m} 
        \end{pmatrix}
        \begin{pmatrix}
            K(\hat{r}_{0},\hat{c}_{0}) & K(\hat{r}_{0},\hat{c}_{1}) & \hdots & K(\hat{r}_{0},\hat{c}_{m-1}) \\
            K(\hat{r}_{1},\hat{c}_{0}) & K(\hat{r}_{1},\hat{c}_{1}) & \hdots & K(\hat{r}_{1},\hat{c}_{m-1}) \\
            \vdots & \vdots & \ddots & \vdots \\
            K(\hat{r}_{m-1},\hat{c}_{0}) & K(\hat{r}_{m-1},\hat{c}_{1}) & \hdots & K(\hat{r}_{m-1},\hat{c}_{m-1})  
        \end{pmatrix}
        \begin{pmatrix}
            V_{1} & & & \\
            & V_{2} & & \\
            & & \ddots & \\
            & & & V_{m} 
        \end{pmatrix} \\
        := & \; USV.
    \end{split}
\end{equation}

\begin{figure}[!ht]
\begin{minipage}{\textwidth}
\centering
\resizebox{3.5cm}{!}{
\begin{tikzpicture}[baseline=-0.5ex]
      \tikzset{every left delimiter/.style={xshift=-1ex},every right delimiter/.style={xshift=1ex}}
      \matrix (mat) [left delimiter=(, right delimiter=)] {
\draw[fill=gray] (32,0) rectangle (0,32);
\\
      };
\end{tikzpicture}
}
$\approx$
\resizebox{1.85cm}{!}{
\begin{tikzpicture}[baseline=-0.5ex]
      \tikzset{every left delimiter/.style={xshift=-1ex},every right delimiter/.style={xshift=1ex}}
      \matrix (mat) [left delimiter=(, right delimiter=)] {
      \draw;
\foreach \i in {0,1,...,15}{
   \draw[fill=gray]  (\i,32-\i-\i) rectangle (\i+1,32-\i-\i-2);
}\\
      };
\end{tikzpicture}
}
\resizebox{1.85cm}{!}{
\begin{tikzpicture}[baseline=-0.5ex]
      \tikzset{every left delimiter/.style={xshift=-1ex},every right delimiter/.style={xshift=1ex}}
      \matrix (mat) [left delimiter=(, right delimiter=)] {
      \draw;
\foreach \i in {0,1,...,15}{
\foreach \j in {0,1,...,15}{
      \draw[fill=gray]  (\i,16-\j) rectangle (\i+1,16-\j-1);
}
}\\
      };
\end{tikzpicture}
}
\resizebox{3.5cm}{!}{
\begin{tikzpicture}[baseline=-0.5ex]
      \tikzset{every left delimiter/.style={xshift=-1ex},every right delimiter/.style={xshift=1ex}}
      \matrix (mat) [left delimiter=(, right delimiter=)] {
      \draw;
\foreach \i in {0,1,...,15}{
   \draw[fill=gray]  (32-\i-\i,\i) rectangle (32-\i-\i-2,\i+1);
}\\
      };
\end{tikzpicture}
}\\
\end{minipage}
\caption{An example of the LRCS in \eqref{facK} of the complementary two-dimensional low-rank kernel matrix $K$ in Figure \ref{fig:partition-K}. Non-zero submatrices in \eqref{facK} are shown in gray areas. }
\label{fig:rlfac}
\end{figure}
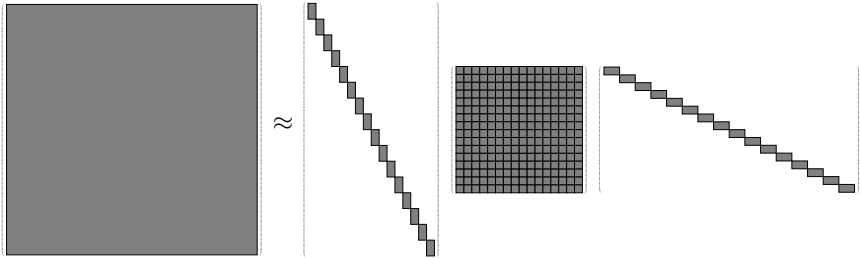

For a concrete example, Figure~\ref{fig:rlfac} illustrates the non-zero pattern of the LRCS  in \eqref{facK} of $K$ in Figure~\ref{fig:partition-K}.

The main contribution of the LRCS is that $M$ and $S$ are only required to be generated and stored via the skeleton of row and column index sets with $\O{\frac{k^3}{n_0}N}$ operations and $\O{\frac{k^2}{n_0}N}$ memories, instead of being computed explicitly, since there are only $2m=\frac{2N}{n_0}$ IDs in total. Notice that the matrix $S$ in $K\approx USV$ is also a complementary low-rank matrix. The row and column quadtrees $\hat{T}_X$ and $\hat{T}_\Omega$ of $S$ are the compressed version of the row and column quadtrees ${T}_X$ and ${T}_\Omega$ of $K$. If we consider $\hat{T}_X$ and $\hat{T}_\Omega$ as quadtrees with one depth less than the leaf level of ${T}_X$ and ${T}_\Omega$, they will be compressible.

\subsection{Matrix splitting with complementary skeletonization (MSCS)}
\label{sec:matrixsplit}

Now we introduce another key idea repeatedly applied in 2D IDBF, the MSCS. According to the nodes of the second level of the row and column quadtrees $T_X$ and $T_\Omega$ (with $m=4^L$ leaves), the complementary 2D low-rank kernel matrix $K$ can be split into a $4\times 4$ block matrix
\begin{equation}
    \label{eq:matrixsplit}
    K = 
    \begin{pmatrix}
        K_{11} & K_{12} & K_{13} & K_{14} \\
        K_{21} & K_{22} & K_{23} & K_{24} \\
        K_{31} & K_{32} & K_{33} & K_{34} \\
        K_{41} & K_{42} & K_{43} & K_{44}
    \end{pmatrix}.
\end{equation}
It is obvious that $K_{ij}$ is complementary low-rank for all $1 \le i,j \le 4$, with row and column quadtrees $T_{X,ij}$ and $T_{\Omega,ij}$ of depth $L-1$ and with $\frac{m}{4}$ leaves.

Suppose that the LRCS of each $K_{ij}$ is $K_{ij} \approx U_{ij}S_{ij}V_{ij}$. Then, according to the LRCS of $K_{ij}$, the matrix splitting with complementary skeletonization (MSCS) of the kernel matrix $K$ can be proposed as:
\begin{equation}
    \label{eq:expressK}
    K \approx USV,
\end{equation}
where
\begin{equation}
    \label{eq:expressU}
    U = 
    \begin{pmatrix}
        U_{1} & U_{2} & U_{3} & U_{4} 
    \end{pmatrix}
    \quad\text{with}\quad
    U_k = 
    \begin{pmatrix}
        U_{1k} & & &  \\
        & U_{2k} & &  \\
        & & U_{3k} & \\
        & & & U_{4k} 
    \end{pmatrix},
\end{equation}
\begin{equation}
    \label{eq:expressS}
    S = 
  \begin{pmatrix}
        \bar{S}_{11} & \bar{S}_{12} & \bar{S}_{13} & \bar{S}_{14} \\
        \bar{S}_{21} & \bar{S}_{22} & \bar{S}_{23} & \bar{S}_{24} \\
        \bar{S}_{31} & \bar{S}_{32} & \bar{S}_{33} & \bar{S}_{34} \\
        \bar{S}_{41} & \bar{S}_{42} & \bar{S}_{43} & \bar{S}_{44}
    \end{pmatrix}
    \quad\text{with}\quad 
    \bar{S}_{ij} \quad \text{as a $4$ by $4$ block matrix with the $(j,i)$-th block as $S_{ji}$,}
\end{equation}
\begin{equation}
    \label{eq:expressV}
    V = 
    \begin{pmatrix}
        V_{1}  \\
        V_{2}  \\
        V_{3} \\
        V_{4} 
    \end{pmatrix}
    \quad \text{with}\quad
    V_k = 
    \begin{pmatrix}
        V_{k1} & & & \\
        & V_{k2} & & \\
        & & V_{k3} & \\
        & & & V_{k4} \\
    \end{pmatrix}.
\end{equation}

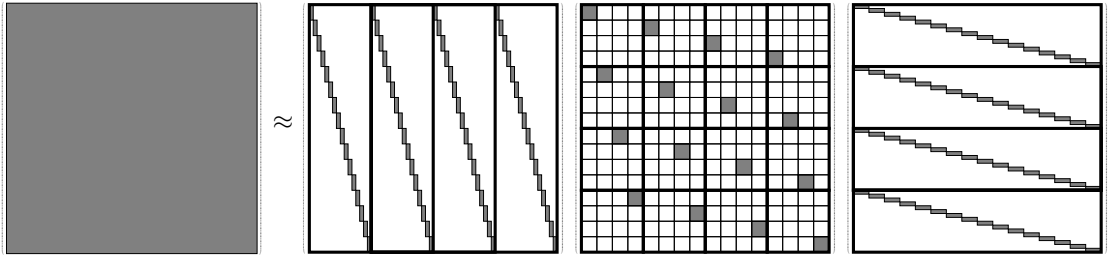
\begin{figure}[!ht]
\begin{minipage}{\textwidth}
\centering
\resizebox{3.5cm}{!}{
\begin{tikzpicture}[baseline=-0.5ex]
      \tikzset{every left delimiter/.style={xshift=-1ex},every right delimiter/.style={xshift=1ex}}
      \matrix (mat) [left delimiter=(, right delimiter=)] {
\draw[fill=gray] (32,0) rectangle (0,32);
\\
      };
\end{tikzpicture}
}
$\approx$
\resizebox{3.5cm}{!}{
\begin{tikzpicture}[baseline=-0.5ex]
      \tikzset{every left delimiter/.style={xshift=-1ex},every right delimiter/.style={xshift=1ex}}
      \matrix (mat) [left delimiter=(, right delimiter=)] {
      \draw;
\foreach \j in {0,1,...,15}{
    \draw[fill=gray] (\j/2,32-\j-\j) rectangle (\j/2+0.5,32-\j-\j-2);
    \draw[fill=gray] (8+\j/2,32-\j-\j) rectangle (8+\j/2+0.5,32-\j-\j-2);
    \draw[fill=gray] (16+\j/2,32-\j-\j) rectangle (16+\j/2+0.5,32-\j-\j-2);
    \draw[fill=gray] (24+\j/2,32-\j-\j) rectangle (24+\j/2+0.5,32-\j-\j-2);
}
      \foreach \i in {0,1,...,4}{
      \draw  [line width=4mm ] (8*\i,0) rectangle (8,32);
}
\\
      };
\end{tikzpicture}
}
\resizebox{3.5cm}{!}{
\begin{tikzpicture}[baseline=-0.5ex]
      \tikzset{every left delimiter/.style={xshift=-1ex},every right delimiter/.style={xshift=1ex}}
      \matrix (mat) [left delimiter=(, right delimiter=)] {
      \draw;
      \draw[fill=gray] (0,30) rectangle (2,32);
      \draw[fill=gray] (8,28) rectangle (10,30);
      \draw[fill=gray] (16,26) rectangle (18,28);
      \draw[fill=gray] (24,24) rectangle (26,26);
      \draw[fill=gray] (2,22) rectangle (4,24);
      \draw[fill=gray] (10,20) rectangle (12,22);
      \draw[fill=gray] (18,18) rectangle (20,20);
      \draw[fill=gray] (26,16) rectangle (28,18);
      \draw[fill=gray] (4,14) rectangle (6,16);
      \draw[fill=gray] (12,12) rectangle (14,14);
      \draw[fill=gray] (20,10) rectangle (22,12);
      \draw[fill=gray] (28,8) rectangle (30,10);
      \draw[fill=gray] (6,6) rectangle (8,8);
      \draw[fill=gray] (14,4) rectangle (16,6);
      \draw[fill=gray] (22,2) rectangle (24,4);
      \draw[fill=gray] (30,0) rectangle (32,2);
      \foreach \i in {0,1,...,15}{
\foreach \j in {0,1,...,15}{
      \draw  [line width=1mm ] (2*\i,32-2*\j) rectangle (2*\i+2,32-2*\j-2);
}
}
      \foreach \i in {0,1,...,3}{
\foreach \j in {0,1,...,3}{
      \draw  [line width=4mm ] (8*\i,32-8*\j) rectangle (8*\i+8,32-8*\j-8);
}
}
\\
      };
\end{tikzpicture}
}
\resizebox{3.5cm}{!}{
\begin{tikzpicture}[baseline=-0.5ex]
      \tikzset{every left delimiter/.style={xshift=-1ex},every right delimiter/.style={xshift=1ex}}
      \matrix (mat) [left delimiter=(, right delimiter=)] {
      \draw;
\foreach \j in {0,1,...,15}{
    \draw[fill=gray] (\j+\j,32-\j/2) rectangle (\j+\j+2,32-\j/2-0.5);
    \draw[fill=gray] (\j+\j,32-\j/2-8) rectangle (\j+\j+2,32-\j/2-8-0.5);
    \draw[fill=gray] (\j+\j,32-\j/2-16) rectangle (\j+\j+2,32-\j/2-16-0.5);
    \draw[fill=gray] (\j+\j,32-\j/2-24) rectangle (\j+\j+2,32-\j/2-24-0.5);
}
      \foreach \i in {0,1,...,4}{
      \draw  [line width=4mm ] (0,8*\i) rectangle (32,8);
}
\\
      };
\end{tikzpicture}
}\\
\end{minipage}
\caption{The illustration of an MSCS of a complementary 2D low-rank kernel matrix $K \approx USV$ with quadtrees of depth $3$ and $16$ leaf nodes in Figure \ref{fig:partition-K}. Non-zero blocks in \eqref{eq:expressU}-\eqref{eq:expressV} are shown in gray areas. $\{U_i\}_{1\leq i\leq 4}$, $\{\bar{S}_{ij}\}_{1\leq i\leq j\leq 4}$, and $\{V_i\}_{1\leq i\leq 4}$ are visualized by large submatrices with wide edges in the middle left, middle right, and right figures, respectively. }
\label{fig:exlfac}
\end{figure}

Recall that the middle factor $S$ is only required to be generated by some entries of the original kernel matrix, forming \eqref{eq:expressK}-\eqref{eq:expressV} will be a linear scaling algorithm as well. Figure~\ref{fig:exlfac} illustrates the MSCS of a complementary 2D low-rank kernel matrix $K$ with quadtrees of depth $3$ and $16$ leaf nodes in Figure~\ref{fig:partition-K}.

\subsection{Recursive MSCS}
\label{sec:recursivefac}
This subsection applies MSCS recursively to obtain the full 2D IDBF of a complementary 2D low-rank kernel matrix $K$.

First, we denote the first level of MSCS of $K$ in \eqref{eq:expressK} as
\begin{equation}
    \label{eq:leaffac}
    {K \approx U^{L}S^{L}V^{L}},
\end{equation}
where $U^{L}, S^{L}, V^{L}$ maintain the same structures as \eqref{eq:expressU}-\eqref{eq:expressV}. Then, the index set $r$ and the column index set $c$ of $K$ can be partitioned into leaves $\{r_i\}_{0\leq i\leq m-1}$ and $\{c_j\}_{0\leq j\leq m-1}$ at the leaf level of the row and column quadtrees as \eqref{eq:partrc}. In addition, the skeleton index sets $\hat{r}_i\subset r_i$ and $\hat{c}_j\subset c_j$ will be obtained by applying the $\rID$s and $\cID$s to the construction of \eqref{eq:leaffac}, and the middle factor $S^{L}$ will be constructed by the non-zero submatrices $S^{L}_{ij}$ for all $1 \le i,j \le 4$ as follows:
\begin{equation}
    S^{L}_{ij} = 
    \begin{pmatrix}
        K(\hat{r}_{(i-1)(m-1)/4+1},\hat{c}_{(j-1)(m-1)/4+1}) & \cdots & K(\hat{r}_{(i-1)(m-1)/4+1},\hat{c}_{j(m-1)/4}) \\
        \vdots & \ddots & \vdots \\
        K(\hat{r}_{i(m-1)/4},\hat{c}_{(j-1)(m-1)/4+1}) & \cdots & K(\hat{r}_{i(m-1)/4},\hat{c}_{j(m-1)/4}) \\
    \end{pmatrix}.
\end{equation}

Since $S^L_{ij}$ consists of the important rows and columns of $K_{ij}$ for all $1 \le i,j \le 4$, it will inherit the complementary low-rank property of $K_{ij}$. Suppose that $T_{X,ij}$ and $T_{\Omega,ij}$ are the quadtrees of the row and column spaces of $K_{ij}$ with $\frac{m}{4}$ leaves and $L-1$ depth. Then, $S^L_{ij}$ has compressible row and column quadtrees $\hat{T}_{X,ij}$ and $\hat{T}_{\Omega,ij}$ with $\frac{m}{16}$ leaves and $L-2$ depth according to Subsection~\ref{sec:lrskeleton}.

Next, a recursive MSCS will be applied to each $S^L_{ij}$. The first step is similar to that of MSCS, we divide each $S^{L}_{ij}$ into a $4\times 4$ block matrix according to the nodes at the second level of its row and column quadtrees:
\begin{equation}
    \label{eq:sl}
    S^{L}_{ij} = 
    \begin{pmatrix}
        (S^{L}_{ij})_{11} & (S^{L}_{ij})_{12} & (S^{L}_{ij})_{13} & (S^{L}_{ij})_{14} \\
        (S^{L}_{ij})_{21} & (S^{L}_{ij})_{22} & (S^{L}_{ij})_{23} & (S^{L}_{ij})_{24} \\
        (S^{L}_{ij})_{31} & (S^{L}_{ij})_{32} & (S^{L}_{ij})_{33} & (S^{L}_{ij})_{34} \\
        (S^{L}_{ij})_{41} & (S^{L}_{ij})_{42} & (S^{L}_{ij})_{43} & (S^{L}_{ij})_{44} 
    \end{pmatrix}.
\end{equation}
For each block $(S^{L}_{ij})_{k\ell}$, the LRCS can be constructed as $(S^{L}_{ij})_{k\ell} \approx (U^{L-1}_{ij})_{k\ell} (S^{L-1}_{ij})_{k\ell} (V^{L-1}_{ij})_{k\ell}$ for all $1 \le k,\ell \le 4$. After that, the MSCS of $S^{L}_{ij}$ will be obtained as follows:
\begin{equation}
    \label{eq:facSL1234}
    S^{L}_{ij} \approx U^{L-1}_{ij}S^{L-1}_{ij}V^{L-1}_{ij},
\end{equation}
where $U^{L-1}_{ij}, S^{L-1}_{ij}, V^{L-1}_{ij}$ are constructed by $(U^{L-1}_{ij})_{k\ell} (S^{L-1}_{ij})_{k\ell} (V^{L-1}_{ij})_{k\ell}$ for all $1 \le k,\ell \le 4$ as in \eqref{eq:expressU}-\eqref{eq:expressV}.

Eventually, the factorization in \eqref{eq:facSL1234} for all $1 \le i,j \le 4$ will be combined to form a factorization of $S^{L}$:
\begin{equation}
    \label{eq:facSL}
    S^{L} \approx U^{L-1}S^{L-1}V^{L-1},
\end{equation}
where 
\begin{equation}
    \label{eq:expressU_L-1}
    U^{L-1} = 
    \begin{pmatrix}
        U^{L-1}_{1} & & & \\
        & U^{L-1}_{2} & &  \\
        & & U^{L-1}_{3} &  \\
        & & & U^{L-1}_{4}  \\
    \end{pmatrix}
    \quad\text{with}\quad
        U^{L-1}_k = 
    \begin{pmatrix}
        U^{L-1}_{1k} & & & \\
        & U^{L-1}_{2k} & &  \\
        & & U^{L-1}_{3k} & \\
        & & & U^{L-1}_{4k}  \\
    \end{pmatrix},
\end{equation}
\begin{equation}
    \label{eq:expressS_L-1}
    S^{L-1} = 
  \begin{pmatrix}
        \bar{S}^{L-1}_{11} & \bar{S}^{L-1}_{12} & \bar{S}^{L-1}_{13} & \bar{S}^{L-1}_{14} \\
        \bar{S}^{L-1}_{21} & \bar{S}^{L-1}_{22} & \bar{S}^{L-1}_{23} & \bar{S}^{L-1}_{24} \\
        \bar{S}^{L-1}_{31} & \bar{S}^{L-1}_{32} & \bar{S}^{L-1}_{33} & \bar{S}^{L-1}_{34} \\
        \bar{S}^{L-1}_{41} & \bar{S}^{L-1}_{42} & \bar{S}^{L-1}_{43} & \bar{S}^{L-1}_{44}
    \end{pmatrix}
\end{equation}
with $\bar{S}^{L-1}_{ij}$ as a $4\times 4$ block matrix with the $(j,i)$-th block as $S_{ji}^{L-1}$,
\begin{equation}
    \label{eq:expressV_L-1}
    V^{L-1} = 
    \begin{pmatrix}
        V^{L-1}_{1} & & &  \\
        & V^{L-1}_{2} & & \\
        & & V^{L-1}_{3} &  \\
        & & & V^{L-1}_{4} \\
    \end{pmatrix}
\quad \text{with}\quad
    V_k^{L-1} = 
    \begin{pmatrix}
        V^{L-1}_{k1} & & & \\
        & V^{L-1}_{k2} & & \\
        & & V^{L-1}_{k3} & \\
        & & & V^{L-1}_{k4} \\
    \end{pmatrix}.
\end{equation}
Hence, the second level factorization of $K$ can be constructed as follows:
\[
K\approx U^L U^{L-1}S^{L-1}V^{L-1}V^L.
\]

Comparing \eqref{eq:leaffac} and \eqref{eq:facSL}, a fractal structure can be found in each level of the middle factor $S^L$ and $S^{L-1}$. For example, $S^L$ and $S^{L-1}$ have the same structure consisting of 16 submatrices as shown in \eqref{eq:expressS} and \eqref{eq:expressS_L-1}. Besides, submatrices $S^{L-1}_{ij}$ can be factorized into a product of three matrices $U^{L-2}_{ij}$, $S^{L-2}_{ij}$, $V^{L-2}_{ij}$ with the same sparsity structure as that of $S^L$ in \eqref{eq:facSL}-\eqref{eq:expressV_L-1}. Thus, the recursive MSCS can be applied repeatedly to each $S^{\ell}$ for $\ell=L$, $L-1$, $\dots$, $\frac{L}{2}$ and the matrix factors can be assembled hierarchically as follows:
\begin{equation}
    \label{eq:finishIDBF}
    K \approx U^{L}U^{L-1}\cdots U^{h} S^h V^{h}\cdots V^{L-1}V^{L},
\end{equation}
where $h=\frac{L}{2}$.

In the $\ell$-th  recursive MSCS, there are $4^{2(L-\ell+1)}$ dense submatrices with compressible row and column quadtrees, which consist $\frac{m}{4^{2(L-\ell+1)}}$ leaves and depth $L-2(L-\ell+1)$, in $S^\ell$. Thus, after $h=\frac{L}{2}$ iterations, the recursive MSCS will stop, since there is not any compressible submatrix in $S^h$. Otherwise, when $S^\ell$ is still compressible, there are $4^{2(L-\ell+1)}\frac{m}{4^{2(L-\ell+1)}}=\frac{N}{n_0}$ low-rank submatrices to be factorized. Linear IDs only require $\O{k^3}$ operations for each low-rank submatrix, and hence at most $\O{\frac{k^3}{n_0} N}$ for each level of factorization, and $\O{\frac{k^3}{n_0} N\L{N}}$ for the whole 2D IDBF.

\subsection{Extensions}\label{sec:ext}

We have introduced the 2D IDBF for a complementary low-rank kernel matrix $K$ in the entire domain $X \times \Omega$. Although we have assumed the uniform grid in $X$ and $\Omega$, the butterfly factorization extends naturally to more general settings. In the case with non-uniform point sets $X$ or $\Omega$, one can still construct a butterfly factorization for $K$ following the same procedure. More specifically, we construct two trees $T_X$ and $T_\Omega$ adaptively via hierarchically partitioning the square domains covering $X$ and $\Omega$.  For non-uniform point sets $X$ and $\Omega$, the numbers of points in $A^\ell_i$ and $B^\ell_j$ are different. If a node does not contain any point inside it, it is simply discarded from the quadtree. We can also extend the 2D IDBF to the 3D case by constructing two octrees $T_X$ and $T_\Omega$ via hierarchically partitioning the cube domains covering $X$ and $\Omega$. Lastly, the numerical rank in all low-rank approximations in the IDBF presented is fixed. It's easy to extend the current version to an adaptive one with an adaptive rank $k_\epsilon$ in IDs depending on a target accuracy $\epsilon$. For example, choose $k_\epsilon = \min \{k: R_1(k,k) \le \epsilon R_1(1,1) \}$ and update $k\leftarrow k_\epsilon$ after the QR in IDs. An adaptive rank leads to a more compressed IDBF while a fixed rank results in a more predictable sparsity pattern in IDBF.

\section{Numerical results} \label{sec:results}

This section presents several numerical examples to demonstrate the efficiency of the proposed framework. All implementations are in MATLAB\textsuperscript{\textregistered} on a server computer with a single thread and 3.2 GHz CPU, and are available in the ButterflyLab (\url{https://github.com/ButterflyLab/ButterflyLab}). 

Let $\{g^d(x),x\in X\}$ and $\{g^b(x),x\in X\}$ denote the results given by the direct matrix-vector multiplication and MIDBF, respectively. The accuracy of applying fast algorithms is estimated by the relative error defined as follows:
\begin{equation}
    \epsilon^b = \sqrt{\cfrac{\sum_{x\in S}|g^b(x)-g^d(x)|^2} {\sum_{x\in S}|g^d(x)|^2}},
\end{equation}
where $S$ is an index set containing $256$ randomly sampled row indices of the kernel matrix $K$. The error for recovering the kernel function is defined as
\begin{equation}
    \epsilon^K = \frac{\|e^{2\pi\i \Phi(S,S)} - e^{2\pi\i U(S,:)V(:,S)^T}\|_2} {\|e^{2\pi\i \Phi(S,S)}\|_2},
\end{equation}
where $\Phi$ is the phase matrix and $UV^T$ is its low-rank recovery. In all of our examples, the tolerance parameter $\epsilon$ is set to $10^{-9}$, the over-sampling parameter $q$ in low-rank phase matrix factorization is set to $2$, the threshold $\tau$ for detecting discontinuity in multidimensional cases is set to $\frac{1}{4}$, the number of points in a leave node $n_0$ in the MIDBF is set to $8^d$, and the over-sampling parameter $t$ in ID in MIDBF is set to $5$. We apply IDs with an adaptive rank and $k$ denotes our empirically estimated rank.

\subsection{Accuracy and scaling of low-rank matrix recovery and MIDBF}

In this part, we present numerical results of several examples to demonstrate the accuracy and asymptotic scaling of the proposed low-rank matrix recovery for phase functions, and MIDBF. With no loss of generality, we only focus on Scenario $2$ of indirect access. {Since there is not any detected discontinuous point in the phase matrices of Example~$1$ and Example~$3$ when $\frac{1}{4} \ge \tau \ge \frac{1}{10}$, we will only address the related discontinuities discussion in Example~$2$. Each experiment will be repeatedly tested for $10$ times.}

\paragraph{Example 1.}

Our first example is to evaluate a 2D generalized Radon transform which is a Fourier integral operator (FIO) \cite{NUFFTorBF} defined as follows:
\begin{equation}
    \label{eq:2D}
    g(x) = \int_{\mathbb{R}}e^{2\pi i \Phi(x,\xi)} \widehat{f}(\xi) d\xi,
\end{equation}
where $\widehat{f}$ is the {Fourier} transform of $f$, and $\Phi(x,\xi)$ is a phase function given by
\begin{equation}
    \begin{split}
        \Phi(x,\xi) & = x \cdot \xi' + \sqrt{c_1^2(x) \cdot \xi_1^2 + c_2^2(x) \cdot \xi_2^2}, \\
        c_1(x) & = (2 + \sin(2\pi x_1) \sin(2\pi x_2))/16, \\
        \text{and}\quad
        c_2(x) & = (2 + \cos(2\pi x_1) \cos(2\pi x_2))/16.
    \end{split}
\end{equation}
The discretization of \eqref{eq:2D} is
\begin{equation}
    \label{eqn:2D}
    g(x) = \sum_{\xi \in \Omega} e^{2\pi i \Phi(x,\xi)} \widehat{f} (\xi),
\quad x \in X,
\end{equation}
where $X$ and $\Omega$ are the sets of $\O{N}$ points uniformly distributed in $[0,1) \times [0,1)$. The computation in \eqref{eqn:2D} approximately integrates over spatially varying ellipses, for which $c_1(x)$ and $c_2(x)$ are the axis lengths of the ellipse centered at the point $x\in X$. The corresponding matrix form of \eqref{eqn:2D} is simply
\begin{equation}
    u=Kg, \quad K = (e^{2\pi i\Phi(x,\xi)})_{x\in X,\xi\in\Omega}.
    \label{eq:uKg}
\end{equation}

\begin{figure}[!ht]
    \centering
    \begin{tabular}{ccccc}
        \includegraphics[height=1.3in,trim={3.7cm 7.2cm 3.5cm 7.5cm}, clip]{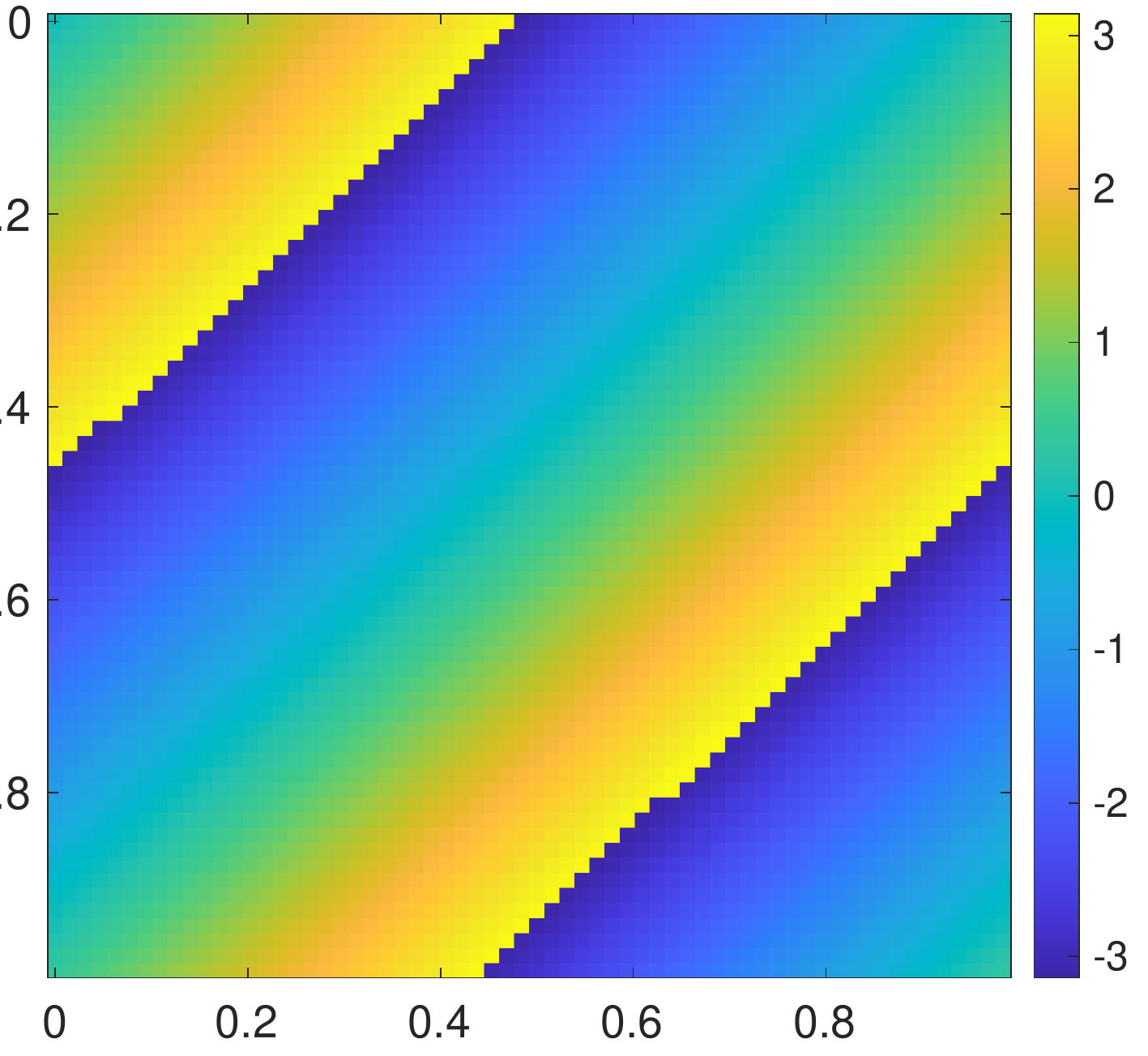} & \includegraphics[height=1.3in,trim={3.7cm 7.2cm 3.5cm 7.5cm}, clip]{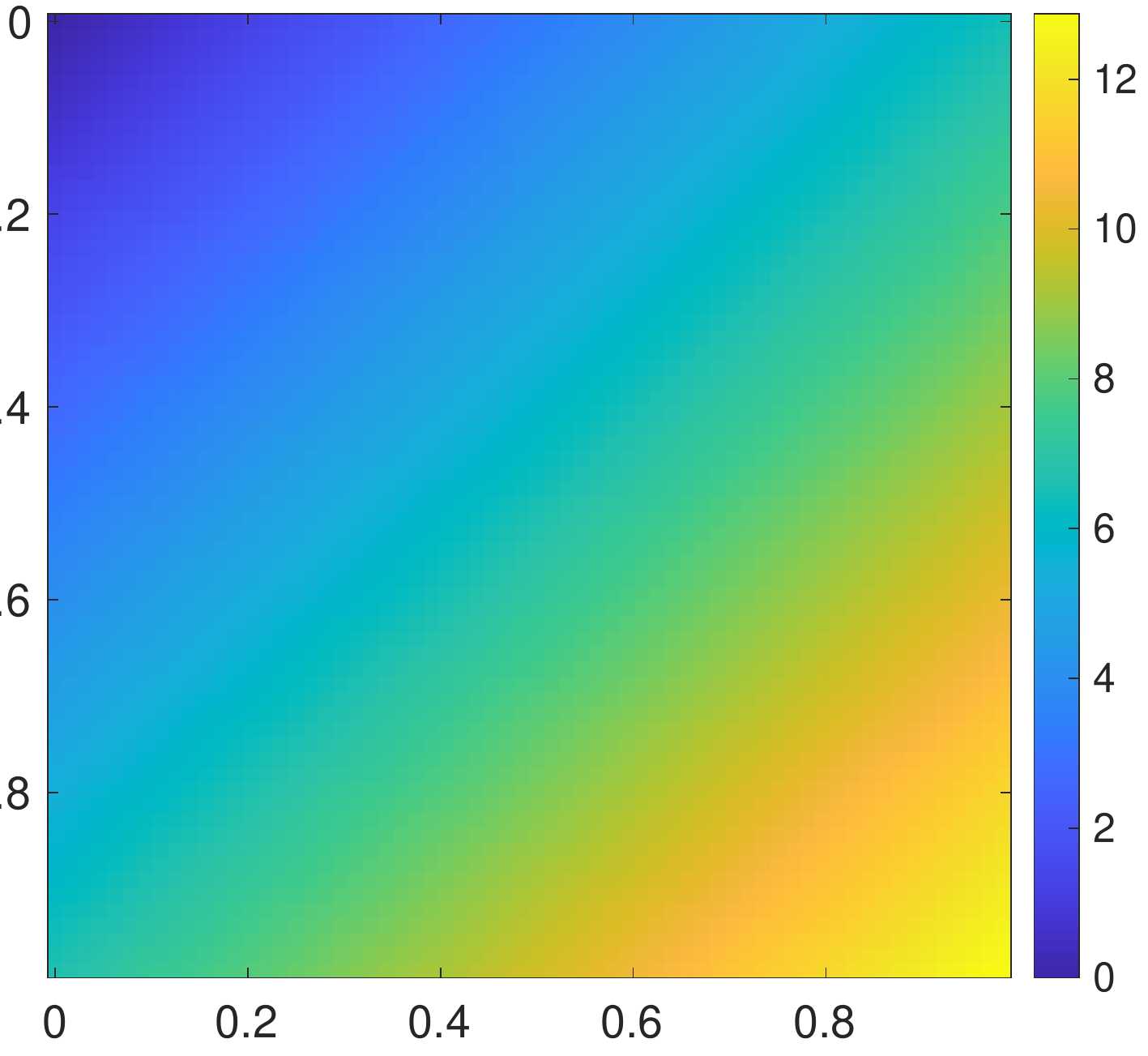} & & \includegraphics[height=1.3in,trim={3.7cm 7.2cm 3.5cm 7.5cm}, clip]{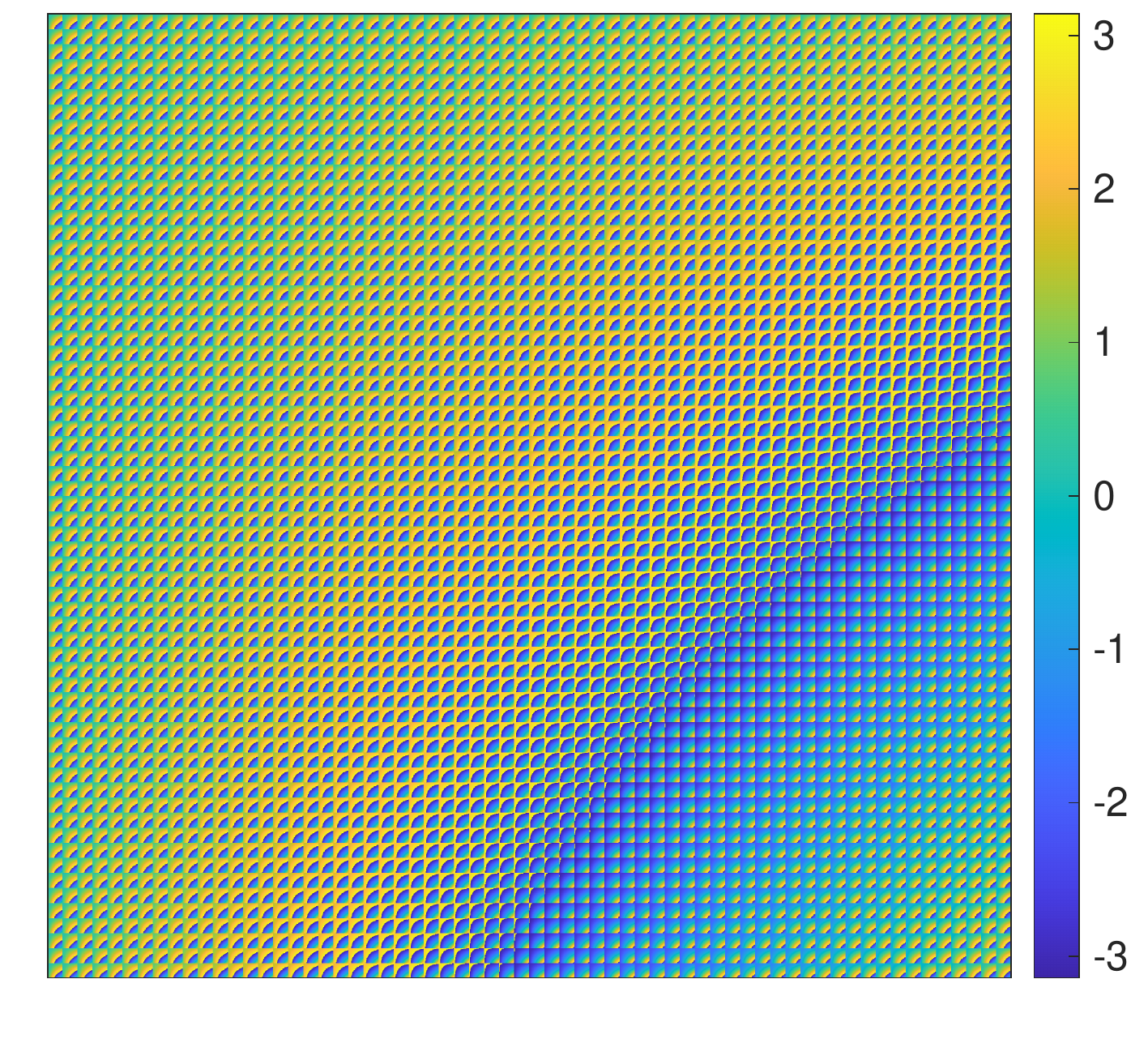} & \includegraphics[height=1.3in,trim={3.7cm 7.2cm 3.5cm 7.5cm}, clip]{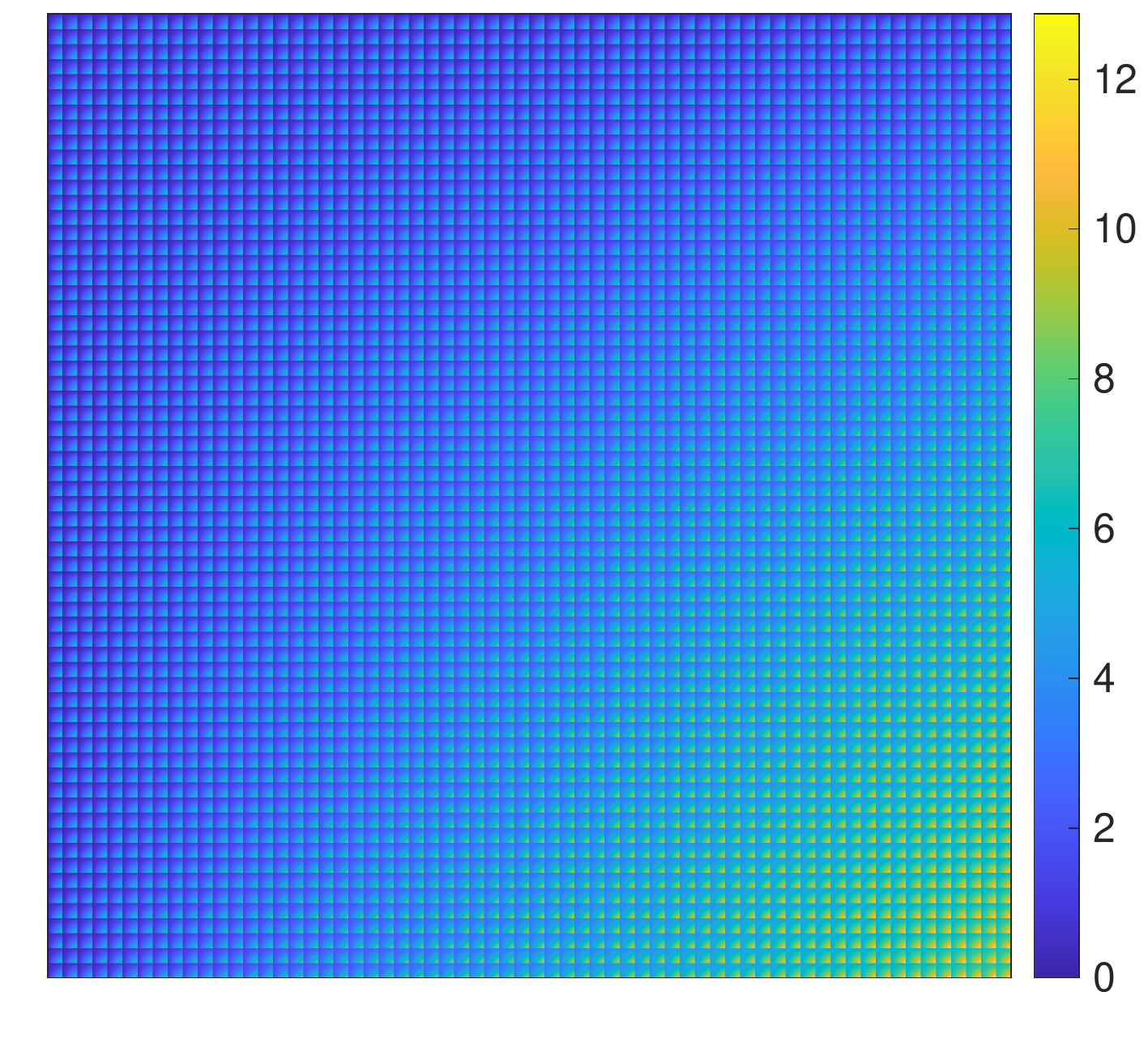} \\
        (a) & (b) & & (c) & (d)
    \end{tabular}
    \caption{{Phase recovery results for the 2D uniform FIO given in \eqref{eqn:2D}. $N = 64^2$ is the size of the phase matrix $\left(\Phi(x,\xi)\right)_{x\in X,\xi\in\Omega}$. (a) A row vector of the phase matrix before recovery and reshaped into a matrix of size $64\times 64$. (b) A recovered row vector of the phase matrix when it is reshaped into a matrix of size $64\times 64$. (c) The phase matrix of size $64^2\times 64^2$ before recovery.  (d) The recovered phase matrix of size $64^2\times 64^2$.}}
    \label{fig:disc}
\end{figure}

The framework is applied to recover the phase functions in the form of low-rank matrix factorization, compute the MIDBF of the kernel function, and apply it to a randomly generated $f$ in \eqref{eq:2D} to obtain $g$. {Figure~\ref{fig:disc} illustrates the results of the recovery step for the phase matrix $\left(\Phi(x,\xi)\right)_{x\in X,\xi\in\Omega}$, the recovered phase matrix in (d) is set as an initial guess for the low-rank factorization step.}

\begin{figure}[!ht]
    \centering
    \begin{tabular}{cccc}
        \includegraphics[height=1.6in,trim={13.5cm 8.5cm 2.2cm 8.5cm}, clip]{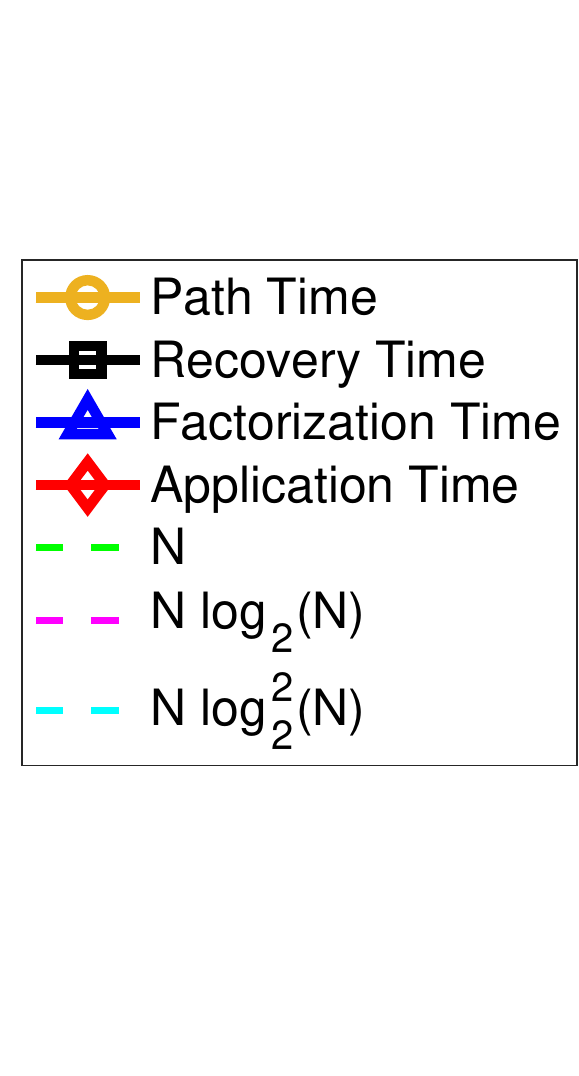} &
        \includegraphics[height=1.6in,trim={3cm 6.3cm 3.5cm 7cm}, clip]{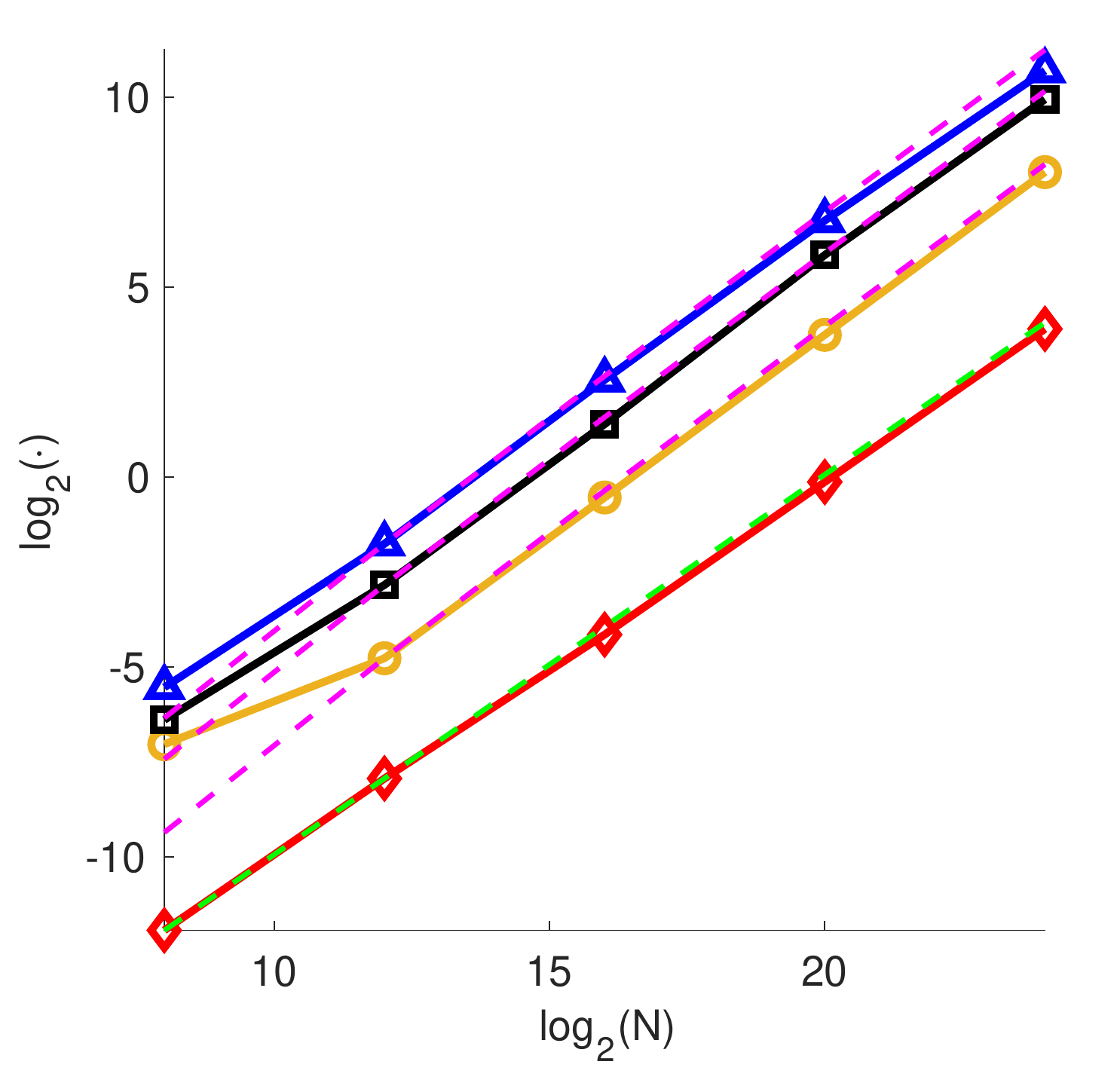} &
        \includegraphics[height=1.6in,trim={3cm 6.3cm 3.5cm 7cm}, clip]{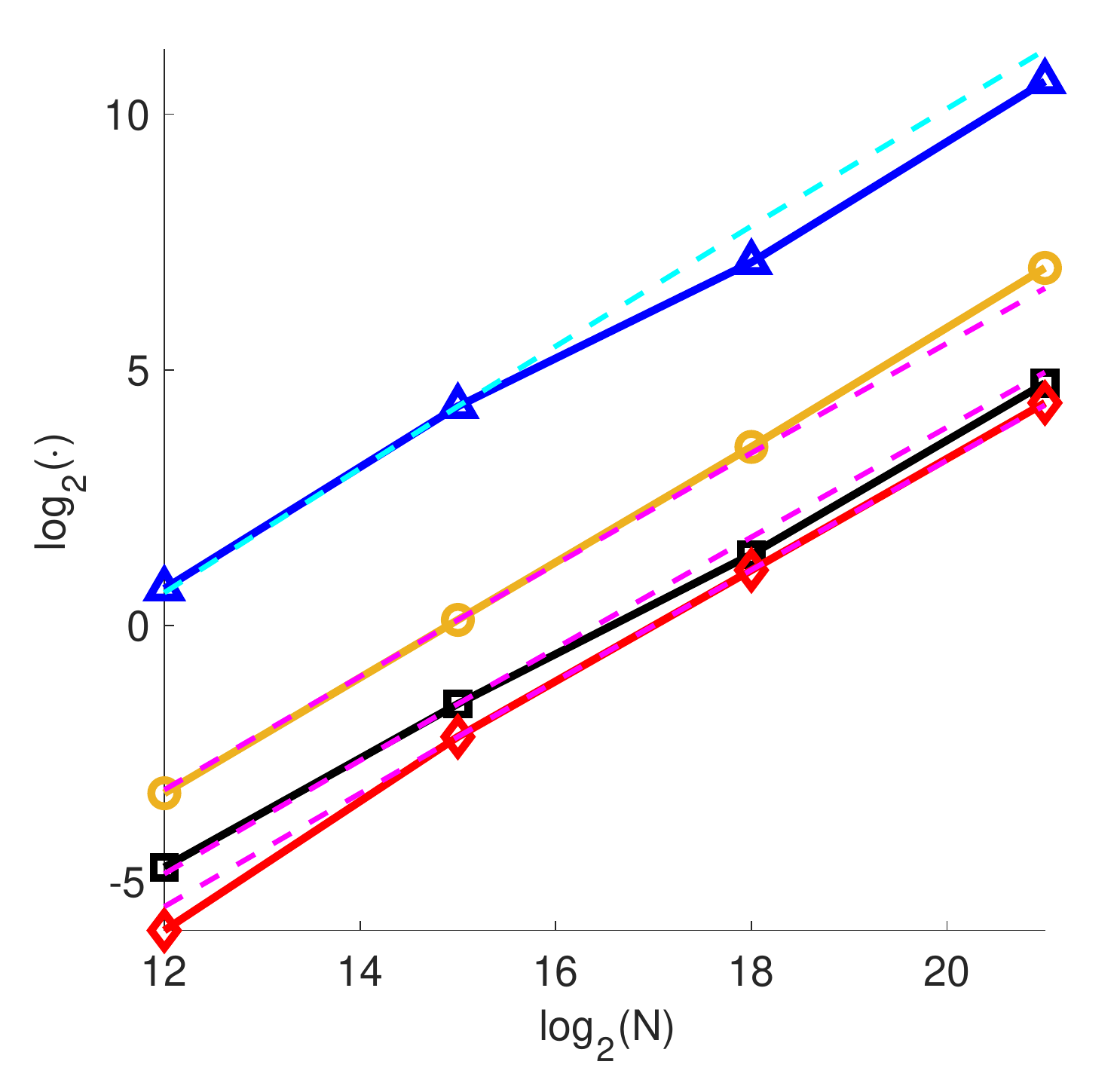} &
        \includegraphics[height=1.6in,trim={3cm 6.3cm 3.5cm 7cm}, clip]{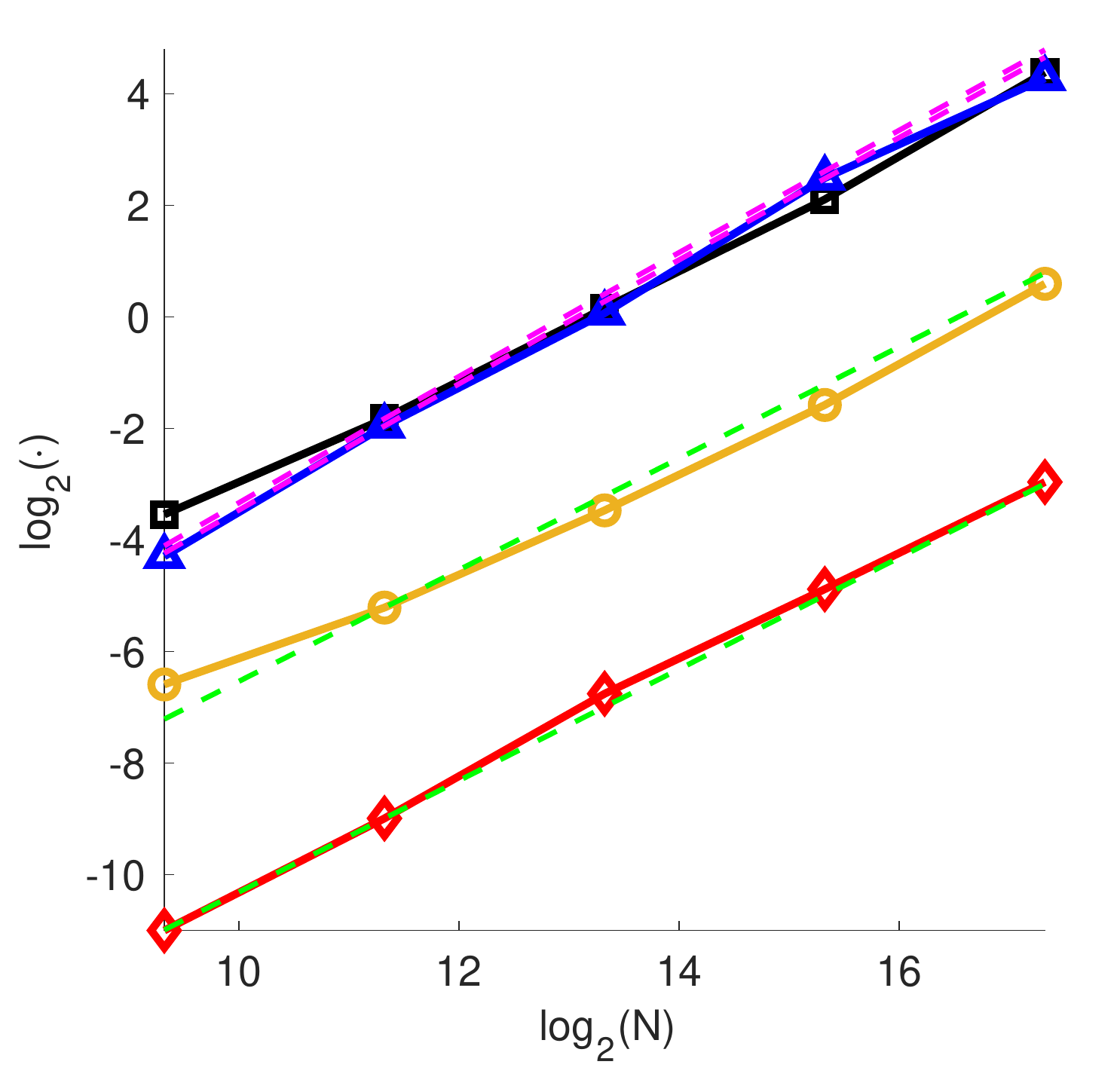} \\
         & (a) & (b) & (c)
    \end{tabular}
    \caption{The visualization of the computational complexity. $N$ is the size of the matrix. (a) the 2D uniform FIO given in \eqref{eqn:2D}. (b) the 3D Fourier transform given in \eqref{eqn:3D}. (c) the example in \eqref{eqn:2D_Hel}.}
    \label{fig:3D}
\end{figure}

Table~\ref{tab:2D} summarizes the results of this example for different grid sizes $N = n^2$ and different rank parameters $r$, $k$. It shows that the accuracy of the low-rank matrix recovery and the MIDBF stay almost of the same order, though the accuracy becomes slightly worse as the problem size increases. The slightly increasing error is due to the randomness of the proposed algorithm. As the problem size increases, the probability of capturing the low-rank matrix with a fixed rank parameter becomes smaller. Otherwise, when the rank parameter $r$ or $k$ increases, the accuracy of results will increase as well. In Figure~\ref{fig:3D} (a), we see that the time for computing recovery path matrix, the reconstruction time of the phase functions, the factorization time and the application time of the MIDBF scale nearly linearly, e.g. when $r=20$ and $k=30$. 

\begin{table}[!ht]
    \centering
    \begin{tabular}{cccccccccc}
        \toprule
        $n, r, k$ & $\epsilon^b$ & $\epsilon^K$ & $T_{path}$ & $T_{rec}$ & $T_{fac}$ & $T_{app}$ & $T_d/T_{app}$\\
        \midrule
        16, 10, 30 & 2.85e-07 & 1.27e-08 & 7.77e-03 & 8.10e-03 & 2.23e-02 & 2.81e-04 & 2.11e+01 \\
        16, 20, 20 & 5.08e-06 & 2.64e-09 & 9.75e-03 & 1.67e-02 & 2.14e-02 & 2.87e-04 & 2.79e+01 \\
        16, 20, 30 & 3.01e-07 & 2.63e-09 & 7.62e-03 & 1.19e-02 & 2.15e-02 & 2.56e-04 & 2.17e+01 \\
        \midrule
        64, 10, 30 & 4.93e-08 & 1.29e-08 & 4.29e-02 & 9.56e-02 & 3.29e-01 & 4.43e-03 & 2.37e+02 \\
        64, 20, 20 & 2.36e-06 & 2.42e-09 & 4.15e-02 & 1.70e-01 & 2.59e-01 & 3.42e-03 & 3.24e+02 \\
        64, 20, 30 & 3.51e-08 & 2.36e-09 & 3.66e-02 & 1.39e-01 & 2.96e-01 & 4.08e-03 & 2.23e+02 \\
        \midrule
        256, 10, 30 & 1.19e-08 & 1.34e-08 & 5.14e-01 & 1.27e+00 & 5.10e+00 & 4.00e-02 & 5.33e+03 \\
        256, 20, 20 & 2.28e-08 & 2.28e-09 & 6.52e-01 & 2.43e+00 & 4.37e+00 & 4.33e-02 & 5.55e+03 \\
        256, 20, 30 & 4.12e-09 & 2.23e-09 & 6.87e-01 & 2.62e+00 & 5.88e+00 & 5.63e-02 & 4.75e+03 \\
        \midrule
        1024, 10, 30 & 1.60e-08 & 1.41e-08 & 1.10e+01 & 2.72e+01 & 8.74e+01 & 6.86e-01 & 1.00e+05 \\
        1024, 20, 20 & 3.29e-09 & 2.29e-09 & 1.42e+01 & 6.01e+01 & 9.29e+01 & 1.09e+00 & 9.21e+04 \\
        1024, 20, 30 & 2.76e-09 & 2.33e-09 & 1.34e+01 & 5.74e+01 & 1.08e+02 & 9.13e-01 & 9.21e+04 \\
        \midrule
        4096, 10, 30 & 1.27e-08 & 1.47e-08 & 2.74e+02 & 6.25e+02 & 1.79e+03 & 1.64e+01 & 2.11e+06 \\
        4096, 20, 20 & 3.16e-09 & 2.23e-09 & 2.62e+02 & 1.02e+03 & 1.39e+03 & 1.47e+01 & 2.24e+06 \\
        4096, 20, 30 & 2.30e-09 & 2.15e-09 & 2.60e+02 & 9.82e+02 & 1.66e+03 & 1.49e+01 & 2.18e+06 \\
        \bottomrule
    \end{tabular}
    \caption{Numerical results for the 2D uniform FIO given in \eqref{eqn:2D}. $r$ is the rank parameter of the low-rank approximation of the phase function. $k$ is the rank parameter of the MIDBF.  $T_{path}$ is the time for computing the recovery path matrix. $T_{rec}$ is the time for recovering the phase functions, $T_{fac}$ is the time for computing the MIDBF, $T_{app}$ is the time for applying the MIDBF, and $T_d$ is the time for a direct summation in \eqref{eqn:2D}.}
    \label{tab:2D}
\end{table}

\paragraph{Example 2.}

In this example, we evaluate a 3D non-uniform Fourier transform:
\begin{equation}
    \label{eqn:3D}
    g(x) = \sum_{\xi \in \Omega} e^{2\pi i x^T \xi} \widehat{f}(\xi),
\end{equation}
where $X$ and $\Omega$ are the sets of $N$ points randomly selected in $[0,1)^3$.

{Table~\ref{tab:disc} shows the relationship between the discontinuity threshold $\tau$ and the number of detected discontinuous points. We set $\tau \le \frac{1}{4}$ in order to guarantee that the intersection of each recovered row and column share the same value. The results show that the numbers of detected discontinuity for rows and columns (denoted as $N_{\mathcal{D}_r}$ and $N_{\mathcal{D}_c}$, respectively) are both bounded in $O(1)$ when the problem size $N$ increases. Therefore, $\tau = \frac{1}{4}$ is an appropriate choice for this example.}

\begin{table}[!ht]
    \centering
    \begin{tabular}{cccccccccccccccc}
        \toprule
        $n$ & $\tau$ & $N_{\mathcal{D}_r}$ & $N_{\mathcal{D}_c}$ & & & $n$ & $\tau$ & $N_{\mathcal{D}_r}$ & $N_{\mathcal{D}_c}$ & & & $n$ & $\tau$ & $N_{\mathcal{D}_r}$ & $N_{\mathcal{D}_c}$ \\
        \cmidrule{1-4} \cmidrule{7-10} \cmidrule{13-16}
        8 & $\frac{1}{4}$ & 0 & 0 & & & 16 & $\frac{1}{4}$ & 0 & 0 & & & 32 & $\frac{1}{4}$ & 0 & 0 \\
        8 & $\frac{1}{6}$ & 3.0 & 2.9 & & & 16 & $\frac{1}{6}$ & 0 & 0 & & & 32 & $\frac{1}{6}$ & 0 & 0 \\
        8 & $\frac{1}{8}$ & 29.1 & 32.0 & & & 16 & $\frac{1}{8}$ & 0.2 & 0.1 & & & 32 & $\frac{1}{8}$ & 0 & 0 \\
        8 & $\frac{1}{10}$ & 82.5 & 82.7 & & & 16 & $\frac{1}{10}$ & 2.4 & 1.6 & & & 32 & $\frac{1}{10}$ & 0 & 0 \\
        \bottomrule
    \end{tabular}
    \caption{{The number of discontinuous points of the 3D non-uniform Fourier transform given in \eqref{eqn:3D}. $N = n^3$ is the size of grid. $\tau$ is the threshold for detecting the discontinuity. $N_{\mathcal{D}_r}$ and $N_{\mathcal{D}_c}$ are the numbers of discontinuous points along recovery rows and columns of the phase matrix, respectively.}}
    \label{tab:disc}
\end{table}

Table~\ref{tab:3D} summarizes the results of this example for different grid sizes $N = n^3$ and different rank parameters $r$ in the low-rank approximation of the phase function. In the MIDBF, the rank parameter $k$ is $80$. The accuracy of the low-rank matrix recovery and the MIDBF stay almost of the same order in Table~\ref{tab:3D}. In Figure~\ref{fig:3D} (b), we see that each part of the whole process scales nearly linearly, e.g., when $r=5$.

\begin{table}[!ht]
    \centering
    \begin{tabular}{cccccccccc}
        \toprule
        $n, r$ & $\epsilon^b$ & $\epsilon^K$ & $T_{path}$ & $T_{rec}$ & $T_{fac}$ & $T_{app}$ & $T_d/T_{app}$\\
        \midrule
        16, 3 & 2.61e-01 & 3.22e-01 & 1.11e-01 & 2.76e-02 & 1.32e+00 & 1.24e-02 & 7.33e+01 \\
        16, 4 & 1.10e-06 & 1.02e-14 & 1.21e-01 & 3.80e-02 & 1.81e+00 & 1.72e-02 & 5.15e+01 \\
        16, 5 & 1.10e-06 & 6.67e-15 & 1.03e-01 & 3.77e-02 & 1.65e+00 & 1.60e-02 & 5.35e+01 \\
        \midrule
        32, 3 & 2.85e-01 & 3.80e-01 & 1.15e+00 & 2.26e-01 & 1.32e+01 & 1.09e-01 & 6.01e+02 \\
        32, 4 & 4.19e-08 & 8.66e-15 & 1.15e+00 & 2.91e-01 & 1.93e+01 & 2.64e-01 & 2.41e+02 \\
        32, 5 & 3.85e-08 & 1.15e-14 & 1.08e+00 & 3.46e-01 & 1.95e+01 & 2.21e-01 & 2.74e+02 \\
        \midrule
        64, 3 & 3.37e-01 & 4.56e-01 & 1.16e+01 & 1.80e+00 & 9.72e+01 & 1.01e+00 & 4.25e+03 \\
        64, 4 & 5.33e-08 & 2.80e-14 & 1.09e+01 & 2.28e+00 & 1.37e+02 & 2.12e+00 & 1.71e+03 \\
        64, 5 & 4.91e-08 & 2.04e-14 & 1.13e+01 & 2.62e+00 & 1.38e+02 & 2.12e+00 & 1.87e+03 \\
        \midrule
        128, 3 & 4.54e-01 & 5.36e-01 & 1.32e+02 & 1.86e+01 & 8.60e+02 & 8.51e+00 & 3.50e+04 \\
        128, 4 & 2.92e-09 & 4.67e-14 & 1.27e+02 & 2.32e+01 & 1.59e+03 & 2.14e+01 & 1.51e+04 \\
        128, 5 & 3.42e-09 & 4.63e-14 & 1.27e+02 & 2.67e+01 & 1.60e+03 & 2.05e+01 & 1.56e+04 \\
        \bottomrule
    \end{tabular}
    \caption{Numerical results for the 3D Fourier transform given in \eqref{eqn:3D}. $T_d$ is the time for a direct summation in \eqref{eqn:3D}.}
    \label{tab:3D}
\end{table}

\paragraph{Example 3.}

The final example is the oscillatory part of the Green's function of a Helmholtz equation \cite{green}:
\begin{equation}
    \label{eqn:2D_Hel}
    g(x) = \sum_{\xi \in \Omega} e^{2\pi i \Phi(x,\xi)} \widehat{f}(\xi), \quad x \in X,
\end{equation}
where $\Phi(x,\xi) = h\cdot \norm{x-\xi}_2$ and $h=\frac{\sqrt{N}}{10} \sim \O{n}$. $X$ and $\Omega$ are the sets of $N$ points generated via a triangular mesh to discretize the surface of a unit sphere. The triangular mesh is generated by uniformly refining an icosahedron and projecting the new mesh nodes, which are the old mesh edge center, onto the sphere. The submatrix of the oscillatory part of the Green's function corresponding to one half of the sphere in X and the other half of the sphere in $\Omega$ is chosen as the matrix to be reconstructed, factorized, and applied to a random vector.

In this example, rank parameters $r = 50$ and $k = 50$. As shown in Table~\ref{tab:2D_Hel}, the accuracy of the low-rank matrix recovery and the MIDBF stay almost of the same order. The result in Figure~\ref{fig:3D} (c) demonstrates the efficiency of the proposed framework.

\begin{table}[!ht]
    \centering
    \begin{tabular}{cccccccccc}
        \toprule
        $n$ & $\epsilon^b$ & $\epsilon^K$ & $T_{path}$ & $T_{rec}$ & $T_{fac}$ & $T_{app}$ & $T_d/T_{app}$\\
        \midrule
        640 & 3.18e-09 & 1.31e-09 & 1.04e-02 & 8.58e-02 & 5.10e-02 & 4.87e-04 & 1.79e+02 \\
        2560 & 8.30e-09 & 4.48e-09 & 2.70e-02 & 2.82e-01 & 2.58e-01 & 1.96e-03 & 3.51e+02 \\
        10240 & 2.79e-08 & 1.12e-08 & 9.03e-02 & 1.12e+00 & 1.05e+00 & 9.24e-03 & 9.58e+02 \\
        40960 & 2.33e-08 & 2.43e-08 & 3.35e-01 & 4.31e+00 & 5.61e+00 & 3.39e-02 & 3.31e+03 \\
        163840 & 5.38e-08 & 5.98e-08 & 1.51e+00 & 2.09e+01 & 1.94e+01 & 1.28e-01 & 1.35e+04 \\
        \bottomrule
    \end{tabular}
    \caption{Numerical results for the case given in \eqref{eqn:2D_Hel}. $T_d$ is the time for a direct summation in \eqref{eqn:2D_Hel}.}
    \label{tab:2D_Hel}
\end{table}

\section{Conclusion} \label{sec:conclusion}

This paper introduced a framework for $\O{N \L{N}}$ evaluation of the multidimensional oscillatory integral transform $g(x) = \int e^{2\pi\i \Phi(x,\xi)} f(\xi)d\xi$. In the case of indirect access of the phase functions, this paper proposed a novel fast algorithm for recovering the phase functions in $\O{N \L{N}}$ operations. Second, a new BF, the multidimensional interpolative decomposition butterfly factorization (MIDBF), for multidimensional kernel matrices in the form of a low-rank factorization is proposed, and it requires only $\O{N \L{N}}$ operations to evaluate the oscillatory integral transform.

\vspace{0.5cm}
{\bf Acknowledgments.} Z. C. was partially supported by the Ministry of Education in Singapore under the grant MOE2018-T2-2-147. H. Y. was partially supported by NSF under the grant award 1945029.

\section{Appendix}
\subsection{Proof of Lemma~\ref{1LRthm}}

\begin{proof}
    First, let one of the block matrices be $\phi$, which is partitioned by discontinuous point sets (corresponding to Line~\ref{alg:rm1:pib} in Algorithm~\ref{alg:rm1}). Then, Line~\ref{alg:rm1:rv11}-\ref{alg:rm1:rv12} in Algorithm~\ref{alg:rm1} can obtain the unique recovery values of the first $3 \times 3$ entries of $\phi$, which are the first three entries in the first three columns. 
    
    Next step, consider the intersection of the fourth row and the fourth column in $\phi$. On one hand, after applying Algorithm~\ref{alg:rv1} in the first column, $\phi(4,1)$ will be obtained by
    \begin{equation}
        \phi(4,1) = \phi(1,1) - 3\phi(2,1) + 3\phi(3,1) + \epsilon_1,
    \end{equation}
    where $\epsilon_1 \in (-\frac{1}{16}, \frac{1}{16})$, according to the property of the first column of $\phi$. Since $\bmod(\phi(4,1),1)$ has been given, the recovery value of $\phi(4,1)$ will be unique.
    
    Similarly, $\phi(4,2)$ and $\phi(4,3)$ can be evaluated by
    \begin{equation}
        \begin{split}
            \phi(4,2) & = \phi(1,2) - 3\phi(2,2) + 3\phi(3,2) + \epsilon_2, \\
            \phi(4,3) & = \phi(1,3) - 3\phi(2,3) + 3\phi(3,3) + \epsilon_3,
        \end{split}
    \end{equation}
    through the second and the third column, where $\epsilon_2, \epsilon_3 \in (-\frac{1}{16}, \frac{1}{16})$.
    
    Next, apply Algorithm~\ref{alg:rv1} to the fourth row to evaluate $\phi(4,4)$:
    \begin{equation}
        \begin{split}
            \phi(4,4) & = \phi(4,1) - 3\phi(4,2) + 3\phi(4,3) + \epsilon_4 \\
            & = \phi(1,1) - 3\phi(2,1) + 3\phi(3,1) + \epsilon_1 - 3\phi(1,2) + 9\phi(2,2) - 9\phi(3,2) - 3\epsilon_2 \\
            & \quad + 3\phi(1,3) - 9\phi(2,3) + 9\phi(3,3) + 3\epsilon_3 + \epsilon_4 \\
            & = C + \epsilon_1 - 3\epsilon_2 + 3\epsilon_3 + \epsilon_4,
        \end{split}
    \end{equation}
    where $\epsilon_4 \in (-\frac{1}{16}, \frac{1}{16})$ and $C = \phi(1,1) - 3\phi(2,1) + 3\phi(3,1) - 3\phi(1,2) + 9\phi(2,2) - 9\phi(3,2) + 3\phi(1,3) - 9\phi(2,3) + 9\phi(3,3)$.
    
    Since $\epsilon_1 - 3\epsilon_2 + 3\epsilon_3 + \epsilon_4 \in (-\frac{1}{2}, \frac{1}{2})$, $\phi(4,4)$ can be obtained by identifying a unique integer $a$, such that
    \begin{equation}
        \label{int_a}
        \bmod(\phi(4,4),1) + a \in (C-\frac{1}{2}, C+\frac{1}{2}).
    \end{equation}
    Then, the recovery value of $\phi(4,4)$ through the fourth row will be unique as $\bmod(\phi(4,4),1) + a$.
    
    On the other hand, the same method can be applied to obtain
    \begin{equation}
        \begin{split}
            \phi(1,4) & = \phi(1,1) - 3\phi(1,2) + 3\phi(1,3) + \epsilon_1^\prime, \\
            \phi(2,4) & = \phi(2,1) - 3\phi(2,2) + 3\phi(2,3) + \epsilon_2^\prime, \\
            \phi(3,4) & = \phi(3,1) - 3\phi(3,2) + 3\phi(3,3) + \epsilon_3^\prime,
        \end{split}
    \end{equation}
    where $\epsilon_1^\prime, \epsilon_2^\prime, \epsilon_3^\prime \in (-\frac{1}{16}, \frac{1}{16})$.
    
    Next, apply Algorithm~\ref{alg:rv1} again to the fourth column to evaluate $\phi(4,4)$ accompanying with a parameter $\epsilon_4^\prime \in (-\frac{1}{16}, \frac{1}{16})$:
    \begin{equation}
        \begin{split}
            \phi(4,4) & = \phi(1,4) - 3\phi(2,4) + 3\phi(3,4) + \epsilon_4^\prime \\
            & = \phi(1,1) - 3\phi(1,2) + 3\phi(1,3) + \epsilon_1^\prime - 3\phi(2,1) + 9\phi(2,2) - 9\phi(2,3) - 3\epsilon_2^\prime \\
            & \quad + 3\phi(3,1) - 9\phi(3,2) + 9\phi(3,3) + 3\epsilon_3^\prime + \epsilon_4^\prime \\
            & = C + \epsilon_1^\prime - 3\epsilon_2^\prime + 3\epsilon_3^\prime + \epsilon_4^\prime \\
            & \in (C-\frac{1}{2}, C+\frac{1}{2}).
        \end{split}
    \end{equation}
    
    Similarly, $\phi(4,4)$ can be obtained by identifying a unique integer $b$, such that 
    \begin{equation}
        \label{int_b}
        \bmod(\phi(4,4),1) + b \in (C-\frac{1}{2}, C+\frac{1}{2}).
    \end{equation}
    
    Combining (\ref{int_a}) and (\ref{int_b}), integers $a,b \in (C-\bmod(\phi(4,4),1)-\frac{1}{2}, C-\bmod(\phi(4,4),1)+\frac{1}{2})$, which is obvious to conclude that $a = b$. Thus, the intersection $\phi(4,4)$ recovered by the fourth row and the fourth column using Algorithm~\ref{alg:rv1} will share the same value.
    
    The same, when the recovered values of the first three entries of the second to fourth columns have been obtained using the previous method, a unique recovery value of $\phi(4,5)$ would be evaluated, which means that the intersection recovered by the fourth row and the fifth column will share the same value.
    
    Furthermore, the unique recovery values of $\phi(4,6), \phi(4,7), \dots, \phi(4,m)$ can also be evaluated. Therefore, when the values of the first three entries of the first three columns have been fixed, any entry in the fourth row as the intersection will share the same value when recovering the corresponding row and column. The method can be applied to prove the same property in the rest rows.
    
    In conclusion, if the nine values of the first three entries of the first three columns have been fixed, any recovered row and column by Algorithm~\ref{alg:rv1} will share the same value at the intersection.
\end{proof}


\bibliographystyle{abbrv}
\bibliography{ref}

\begin{thebibliography}{10}

\bibitem{Gang}
G.~Bao and W.~W. Symes.
\newblock {Computation of Pseudo-Differential Operators}.
\newblock {\em SIAM Journal on Scientific Computing}, 17(2):416--429, 1996.

\bibitem{Mock2}
J.~P. Boyd and F.~Xu.
\newblock {Divergence ({R}unge {P}henomenon) for least-squares polynomial
  approximation on an equispaced grid and {M}ock {C}hebyshev subset
  interpolation}.
\newblock {\em Applied Mathematics and Computation}, 210(1):158 -- 168, 2009.

\bibitem{James:2017}
J.~Bremer.
\newblock {An algorithm for the rapid numerical evaluation of {B}essel
  functions of real orders and arguments}.
\newblock {\em arXiv:1705.07820 [math.NA]}, 2017.

\bibitem{Bremer201815}
J.~Bremer.
\newblock {An algorithm for the numerical evaluation of the associated
  {L}egendre functions that runs in time independent of degree and order}.
\newblock {\em Journal of Computational Physics}, 360:15 -- 38, 2018.

\bibitem{dt1}
K.~{Buchin} and W.~{Mulzer}.
\newblock {Delaunay Triangulations in $O(\text{sort}(n))$ Time and More}.
\newblock In {\em 2009 50th Annual IEEE Symposium on Foundations of Computer
  Science}, pages 139--148, Oct 2009.

\bibitem{FIO09}
E.~J. Cand{\`e}s, L.~Demanet, and L.~Ying.
\newblock {A Fast Butterfly Algorithm for the Computation of {{F}ourier}
  Integral Operators.}
\newblock {\em Multiscale Modeling and Simulation}, 7(4):1727--1750, 2009.

\bibitem{Unw2}
M.~Costantin, A.~Farina, and F.~Zirilli.
\newblock {A Fast Phase Unwrapping Algorithm for SAR Interferometry}.
\newblock {\em IEEE TRANSACTIONS ON GEOSCIENCE AND REMOTE SENSING}, 37(1),
  1999.

\bibitem{green}
B.~Davies.
\newblock {\em {Green's Functions}}, pages 163--179.
\newblock Springer New York, New York, NY, 2002.

\bibitem{DTMST}
M.~de~Berg, O.~Cheong, M.~van Kreveld, and M.~Overmars.
\newblock {\em {Delaunay Triangulations}}, pages 191--218.
\newblock Springer Berlin Heidelberg, Berlin, Heidelberg, 2008.

\bibitem{Yingwave}
L.~Demanet and L.~Ying.
\newblock {Fast wave computation via {{F}ourier} integral operators.}
\newblock {\em Math. Comput.}, 81(279), 2012.

\bibitem{RadiusSearch}
M.~T. Dickerson and R.~Drysdale.
\newblock {Fixed-radius near neighbors search algorithms for points and
  segments}.
\newblock {\em Information Processing Letters}, 35(5):269--273, 1990.

\bibitem{Butterfly4}
B.~Engquist and L.~Ying.
\newblock {A fast directional algorithm for high frequency acoustic scattering
  in two dimensions}.
\newblock {\em Communications in Mathematical Sciences}, 7(2):327--345, 06
  2009.

\bibitem{NUFFT}
L.~Greengard and J.-Y. Lee.
\newblock {Accelerating the {N}onuniform {F}ast {F}ourier {T}ransform}.
\newblock {\em SIAM Review}, 46(3):443--454, 2004.

\bibitem{LUBF}
H.~Guo, Y.~Liu, J.~Hu, and E.~Michielssen.
\newblock {A Butterfly-Based Direct Integral-Equation Solver Using Hierarchical
  {LU} Factorization for Analyzing Scattering From Electrically Large
  Conducting Objects}.
\newblock {\em IEEE Transactions on Antennas and Propagation},
  65(9):4742--4750, Sept 2017.

\bibitem{randomSVD}
N.~Halko, P.-G. Martinsson, and J.~A. Tropp.
\newblock {Finding structure with randomness: Probabilistic algorithms for
  constructing approximate matrix decompositions}.
\newblock {\em SIAM review}, 53(2):217--288, 2011.

\bibitem{Mock1}
P.~Hoffman and K.~Reddy.
\newblock {Numerical Differentiation by High Order Interpolation}.
\newblock {\em SIAM Journal on Scientific and Statistical Computing},
  8(6):979--987, 1987.

\bibitem{precon3}
H.~Isozaki and J.~L. Rousseau.
\newblock {Pseudodifferential Multi-Product Representation of the Solution
  Operator of a Parabolic Equation}.
\newblock {\em Communications in Partial Differential Equations},
  34(7):625--655, 2009.

\bibitem{JIANCHUN1995162}
L.~Jianchun, G.~A. Pope, and K.~Sepehrnoori.
\newblock {A high-resolution finite-difference scheme for nonuniform grids}.
\newblock {\em Applied Mathematical Modelling}, 19(3):162 -- 172, 1995.

\bibitem{BFS}
C.~Y. {Lee}.
\newblock {An Algorithm for Path Connections and Its Applications}.
\newblock {\em IRE Transactions on Electronic Computers}, EC-10(3):346--365,
  Sep. 1961.

\bibitem{IBF}
Y.~Li and H.~Yang.
\newblock {Interpolative Butterfly Factorization}.
\newblock {\em SIAM Journal on Scientific Computing}, 39(2):A503--A531, 2017.

\bibitem{BF}
Y.~Li, H.~Yang, E.~R. Martin, K.~L. Ho, and L.~Ying.
\newblock {Butterfly {F}actorization}.
\newblock {\em Multiscale Modeling \& Simulation}, 13(2):714--732, 2015.

\bibitem{MBF}
Y.~Li, H.~Yang, and L.~Ying.
\newblock {Multidimensional butterfly factorization}.
\newblock {\em Applied and Computational Harmonic Analysis}, 2017.

\bibitem{HSSBF}
Y.~Liu, H.~Guo, and E.~Michielssen.
\newblock {An {HSS} Matrix-Inspired Butterfly-Based Direct Solver for Analyzing
  Scattering From Two-Dimensional Objects}.
\newblock {\em IEEE Antennas and Wireless Propagation Letters}, 16:1179--1183,
  2017.

\bibitem{dt2}
S.~Lo.
\newblock {Parallel Delaunay triangulation in three dimensions}.
\newblock {\em Computer Methods in Applied Mechanics and Engineering},
  237-240:88 -- 106, 2012.

\bibitem{Butterfly1}
E.~Michielssen and A.~Boag.
\newblock {A multilevel matrix decomposition algorithm for analyzing scattering
  from large structures}.
\newblock {\em Antennas and Propagation, IEEE Transactions on},
  44(8):1086--1093, Aug 1996.

\bibitem{Unw3}
G.~Nico, G.~Palubinskas, and M.~Datcu.
\newblock {Bayesian Approaches to Phase Unwrapping: Theoretical Study}.
\newblock {\em IEEE TRANSACTIONS ON SIGNAL PROCESSING}, 48(9), 2000.

\bibitem{Butterfly2}
M.~O'Neil, F.~Woolfe, and V.~Rokhlin.
\newblock {An algorithm for the rapid evaluation of special function
  transforms}.
\newblock {\em Appl. Comput. Harmon. Anal.}, 28(2):203--226, 2010.

\bibitem{IDBF}
Q.~Pang, K.~L. Ho, and H.~Yang.
\newblock {Interpolative Decomposition Butterfly Factorization}.
\newblock {\em arXiv:1809.10573 [math.NA]}, 2018.

\bibitem{Prim}
R.~Prim.
\newblock {Shortest Connection Networks And Some Generalizations}.
\newblock {\em Bell System Technical Journal}, 36:1389--1401, 11 1957.

\bibitem{precon1}
J.~L. Rousseau.
\newblock {Fourier-Integral-Operator Approximation of Solutions to First-Order
  Hyperbolic Pseudodifferential Equations {I}: Convergence in Sobolev Spaces}.
\newblock {\em Communications in Partial Differential Equations},
  31(6):867--906, 2006.

\bibitem{precon2}
J.~L. Rousseau and G.~H{\"o}rmann.
\newblock {Fourier-integral-operator approximation of solutions to first-order
  hyperbolic pseudodifferential equations {II}: Microlocal analysis}.
\newblock {\em Journal de Mathématiques Pures et Appliqu{\'e}es}, 86(5):403 --
  426, 2006.

\bibitem{Alex2}
D.~Ruiz-Antolín and A.~Townsend.
\newblock {A Nonuniform Fast Fourier Transform Based on Low Rank
  Approximation}.
\newblock {\em SIAM Journal on Scientific Computing}, 40(1):A529--A547, 2018.

\bibitem{dt3}
M.~Smid.
\newblock {The Well-Separated Pair Decomposition and Its Applications}.
\newblock In {\em Handbook of Approximation Algorithms and Metaheuristics},
  2007.

\bibitem{Unw1}
E.~Trouv{\'e}, J.-M. Nicolas, and H.~Ma{\^i}tre.
\newblock {Improving Phase Unwrapping Techniques by the Use of Local Frequency
  Estimates}.
\newblock {\em IEEE TRANSACTIONS ON GEOSCIENCE AND REMOTE SENSING}, 36(6),
  1998.

\bibitem{KNN}
P.~Vaidya.
\newblock {An $O(n \log n)$ Algorithm for the All-Nearest-Neighbors Problem.}
\newblock {\em Discrete and computational geometry}, 4(2):101--116, 1989.

\bibitem{doi:10.1137/1.9781611970999}
C.~Van~Loan.
\newblock {\em {Computational Frameworks for the Fast Fourier Transform}}.
\newblock Society for Industrial and Applied Mathematics, 1992.

\bibitem{VASILYEV2000746}
O.~V. Vasilyev.
\newblock {High Order Finite Difference Schemes on Non-uniform Meshes with Good
  Conservation Properties}.
\newblock {\em Journal of Computational Physics}, 157(2):746 -- 761, 2000.

\bibitem{NUFFTorBF}
H.~{Yang}.
\newblock {A unified framework for oscillatory integral transforms: When to use
  NUFFT or butterfly factorization?}
\newblock {\em Journal of Computational Physics}, 388:103--122, Jul 2019.

\end{thebibliography}

\end{document}